\PassOptionsToPackage{unicode}{hyperref}
\PassOptionsToPackage{hyphens}{url}
\PassOptionsToPackage{dvipsnames,svgnames,x11names}{xcolor}
\documentclass[
  12pt]{article}

\usepackage{amsmath,amssymb}
\usepackage{iftex}
\ifPDFTeX
  \usepackage[T1]{fontenc}
  \usepackage[utf8]{inputenc}
  \usepackage{textcomp} 
\else 
  \usepackage{unicode-math}
  \defaultfontfeatures{Scale=MatchLowercase}
  \defaultfontfeatures[\rmfamily]{Ligatures=TeX,Scale=1}
\fi
\usepackage{lmodern}
\ifPDFTeX\else  
\fi
\IfFileExists{upquote.sty}{\usepackage{upquote}}{}
\IfFileExists{microtype.sty}{
  \usepackage[]{microtype}
  \UseMicrotypeSet[protrusion]{basicmath} 
}{}
\makeatletter
\@ifundefined{KOMAClassName}{
  \IfFileExists{parskip.sty}{%
    \usepackage{parskip}
  }{
    \setlength{\parindent}{0pt}
    \setlength{\parskip}{6pt plus 2pt minus 1pt}}
}{
  \KOMAoptions{parskip=half}}
\makeatother
\usepackage{xcolor}
\makeatletter
\ifx\paragraph\undefined\else
  \let\oldparagraph\paragraph
  \renewcommand{\paragraph}{
    \@ifstar
      \xxxParagraphStar
      \xxxParagraphNoStar
  }
  \newcommand{\xxxParagraphStar}[1]{\oldparagraph*{#1}\mbox{}}
  \newcommand{\xxxParagraphNoStar}[1]{\oldparagraph{#1}\mbox{}}
\fi
\ifx\subparagraph\undefined\else
  \let\oldsubparagraph\subparagraph
  \renewcommand{\subparagraph}{
    \@ifstar
      \xxxSubParagraphStar
      \xxxSubParagraphNoStar
  }
  \newcommand{\xxxSubParagraphStar}[1]{\oldsubparagraph*{#1}\mbox{}}
  \newcommand{\xxxSubParagraphNoStar}[1]{\oldsubparagraph{#1}\mbox{}}
\fi
\makeatother

\usepackage{longtable,booktabs,array}
\usepackage{calc} 
\usepackage{etoolbox}
\makeatletter
\patchcmd\longtable{\par}{\if@noskipsec\mbox{}\fi\par}{}{}
\makeatother
\IfFileExists{footnotehyper.sty}{\usepackage{footnotehyper}}{\usepackage{footnote}}
\makesavenoteenv{longtable}
\usepackage{graphicx}
\makeatletter
\def\maxwidth{\ifdim\Gin@nat@width>\linewidth\linewidth\else\Gin@nat@width\fi}
\def\maxheight{\ifdim\Gin@nat@height>\textheight\textheight\else\Gin@nat@height\fi}
\makeatother
\setkeys{Gin}{width=\maxwidth,height=\maxheight,keepaspectratio}
\makeatletter
\def\fps@figure{htbp}
\makeatother

\makeatletter
\@ifpackageloaded{caption}{}{\usepackage{caption}}
\AtBeginDocument{%
\ifdefined\contentsname
  \renewcommand*\contentsname{Table of contents}
\else
  \newcommand\contentsname{Table of contents}
\fi
\ifdefined\listfigurename
  \renewcommand*\listfigurename{List of Figures}
\else
  \newcommand\listfigurename{List of Figures}
\fi
\ifdefined\listtablename
  \renewcommand*\listtablename{List of Tables}
\else
  \newcommand\listtablename{List of Tables}
\fi
\ifdefined\figurename
  \renewcommand*\figurename{Figure}
\else
  \newcommand\figurename{Figure}
\fi
\ifdefined\tablename
  \renewcommand*\tablename{Table}
\else
  \newcommand\tablename{Table}
\fi
}
\@ifpackageloaded{float}{}{\usepackage{float}}
\floatstyle{ruled}
\@ifundefined{c@chapter}{\newfloat{codelisting}{h}{lop}}{\newfloat{codelisting}{h}{lop}[chapter]}
\floatname{codelisting}{Listing}

\makeatother
\makeatletter
\@ifpackageloaded{caption}{}{\usepackage{caption}}
\@ifpackageloaded{subcaption}{}{\usepackage{subcaption}}
\makeatother

\ifLuaTeX
  \usepackage{selnolig}  
\fi
\usepackage[]{natbib}
\usepackage{bookmark}

\IfFileExists{xurl.sty}{\usepackage{xurl}}{} 
\hypersetup{
  pdftitle={Title},
  pdfauthor={Author 1; Author 2},
  pdfkeywords={3 to 6 keywords, that do not appear in the title},
  colorlinks=true,
  linkcolor={blue},
  filecolor={Maroon},
  citecolor={Blue},
  urlcolor={Blue},
  pdfcreator={LaTeX via pandoc}}

\newcommand{\anon}{1}

\usepackage{appendix}
\usepackage{amsthm,mathrsfs}
\usepackage{tikz,bm,xr}
\usepackage{multicol,multirow}
\usepackage{algorithm,algpseudocode}
\usepackage{graphicx}
\usepackage{booktabs}

\newtheorem{theorem}{Theorem}[section]
\newtheorem{lemma}[theorem]{Lemma}

\theoremstyle{definition}
\newtheorem{definition}[theorem]{Definition}

\newcommand{\Abs}[2]{{\Big|#1\Big|}_{#2}}
\newcommand{\norm}[2]{{\|#1\|}_{#2}}
\newcommand{\Norm}[2]{{\Big\|#1\Big\|}_{#2}}

\newcommand{\E}{\normalfont \mathbb{E}}

\newcommand{\Prob}{\normalfont \mathbb{P}}

\newcommand{\bbR}{\mathbb{R}}

\newcommand{\frakm}{\mathfrak{m}}

\newcommand{\frakb}{\mathfrak{b}}

\newcommand{\tl}{\mathrm{TL}}
\newcommand{\mbf}[1]{\mathbf{#1}}

\newcommand{\XtT}{X_{t,T}}
\newcommand{\ZtT}{Z_{t,T}}
\newcommand{\RtT}{R_{t,T}}

\newcommand{\tT}{\frac{t}{T}}
\newcommand{\tTz}{\frac{t_0}{T_0}}

\begin{document}

\def\spacingset#1{\renewcommand{\baselinestretch}%
{#1}\small\normalsize} \spacingset{1}


\if1\anon
{
    \title{\bf Transfer Learning and Locally Linear Regression for Locally Stationary Time Series}
    \author{Jinwoo Park\\
    Department of Statistics, Seoul National University}
    \maketitle
} \fi

\if0\anon
{
  \bigskip
  \bigskip
  \bigskip
  \begin{center}
    {\LARGE\bf Title}
\end{center}
  \medskip
} \fi

\bigskip

\begin{abstract}
This paper investigates locally linear regression for locally stationary time series and develops theoretical results for locally linear smoothing and transfer learning. Existing analyses have focused on local constant estimators and given samples, leaving the principles of transferring knowledge from auxiliary sources across heterogeneous time-varying domains insufficiently established. 
We derive uniform convergence for multivariate locally linear estimators under strong mixing. The resulting error expansion decomposes stochastic variation, smoothing bias, and a term induced by local stationarity. This additional term, originating from the locally stationary structure, has smaller order than in the Nadaraya-Watson benchmark, explaining the improved local linear performance. 
Building on these results, we propose bias-corrected transfer learned estimators that connect a sparsely observed series with densely observed related sources through a smoothly varying bias function defined over rescaled time and covariates. An additional refinement shows how local temporal adjustment of this bias enhances stability and enables efficient information borrowing across domains. Simulation studies and an empirical analysis of international fuel prices support the theoretical predictions and demonstrate the practical advantages of transfer learning.
\end{abstract}

\noindent%
\vfill

\newpage
\section{Introduction}\label{sec-intro}
Classical time series analysis has traditionally assumed stationarity, meaning that the mean, variance, and dependence structure of a process remain fixed over time. This assumption enables the use of ergodic and mixing arguments and underlies the development of autoregressive, spectral, and state space methods. Yet empirical data seldom satisfy such stability. Economic and financial series evolve under shifting regimes, policy changes, and technological progress. Volatility levels move, and cross correlations respond to market conditions. When these forms of temporal change arise, stationary models become misspecified, and additional observations no longer provide information of the same nature. As a result, inference based on stationary approximations loses reliability.

To address the shortcomings of global stationarity, the notion of \emph{locally stationary processes} was introduced to represent nonstationarity as a smooth evolution of the data generating mechanism (\cite{DAHLHAUS1996139}; \cite{Dahlhaus.R.1997}). Allowing the dependence structure to vary continuously with rescaled time preserves local stability while remaining amenable to asymptotic analysis. This perspective has proved useful in many settings, including time varying ARCH and GARCH models (\cite{DahlhausRainer2006SIfT}), tests for second-order stationarity in multivariate processes (\cite{JENTSCH2015124}), nonlinear locally stationary processes (\cite{Dalhaus.Richter.Wu_2019}), high-dimensional covariance and spectral density estimation (\cite{Zhang.Wu.2021}), nonparametric regression for functional time series (\cite{Kurisu_2022}), and graphical modeling for nonstationary series (\cite{Basu.SubbaRao_2023}). More recent work includes Gaussian approximation with optimal rates (\cite{Bonnerjee.Karmakar.Wu.2024}) and time varying network estimation for high-dimensional data (\cite{CHEN2025105941}). Together, these contributions show that many empirical forms of temporal dependence can be represented as smoothly changing, providing a tractable framework that remains closely aligned with observed dynamics.

Within this setting, nonparametric regression offers a flexible way to capture evolving conditional relationships. The first rigorous asymptotic analysis was provided by \cite{Vogt_2012}, who studied estimation of the regression function \(m\) using the Nadaraya--Watson estimator. That work extended stationary kernel regression theories from \cite{MasryElias1996MLPR}, \cite{Bosq1998}, and \cite{Hansen_2008}, and established that the estimation error consists of the usual stochastic fluctuation and smoothing bias, together with an additional term generated by temporal drift. This remainder reflects how local stationarity blends dependence control via mixing with smooth time variation in the underlying process.

The Nadaraya-Watson estimator, however, is prone to boundary bias and instability near the edges of the time domain. As kernel weights become asymmetric, the estimator understates local slopes and exhibits inflated variance, particularly near points of structural change. This limitation motivates the use of higher--order local approximations.

Our first contribution extends the analysis of \cite{Vogt_2012} by developing a multivariate locally linear estimator in both rescaled time \(u\) and covariates \(x\). The method applies a first-order Taylor expansion within each local window so that kernel weighting corrects boundary bias and stabilizes the local design. While preserving the same stochastic order as the Nadaraya--Watson estimator, the locally linear approach yields sharper deterministic bias control and uniform risk bounds over \([0,1]\times S\). These improvements ensure stable performance even near temporal boundaries, where standard kernel methods typically deteriorate.

While the locally linear refinement addresses key limitations of kernel smoothing, many contemporary datasets present a different challenge: the target series is often sparse, whereas related auxiliary records are dense and highly informative. Examples include domestic weekly data paired with international daily measurements, or coarse macroeconomic indicators supplemented by high frequency covariates. Such imbalanced structures naturally motivate \emph{transfer learning}, which seeks to improve estimation in a small \emph{target} domain by using information from one or more related \emph{sources}. This raises a basic question of when information can be reliably transferred across heterogeneous domains.

A substantial literature has examined this problem under independent sampling. In nonparametric classification, minimax rates were characterized using a \emph{transfer exponent} that measures the discrepancy between source and target distributions (\cite{CaiT.Tony2021TRANSFERLEARNINGFORNONPARAMETRICCLASSIFICATION}; \cite{ReeveAdaptivetransferlearning_2021.}), and the framework now includes adversarially contaminated or unreliable sources (\cite{FanJianqing2025Rtlw}). For nonparametric regression, recent work established nonasymptotic minimax rates and adaptive procedures describing when transfer improves accuracy (\cite{Tony_Cai_2024}). In high dimensional linear and generalized linear models, joint estimators achieving optimal efficiency under sparse coefficient differences have been developed (\cite{LiSai2022Tlfh}; \cite{LiSaiTonyCai2024}; \cite{Tian02102023}). Further contributions include augmented estimators with semiparametric nuisance components (\cite{LiuMolei2023ATRL}) and model averaging methods attaining semiparametric efficiency bounds (\cite{HuXiaonan2023OPLb}). Extensions to functional and robust settings have also appeared, such as functional classification (\cite{QinCaihong2025AATL}) and contaminated source environments (\cite{FanJianqing2025Rtlw}). These results form a solid theoretical basis for transfer learning in independently sampled data.

Most theoretical work, however, still relies on independence across observations. Time series data feature serial correlation, nonstationarity, and gradual temporal drift that alter the standard bias-variance tradeoff. It is therefore unclear whether transfer learning can deliver similar benefits in dependent environments. This paper develops a framework for transfer learning with temporal dependence by integrating locally stationary process theory with nonparametric regression. We first establish uniform convergence for a multivariate locally linear estimator under standard mixing and smoothness conditions, achieving sharper bias and variance control than the classical Nadaraya-Watson estimator. We then introduce a bias corrected transfer estimator that links a short target sample and a long related source sample through a smooth bias function, allowing a locally smoothed adjustment of cross-domain differences. The resulting asymptotic risk bounds exhibit a phase transition: transfer yields improvement when the discrepancy evolves smoothly relative to the bandwidth, whereas substantial heterogeneity limits or nullifies the gain. To our knowledge, these results provide the first theoretical guarantees for transfer learning in nonparametric regression with dependent and locally stationary time series.

In summary, this paper develops a unified approach to transfer learning for nonparametric regression under temporal dependence. We first establish asymptotic results for a multivariate locally linear estimator that sharpens the analysis of \cite{Vogt_2012}. On this basis we construct a bias corrected transfer estimator that connects a short target series and a long related source series through a smooth bias function and permits adaptive borrowing of information across heterogeneous domains. The derived risk bounds identify when transfer learning is effective and demonstrate a distinct phase transition driven by the smoothness of cross domain discrepancy. Simulation studies and an empirical analysis of international fuel prices support the theoretical conclusions and illustrate the practical relevance of the method. Section 2 introduces the locally stationary regression setting, Section 3 develops the asymptotic theory, Section 4 presents the transfer framework and main results, and Sections 5 and 6 report simulation and empirical analyses. Technical proofs appear in the Appendix.

\section{Kernel estimation for locally stationary time series}
\label{sec-kernel}
\subsection{Local Stationarity}
A key element of the locally stationary framework is the rescaling of calendar time by the sample size, \(u = t/T \in [0,1]\). This transformation plays several roles. It compresses an expanding time index onto a fixed interval, allowing asymptotic arguments to be formulated on a compact domain. It also aligns the time dimension with kernel smoothing, since neighborhoods are defined around points \(u\) and a bandwidth \(h\) corresponds to roughly \(Th\) observations in actual time. The rescaling further provides the scale at which Riemann-sum approximations hold, so that quantities of the form
\[
    \frac{1}{T}\sum_{t=1}^T g\!\left(\tfrac{t}{T}\right)
\]
converge to meaningful integrals. Without this normalization, kernel weights and convergence rates would
not behave consistently as the sample grows.

Local stationarity is the second component of the framework. Rather than requiring distributional features to remain constant over the entire sample, one allows the process to be approximated over short segments of rescaled time by stationary laws that vary smoothly with \(u\). This reflects the empirical fact that means, variances, and dependence patterns evolve gradually as economic or technological conditions shift. The smooth drift along \(u\) provides the flexibility needed to capture long-run dynamics while retaining the
analytical structure required for local approximations.

\begin{definition}[Local Stationarity; cf.\ \cite{Vogt_2012}]
The $d$-dimensional process $\{X_{t,T}: t = 1,2,\dots,T\}_{T=1}^{\infty}$ is locally stationary if for every $u\in[0,1]$ there exists a strictly stationary process $\{X_t(u)\}$ with density $f_{X_t}(u)$ such that
\begin{equation}\label{locstat}
    \|X_{t,T}-X_t(u)\| \;\le\; \Big(\Big|t/T-u\Big|+1/T\Big)\,U_{t,T}(u)\quad\text{a.s.},
\end{equation}
where $\{U_{t,T}(u)\}$ are positive random variables with uniformly bounded $\rho$th moments, $\E[(U_{t,T}(u))^{\rho}] < C$ for some $\rho>0$ and constant $C < \infty$. The norm $\| \cdot\|$ is arbitrary on $\bbR^d$. 
\end{definition}

The parameter $\rho$ quantifies how sharply the stationary approximation concentrates around $u$; larger values correspond to tighter local fits. The combination of rescaled time and local stationarity thus provides both the mathematical tractability to derive uniform convergence rates to capture gradual structural change. 

\section{Multivariate Locally Linear Estimation for locally stationary time series}\label{sec-multill}
In this section, we consider multivariate locally linear estimation in the general nonparametric time-varying model
\begin{align}\label{eq:model}
    Y_{t,T} = m\!\left( t/T , X_{t,T} \right) + \varepsilon_{t,T},
\end{align}
with $\mathbb{E}[\varepsilon_{t,T}\mid X_{t,T}] = 0$, where $Y_{t,T}$ and $X_{t,T}$ are random variables of dimension $1$ and $d$, respectively. The covariate process $X_{t,T}$ is assumed to be locally stationary, and the regression function $m$ is allowed to vary smoothly over time. Throughout, $m$ depends on the rescaled time $u=t/T$ rather than on the calendar time $t$. 

Classical time series analysis often rests on the assumption of global stationarity, although many empirical series exhibit gradual changes in their distributional characteristics. A locally stationary formulation provides a more realistic description by assuming that over short windows of rescaled time ($u = t/T \in [0,1]$) the process can be approximated by a stationary one, while its law and the regression function ($m(u,x)$) may evolve smoothly in ($u$). This perspective has been developed for parametric time varying AR and ARCH models in \cite{DahlhausRainer2006SIfT} and for nonparametric regression settings in \cite{Vogt_2012}.

\subsection{Estimation Procedure}
We observe the triangular array $\{(Y_{t,T},X_{t,T}): t=1,\dots,T\}$ from model \eqref{eq:model}, where $X_{t,T}=(X_{t,T}^1,\dots,X_{t,T}^d)^\top\in\mathbb{R}^d$ and $m(u,x)$ may vary smoothly with the rescaled time $u\in[0,1]$ and covariates $x=(x^1,\dots,x^d)^\top \in \mathbb{R}^d$. 
Following the kernel smoothing framework of \cite{Vogt_2012}, the key idea is to smooth jointly in time and covariates, so that local information around $(u,x)$ is borrowed both across nearby time points and nearby covariate values. 

For later use, denote by $\nabla_0(m)(u,x)$ the partial derivative of $m$ with respect to the time coordinate, and by $\nabla_j(m)(u,x)$ the partial derivative with respect to the $j$-th coordinate of $x$ $(j=1,\dots,d)$:
\[
    \nabla_0(m)(u,x) := \frac{\partial}{\partial u} m(u,x), 
    \qquad 
    \nabla_j(m)(u,x) := \frac{\partial}{\partial x^j} m(u,x).
\]

$K$ denotes a one-dimensional kernel function which is bounded, symmetric, lipschitz continuous kernel with compact support and write $K_ h(v)= \frac{1}{h} K(v/h)$ for bandwidth $h>0$. We use a product kernel $K_h (u - t/T) \cdot\prod_{j = 1}^d K_{h}(x^j - \XtT^j)$, for smoothing over $(t/T, \, \XtT)$ around $(u,\,x)$. For convenience, assume that the bandwidth $h$ is the same in each direction. This results can be modified to allow different bandwidths for different direction. We smooth in the direction of the covariates $\XtT$ but also in the time direction. 

Then, the multivariate locally linear estimator of $m, \nabla_0(m), \cdots \nabla_d{(m)}$ is given by 
\begin{equation*}
    \begin{aligned}
        \hat{\frakm}(u,x)&:=\left(\hat{m} (u, x), h \cdot\hat{\nabla}_0 (m)(u,x), h \cdot\hat{\nabla}_1 (m)(u,x),\, \dots \, ,  h \cdot\hat{\nabla}_d (m)(u,x) \right)^{\top} \\& =  \underset{\mbf{a \in \bbR^{d+2}}}{\mathrm{argmin}} \Big( \frac {1} {T} \Big[ \sum_{t =1}^{T} \Bigg(Y_{t ,T} - a_0  - a_1 \left(\frac{t/T - u}{h} \right)- a_2\left(\frac{X_{t, T}^1 - x^1}{h}\right) - \\ & \qquad \qquad \qquad \qquad  \cdots - a_{d+1}  \left(\frac{X_{t, T}^d - x^d}{h}\right)\Bigg)^2 K_{h}(u - \frac{t} {T}) \prod_{j=1}^d K_{h}(x - X_{t, T}) \Big]\Big).
    \end{aligned}
\end{equation*}
In contrast to the Nadaraya–Watson estimator, the locally linear (LL) estimator requires explicit $h$-scaling. This arises from its construction via a local Taylor expansion, where slope terms are rescaled by $(x_j - X_{j,t,T})/h$ to match the local neighborhood defined by the bandwidth. The scaling matches the leading order of the bias of $m$ and $\nabla_j m$. Hence, the $h$-scaling is intrinsic to the LL estimator and explains its improved bias properties and boundary adaptation. 

\subsection{Assumptions}
We assume the following six conditions.

(C1) Local Stationarity of $X_{t,T}$:\\
For each $u\in[0,1]$, there exists a stationary approximation $\{X_t(u)\}$ that satisfies \eqref{locstat} with $E[(U_{t,T}(u))^\rho]\le C_U$ for some $\rho >0$. 

(C2) Smoothness of the Design Density:\\
The densities $f(u,x):= f_{X_t(u)}(x)$ depending in $t$ is differentiable with respect to $u$ and each $x^j$, with continuous partial derivative $\nabla_0 f(u,x) = \frac{\partial}{\partial u}f(u,x)$, $\nabla_{j} f(u,x):= \frac{\partial}{\partial x_j } f(u,x)$ for $j= 1, 2, \dots, d$.

(C3) Strong mixing (uniform polynomial decay).
Let \(W_{t,T}:=(X_{t,T},\varepsilon_{t,T})\) and define, for \(k\ge1\),
\[
\alpha(k)
:=\sup_{T\ge1}\ \sup_{1\le t\le T-k}
\alpha\!\Big(\sigma(W_{s,T}:1\le s\le t),\ \sigma(W_{s,T}:t+k\le s\le T)\Big).
\]
where \(\alpha(\mathcal{B},\mathcal{C}) = \sup_{B\in\mathcal{B}, C\in\mathcal{C}} | \Prob(B \cap C) - \Prob(B)\Prob(C) |\).
Assume the triangular array \(\{W_{t,T}\}\) is strongly mixing with $\alpha(k)\ \le\ A\,k^{-\beta}$ where $(k\ge1)$ for some constants \(A<\infty\) and \(\beta>0\) large enough for the results below.  
The same bound holds uniformly for the stationary proxies \(\{W_t(u)\}\), \(u\in[0,1]\), and for any measurable transforms of \(W_{t,T}\) used in the analysis.

(C4) Smoothness of  $m(u,x)$:\\
The regression function $m(u,x)$ is twice continuously partially differentiable in $(u,x)$; that is, its first and second partial derivatives exist and are continuous.

(C5) Kernel Properties:\\
The kernel $K:\mathbb{R}\to\mathbb{R}$ is symmetric about $0$, bounded, Lipschitz continuous, nonnegative, and has compact support: there exists $C_1<\infty$ such that $K(v)=0$ for $|v|>C_1$. In particular, there exist constants $C_K, C_{K,\mathrm{Lip}}<\infty$ such that
\begin{equation}\label{assKernel}
    \sup_{v\in\mathbb{R}} |K(v)| \le C_K, \qquad |K(v)-K(w)| \le C_{K,\mathrm{Lip}}\,|v-w| \quad \forall v,w\in\mathbb{R}.
\end{equation}
Assume $S\subset\mathbb{R}^d$ be compact. There exist constants $c_S>0$ such that
\begin{equation}\label{CS}
    \inf_{x\in S}\, \liminf_{h\downarrow 0}\, h^{-d}\,\mathrm{Leb}\!\Big(S\cap\big(x+(-C_1 h,\,C_1 h)^d\big)\Big)\ \ge\ c_S.
\end{equation}
where Leb denotes the $d$-dimensional Lebesgue measure.

We restrict the covariate domain to an arbitrary but fixed compact set $S$ and state all sup-norm results uniformly over $[0,1]\times S$. This compact-design assumption is standard in locally stationary kernel regression and is used to derive uniform bounds for kernel averages in $[0,1] \times S$ and to ensure well-behaved design densities in \cite{Vogt_2012}. 

Finally, throughout this paper, the bandwidth $h$ is assumed to converge to zero at least at polynomial rate, which mean there exists a small $\xi > 0 $ such that $ h \le CT^{-\xi}$    for some constant  $C \, > \, 0$.\\

\subsection{Uniform Convergence rates for kernel averages}
Define the $T \times (d+2)$ matrix $\mathbf{D} = (D_{t,j})_{1 \le t \le T,\; 1 \le j \le d+2}$ by
\begin{equation}\label{defD}
    D_{t,1} = 1, \qquad D_{t,2} = \frac{t/T - u}{h}, \qquad D_{t,j+2} = \frac{X^{\,j}_{t,T} - x_j}{h}, \quad j = 1,\ldots,d.
\end{equation}
Each row of $\mathbf{D}$ corresponds to the local regressors associated with the observation $(t/T, X_{t,T})$, scaled by the bandwidth $h$.

Define further the diagonal weight matrix
\begin{equation}\label{defW}
    {\mbf{W}} = \mathrm{diag} \Big( K_{h}(u - t/T) \prod_{j=1}^{d} K_{h}(x^j - X_{t, T}^j) \Big).
\end{equation}
which assigns kernel weights to observations according to their proximity to $(u,x)$ in both time and covariate space.

Then, writing \(\mbf{Y}_T = (Y_{1, T}, \dots , Y_{T, T})^{\top},\; \mbf{X}_T = (X_{1, T}, \dots , X_{T, T})^{\top} \), where $Y_{j,T} \in \bbR$, $X_{j,T} \in \bbR^d$. 
\begin{align}\label{frakm}
    \hat{\frakm}(u,x) &:= \left( \hat{m} (u, x), h \cdot \hat{\nabla}_0{m}(u,x), h \cdot \hat{\nabla}_1 m(u,x), \dots , h \cdot \hat{\nabla}_d m(u,x) \right)^{\top} \\
    &= (\mbf{D}^{\top} \mbf{W} \mbf{D})^{-1} \mbf{D}^{\top} \mbf{W} \mbf{Y}_T.
\end{align}

Define \(\frakm(u,x)\) and \(\mbf{m}^{\top}\)by 
\begin{equation*}
    \begin{aligned}
        &\frakm(u,x) := \left( {m} (u, x), h \cdot {\nabla}_0{m}(u,x), h \cdot {\nabla}_1 m(u,x), \dots , h \cdot {\nabla}_d m(u,x) \right)^{\top}\\ 
        &\mbf{m} = \left(m\left( {1}/{T},\,X_{1, T}\right), \dots , m\left({T}/{T},\, X_{T, T}\right)\right)^{\top}
    \end{aligned}
\end{equation*}
First, we examine kernel averages of the general form
\begin{align}
    \hat{\Psi}(u,x) = \frac 1 {T} \mbf{D}^{\top} \mbf{W} \mbf{R}_{T}
\end{align}
with $\mbf{R}_{T}$ being an array of one-dimensional random variables $\mbf{R}_T = (R_{1, T}, \dots , R_{T, T})^{\top}$. 

We now derive the uniform convergence rate of $\hat{\Psi}(u,x)-\E[\hat{\Psi}(u,x)]$. 
To do so, we make the following assumptions on the components in $R_{t,T}$. 

(K1) It holds that $E\Big[|R_{t,T}|^s\Big] \le C$ for some $s>2$ and a constant $C < \infty$. 
    
(K2) The array $\{X_{t,T},R_{t,T}\}$ is $\alpha$–mixing with mixing coefficients $\alpha$ satisfying $\alpha(k) \le A\,k^{-\beta}$ for some constant $A < \infty$ and  $\beta > \frac{2s-2}{s-2}$ 
    
(K3) Let $f_{X_{t,T}}$ be the density of $X_{t,T}$ and $f_{X_{t,T},X_{t+l,T}}$ be the joint density of $(X_{t,T}, X_{t+l,T})$. For any compact set $S\subseteq \bbR^d$ there exists a constant depending on $S$, $C = C(S)$ such that $\sup_{t,T} \sup_{x\in S} f_{X_{t,T}}(x) \le C$, $\sup_{t,T}\sup_{x\in S} E\Big[|R_{t,T}|^s \mid X_{t,T}=x\Big] f_{X_{t,T}}(x) \le C$. Moreover, there exists an integer $l^* < \infty$ such that for all $l \ge l^*$, $$\sup_{t,T} \sup_{x,x'\in S} E\Big[|R_{t,T}\|R_{t+l,T}|\mid X_{t,T}=x,\,X_{t+l,T}=x'\Big] f_{X_{t,T},X_{t+l,T}}(x,x') \le C$$

The following theorem generalizes uniform convergence results of \cite{Hansen_2008} and \cite{Vogt_2012}.
\begin{theorem}\label{Thm3.1}
Assume that (K1)–(K3) are satisfied with
\begin{align}\label{assbeta}
    \beta > \frac{2 +s(1+(d+1))}{s-2} \qquad \text{implying} \; \beta>d+2
\end{align}
and that the kernel $K$ fulfills (C5). Moreover, assume the bandwidth $h$ satisfy
\begin{align}\label{asstheta}
    \frac{ \log\log T \cdot \log T}{T^{\theta} \, h^{d+1}} = o(1), \qquad  \theta = \frac{\beta(1- 2/ s ) -  2 / s -1 - (d+1)}{\beta + 3-(d+1)} >0.
\end{align}
Finally, let $S$ be a compact subset of $\mathbb{R}^d$ that satisfies \eqref{CS}. Then it holds that
\begin{align}
    \sup_{u\in[0,1],\,x\in S} \Norm{\hat{\Psi}(u,x)-\E[\hat{\Psi}(u,x)]}{2} = O_P\Big(\sqrt{\frac{\log T}{T\,h^{d+1}}}\Big).
\end{align}
where $\|\cdot\|_{2}$ is euclidean norm defined in $\bbR^{d+2}$.
\end{theorem}

\subsection{Uniform convergence rates for multivariate locally linear estimates}

\begin{theorem}\label{Thm3.2}
Assume (C1)–(C5) and that (K1)–(K3) hold both for $R_{t,T}=1$ and ${R}_{t,T}=\varepsilon_{t,T}$. Assume $\beta$ satisfy the inequality \eqref{assbeta}. $S$ be a compact set that satisfies \eqref{CS}. Suppose that $\inf_{u\in[0,1],\,x\in S} f(u,x) > 0$. Moreover, assume that the bandwidth $h$ satisfies
\begin{align}\label{assh}
    \frac{\log \log T \,\log T}{T^{\theta}\,h^{d+1}} = o(1) \quad \text{and} \quad \frac{1}{T^r\,h^{d+r}} = o(1)
\end{align}
with $r = \min\{\rho,1\}$. Then,
\begin{align}
    \sup_{u\in [0,1],\,x\in S} \Norm{\hat{\frakm}(u,x) - \frakm(u,x)}{2}
    = O_P\Big(h^2 + \sqrt{\frac{\log T}{T\,h^{d+1}}} + \frac{1}{T^r h^{d-1}}\Big).
\end{align}
\end{theorem}

\medskip \noindent
Theorem~3.2 clarifies how local stationarity influences the convergence behavior of nonparametric regression estimators. In addition to the usual stochastic and smoothing components, the estimation error contains a remainder of order \(O(T^{-r}h^{-(d-1)})\), generated by the gradual evolution of the data--generating mechanism. This term reflects the extent to which the process can be approximated by a stationary law within shrinking neighborhoods of rescaled time. When the local approximation is accurate (large \(\rho\)), the remainder is dominated by the standard kernel terms and the estimator behaves much as it does under stationarity. When the approximation weakens (small \(\rho\)), the remainder governs the overall rate and the estimator converges more slowly. The result therefore quantifies how the smoothness of temporal variation determines the effective sample size in locally stationary settings, linking classical stationary theory with nonstationary dynamics within a unified asymptotic framework.

\section{Transfer learning framework}\label{sec-TL}
\subsection{Model Setup}
We consider a transfer learning problem in nonparametric regression with two triangular arrays of observations: one from a small \emph{target} sample and the other from a large \emph{source} sample. For each domain $\ell \in \{0,1\}$, let $(Y_{t_\ell,T_\ell}^{(\ell)}, X_{t_\ell,T_\ell}^{(\ell)})$, $t_\ell = 1, \dots, T_\ell$, denote the observations generated from
\begin{align}
    Y_{t_\ell, T_\ell}^{(\ell)} 
    &= m^{(\ell)}\!\left(\frac{t_\ell}{T_\ell},\, X_{t_\ell, T_\ell}^{(\ell)}\right) + \varepsilon_{t_\ell, T_\ell}^{(\ell)}, 
    \qquad \E\!\left[\varepsilon_{t_\ell, T_\ell}^{(\ell)} \mid X_{t_\ell, T_\ell}^{(\ell)}\right] = 0, 
    \quad t_\ell = 1, \dots, T_\ell,
    \label{eq:model_transfer}
\end{align}
where $m^{(\ell)}$ denotes the domain-specific regression function, $X_{t_\ell,T_\ell}^{(\ell)} \in \mathbb{R}^d$ is a $d$-dimensional covariate vector, and $\varepsilon_{t_\ell,T_\ell}^{(\ell)}$ is a mean-zero random error.

For notational convenience, define the vectors
\[
\mbf{Y}_{T_\ell}^{(\ell)} = \big(Y_{1,T_\ell}^{(\ell)}, \dots, Y_{T_\ell,T_\ell}^{(\ell)}\big)^{\top},
\quad
\mbf{X}_{T_\ell}^{(\ell)} = \big(X_{1,T_\ell}^{(\ell)}, \dots, X_{T_\ell,T_\ell}^{(\ell)}\big)^{\top},
\quad
\mbf{m}^{(\ell)} = \big(m^{(\ell)}(\tfrac{t_\ell}{T_\ell}, X_{t_\ell,T_\ell}^{(\ell)})\big)_{t_\ell=1}^{T_\ell}.
\]
We also let $\hat{\mbf{m}}^{(1)} = \big(\hat{m}^{(1)}({t_0}/{T_0}, X_{t_0,T_0}^{(0)})\big)_{t_0=1}^{T_0}$ denote the fitted source-domain regression function evaluated on the target sample for subsequent bias correction.

The goal is to estimate the target regression function $m^{(0)}(u,x)$ from a small sample of size $T_0$ by leveraging information from a related source sample of size $T_1 \gg T_0$. 
To characterize the cross-domain discrepancy, define the bias function
\begin{equation}\label{defb}
    b(u,x) = m^{(0)}(u,x) - m^{(1)}(u,x),
\end{equation}
which measures the difference between the target and source regression surfaces. 
We further quantify its local smoothness by
\begin{equation}\label{defeta}
    \sup_{u,x} \|\nabla b(u,x)\|_{2} = \eta_{1,b}, \qquad \sup_{u,x} \|\nabla^{2} b(u,x)\|_{F} = \eta_{2,b},
\end{equation}
where $\nabla b(u,x)$ and $\nabla^{2} b(u,x)$ denote the gradient and Hessian of $b(u,x)$ with respect to $(u,x)$, and $\eta_{1,b},\, \eta_{2,b} \le 1$ represent bounded smoothness constants controlling the first- and second-order variability of the bias function.

Both datasets are observed over the same rescaled time interval $[0,1]$ but differ in sampling resolution. 
The target sample $\{(X^{(0)}_{t_0,T_0}, Y^{(0)}_{t_0,T_0})\}_{t_0=1}^{T_0}$ is recorded at coarse grid points $\{t_0/T_0\}_{t_0=1}^{T_0}$, whereas the source sample $\{(X^{(1)}_{t_1,T_1}, Y^{(1)}_{t_1,T_1})\}_{t_1=1}^{T_1}$ is observed more densely at $\{t_1/T_1\}_{t_1=1}^{T_1}$ with $T_1 \gg T_0$. 
The finer temporal resolution of the source domain provides auxiliary information that can be exploited to improve estimation of the target regression function $m^{(0)}(u,x)$.

Throughout, the covariates $X^{(\ell)}_{t_\ell,T_\ell}$ are assumed to lie in a compact subset of $\mathbb{R}^d$, and the regression functions $m^{(\ell)}(u,x)$ are smooth in both arguments for $\ell = 0,1$. 
Detailed smoothness and alignment conditions linking $m^{(0)}$ and $m^{(1)}$ are introduced in the next section.

\subsection{Assumptions for transfer learning}
For the transfer learning framework, we assume that both the target and source processes $\{(X^{(0)}_{t_0,T_0}, \varepsilon^{(0)}_{t_0,T_0})\}$ and $\{(X^{(1)}_{t_1,T_1}, \varepsilon^{(1)}_{t_1,T_1})\}$ satisfy the regularity conditions (C1)–(C5) introduced in Section~3.2. 
In particular, each covariate process $X^{(\ell)}_{t_\ell,T_\ell}$, $\ell \in \{0,1\}$, is locally stationary with stationary approximation $\{X^{(\ell)}_{t_\ell}(u)\}$ obeying
\[
\mathbb{E}\big[(U^{(\ell)}_{t_\ell,T_\ell}(u))^{\rho_\ell}\big] < C, 
\qquad \rho_\ell > 0, \; C < \infty,
\]
uniformly in $u, t_\ell, T_\ell$. 
The corresponding design densities $f^{(\ell)}(u,x)$ are continuously differentiable in both $u$ and $x$, and the arrays are uniformly $\alpha$–mixing with polynomially decaying coefficients. 
The regression functions $m^{(0)}(u,x)$ and $m^{(1)}(u,x)$ are twice continuously differentiable in $(u,x)$, and the kernel $K$ satisfies the boundedness, Lipschitz continuity, and compact-support properties specified in (C5). 

We further assume that the pairs $(X_{t,T}, R_{t,T}) \in \{(X^{(0)}_{t_0,T_0},1),\, (X^{(0)}_{t_0,T_0},\varepsilon^{(0)}_{t_0,T_0}),\, (X^{(1)}_{t_1,T_1},1),\, (X^{(1)}_{t_1,T_1},\varepsilon^{(1)}_{t_1,T_1})\}$ satisfy conditions (K1)–(K3) stated in Section~3. 
In particular, the arrays are uniformly $\alpha$–mixing with polynomially decaying coefficients as in~\eqref{assbeta}. 
To ensure well-behaved kernel averages, we assume the design densities are strictly positive on the support,
\begin{equation}\label{assden}
    \inf_{u\in[0,1],\,x\in S} f(u,x) > 0, \qquad \inf_{u\in[0,1],\,x\in S} g(u,x) > 0,
\end{equation}
where $S \subset \mathbb{R}^d$ is a compact set satisfying~\eqref{CS}.  

Bandwidth sequences $(h_0,h_1,h_{\tl})$ obey the rate conditions in~\eqref{assh} and~\eqref{asstheta}, with $r = \min\{\rho_1,\rho_2,1\}$. 
Each bandwidth is assumed to shrink at least polynomially, that is, there exist constants $\xi > 0$ and $C < \infty$ such that
\[
h_0,\, h_{\tl} \le C T_0^{-\xi},
\qquad 
h_1 \le C T_1^{-\xi}.
\]
These assumptions guarantee uniform convergence of kernel estimators under locally stationary dependence across both domains.

\subsection{Transfer Learning for NW estimates}
In this section, we introduce a transfer learning framework for the Nadaraya–Watson estimator in a time-varying nonparametric regression setting. 
The analysis is conducted under the assumption of local stationarity, where smoothing is performed jointly in rescaled time $u = t/T$ and covariates $x$. 
The data consist of a target sample $\{(Y^{(0)}_{t_0,T_0}, X^{(0)}_{t_0,T_0})\}$ and an auxiliary source sample $\{(Y^{(1)}_{t_1,T_1}, X^{(1)}_{t_1,T_1})\}$. 
We begin by constructing a time-varying Nadaraya–Watson estimator on the source domain, then locally estimate the cross-domain bias function $b(u,x) = m^{(0)}(u,x) - m^{(1)}(u,x)$ using residuals from the target data, and finally combine the two through bias correction. 
This procedure retains the local temporal structure of the target process while reducing estimation variance. 
Under assumptions (C1)–(C5) and (K1)–(K3), the resulting estimator achieves faster convergence rates than a target-only Nadaraya–Watson estimator. 
The following subsections formalize the estimator, outline the required conditions, and analyze its asymptotic properties to delineate the situations in which transfer learning yields improvement and those in which it does not.

\subsubsection{Estimation Procedure}
\noindent\textbf{Step 1. Source-domain estimation.}
On the large source sample $\{(X^{(1)}_{t_1,T_1}, Y^{(1)}_{t_1,T_1})\}_{t_1=1}^{T_1}$, the time-varying Nadaraya–Watson estimator of the source regression function $m^{(1)}(u,x)$ is defined as
\begin{equation}\label{eq:nw_source}
    \hat{m}^{(1)}(u,x) = 
    \frac{\sum_{t_1=1}^{T_1} K_{h_1}\!\left(u - \tfrac{t_1}{T_1}\right) \prod_{j=1}^{d} K_{h_1}\!\left(x^{j} - X^{j,(1)}_{t_1,T_1}\right) Y^{(1)}_{t_1,T_1}}{\sum_{t_1=1}^{T_1} K_{h_1}\!\left(u - \tfrac{t_1}{T_1}\right) \prod_{j=1}^{d} K_{h_1}\!\left(x^{j} - X^{j,(1)}_{t_1,T_1}\right)}.
\end{equation}

\noindent\textbf{Step 2. Bias estimation.}
The cross-domain bias $b(u,x) = m^{(0)}(u,x) - m^{(1)}(u,x)$ is estimated by locally smoothing the target residuals relative to $\hat{m}^{(1)}$, namely,
\begin{equation}\label{eq:bias_est}
    \hat{b}(u,x) =
    \frac{\sum_{t_0=1}^{T_0} K_{h_{\mathrm{TL}}}\!\left(u - \tfrac{t_0}{T_0}\right) \prod_{j=1}^{d} K_{h_{\mathrm{TL}}}\!\left(x^{j} - X^{j,(0)}_{t_0,T_0}\right) \big(Y^{(0)}_{t_0,T_0} - \hat{m}^{(1)}(\tfrac{t_0}{T_0}, X^{(0)}_{t_0,T_0})\big)} {\sum_{t_0=1}^{T_0} K_{h_{\mathrm{TL}}}\!\left(u - \tfrac{t_0}{T_0}\right) \prod_{j=1}^{d} K_{h_{\mathrm{TL}}}\!\left(x^{j} - X^{j,(0)}_{t_0,T_0}\right)}.
\end{equation}

\noindent\textbf{Step 3. Bias-corrected estimation.}
The transfer learning estimator combines the source fit with the estimated bias correction:
\begin{equation}\label{eq:mtl_final}
    \hat{m}^{\mathrm{TL}} (u,x) = \hat{m}^{(1)}(u,x) + \hat{b}(u,x).
\end{equation}
This construction preserves the local temporal structure of the target process while leveraging information from the source domain through a data-driven bias adjustment.

\subsubsection{Convergence Rate Analysis}
\begin{theorem}\label{Thm4.1}
Assume that conditions \textup{(C1)}–\textup{(C5)} and \textup{(K1)}–\textup{(K3)} hold for 
\[
(X_{t,T}, R_{t,T}) \in 
\{(X^{(0)}_{t_0,T_0},1),\, (X^{(0)}_{t_0,T_0},\varepsilon^{(0)}_{t_0,T_0}),\, 
  (X^{(1)}_{t_1,T_1},1),\, (X^{(1)}_{t_1,T_1},\varepsilon^{(1)}_{t_1,T_1})\}.
\]
Let $r = \min(\rho_1,\rho_2,1)$ as in \textup{(C1)}, and suppose that the bandwidths $h_0$ and $h_{\mathrm{TL}}$ satisfy \eqref{assh} for some $\theta>0$. Assume that $m^{(0)}$, $m^{(1)}$, and $b(u,x) = m^{(0)}(u,x) - m^{(1)}(u,x)$ are twice continuously differentiable with bounded derivatives. Let $I_h = [C_1h_1, 1 - C_1h_1] \cap [C_1h_{\mathrm{TL}}, 1 - C_1h_{\mathrm{TL}}]$, and suppose that $\nabla^{3}b(u,x)$ is bounded. Assume further that the sample sizes satisfy $T_1 \gg T_0$, and that $S \subset \mathbb{R}^d$ is compact and satisfies~\eqref{CS}. Then the transfer learning estimator obeys the uniform bound
\begin{equation}\label{eqthm4.1} 
    \begin{aligned} 
        \sup_{u \in I_h,\,x \in S}& \Big|\hat{m}^{\tl}(u,x) - m^{(0)}(u,x)\Big| \\ &= O_P\left(\sqrt{\frac{\log T_0}{T_0 h_{\tl}^{d+1}}} + \frac{\eta_{1, b}}{T_0^r h_{\tl}^d} +\eta_{2, b} h_{\tl}^2 + \sqrt{\frac{\log T_1}{T_1 h_1^{d+1}}} + \frac{1}{T_1^r h_1^d} + h_1^2 \right). 
    \end{aligned} 
\end{equation}
\end{theorem}


Theorem~\ref{Thm4.1} describes how local stationarity and cross-domain similarity together determine the convergence behavior of the transfer learning estimator. In addition to the usual stochastic and smoothing components, the rate includes remainder terms of orders \(O(T_0^{-r}h_{\mathrm{TL}}^{-d})\) and \(O(T_1^{-r}h_1^{-d})\), reflecting the fact that both domains satisfy only local stationarity. These terms represent the cost of temporal evolution. When each process admits a sharp stationary approximation (large \(\rho_\ell\)), the remainder becomes negligible and the estimator attains the stationary benchmark rate; when the approximation deteriorates, the remainder governs the attainable accuracy.

The smoothness of the cross--domain bias function, summarized by \((\eta_{1,b}, \eta_{2,b})\), determines how effectively information can be transferred. Strong alignment in first and second order enables variance reduction without introducing additional bias. The theorem therefore clarifies how temporal dependence and smooth structural drift jointly affect the benefits of transfer learning and yields a unified rate characterization for dependent and nonstationary settings.

\subsection{Transfer Learning for Multivariate LL estimates}
\subsubsection{Estimation Procedure}
Building on the notation in~\eqref{defD}–\eqref{defW}, we now formalize the multivariate locally linear framework for transfer learning.
For source domain with sample size~$T_1$ and bandwidth~$h_1$, we collect the local regressors and kernel weights into the design and weight matrices $\mathbf{D}_1$ and~$\mathbf{W}_1$.
Specifically, for the source domains, the design matrix
\begin{equation}\label{defD1}
D_{t_1,1}^{(1)} = 1, \qquad
D_{t_1,2}^{(1)} = \frac{t_1/T_1 - u}{h_1}, \qquad
D_{t_1,j+2}^{(1)} = \frac{X_{t_1,T_1}^{j,(1)} - x_j}{h_1},
\quad j = 1,\ldots,d,
\end{equation}
contains in each row the local regressors associated with $(t_1/T_1,,X_{t_1,T_1}^{(1)})$.
The corresponding diagonal weight matrix
\begin{equation}\label{defW1}
\mathbf{W}_1 = \mathrm{diag}\left(
K_{h_1} \left(u - \tfrac{t_1}{T_1}\right)
\prod_{j=1}^{d} K_{h_1}\left(x^j - X_{t_1,T_1}^{j,(1)}\right)
: t_1 = 1,\ldots,T_1 \right)
\end{equation}
assigns product–kernel weights to each observation according to its proximity to $(u,x)$ in both rescaled time and covariate space.

\medskip
For the transfer–learning fit, we define analogous quantities $\mathbf{D}_{\tl}$ and~$\mathbf{W}_{\tl}$ constructed on the target sample with bandwidth~$h_{\tl}$ and sample size~$T_0$.
The design matrix is given by
\begin{equation}\label{defDtl}
D_{t_0,1}^{\tl} = 1, \qquad
D_{t_0,2}^{\tl} = \frac{t_0/T_0 - u}{h_{\tl}}, \qquad
D_{t_0,j+2}^{\tl} = \frac{X_{t_0,T_0}^{j,\tl} - x_j}{h_{\tl}},
\quad j = 1,\ldots,d,
\end{equation}
and the corresponding weight matrix
\begin{equation}\label{defWtl}
\mathbf{W}_{\tl} = \mathrm{diag}\left(
K_{h_{\tl}} \left(u - \tfrac{t_0}{T_0}\right)
\prod_{j=1}^{d} K_{h_{\tl}}\left(x^j - X_{t_0,T_0}^{j,\tl}\right)
: t_0 = 1,\ldots,T_0 \right).
\end{equation}
Each row of $\mathbf{D}{\tl}$ encapsulates the local linear regressors used for bias correction, and $\mathbf{W}{\tl}$ provides the corresponding adaptive kernel weights.
Together, these quantities constitute the fundamental building blocks for the transfer–learning locally linear estimator introduced below.

For each domain $\ell \in \{0,1\}$, define
\[
    \mathfrak{m}^{(\ell)}(u,x) = \big(m^{(\ell)}(u,x),\; h_\ell \nabla_0 m^{(\ell)}(u,x),\; \dots,\; h_\ell \nabla_d m^{(\ell)}(u,x) \big)^{\top},
\]
which collects the local intercept and the scaled first-order derivatives of the regression function $m^{(\ell)}(u,x)$ with respect to the rescaled time $u$ and covariates $x$. 
In particular, $\mathfrak{m}^{(0)}(u,x)$ and $\mathfrak{m}^{(1)}(u,x)$ correspond to the target and source domains, respectively.

\medskip \noindent \textbf{Step 1. Source-domain estimation.}
Using the source sample $\{(X_{t_1,T_1}^{(1)}, Y_{t_1,T_1}^{(1)})\}_{t_1=1}^{T_1}$, the local linear estimator of $m^{(1)}$ and its first-order derivatives is obtained as
\begin{equation}\label{eq:sourcefit}
    \hat{\mathfrak{m}}^{(1)}(u,x) = (\mathbf{D}_1^{\top}\mathbf{W}_1\mathbf{D}_1)^{-1}\mathbf{D}_1^{\top}\mathbf{W}_1\mathbf{Y}_{T_1}^{(1)}.
\end{equation}

\medskip
\noindent\textbf{Step 2. Bias estimation.}
Define the cross-domain bias function $b(u,x)=m^{(0)}(u,x)-m^{(1)}(u,x)$ and its augmented form 
\[
    \mathfrak{b}(u,x) = \big(b(u,x),\; h_{\tl}\nabla_0 b(u,x),\; \dots,\; h_{\tl}\nabla_d b(u,x) \big)^{\top}.
\]
The estimator of $\mathfrak{b}(u,x)$ is obtained by local smoothing of target residuals:
\begin{equation}\label{eq:biasfit}
    \hat{\mathfrak{b}}(u,x) = (\mathbf{D}_{\mathrm{TL}}^{\top}\mathbf{W}_{\mathrm{TL}}\mathbf{D}_{\mathrm{TL}})^{-1} \mathbf{D}_{\mathrm{TL}}^{\top}\mathbf{W}_{\mathrm{TL}} \big(\mathbf{Y}_{T_0}^{(0)} - \hat{\mathbf{m}}^{(1)}\big),
\end{equation}
where $\hat{\mathbf{m}}^{(1)}$ denotes the source-domain fit evaluated on the target sample.

\medskip
\noindent\textbf{Step 3. Bias-corrected estimation.}
Finally, the transfer learning estimator combines the bias-corrected components as
\begin{equation}\label{eq:TLest}
    \hat{\mathfrak{m}}^{\mathrm{TL}}(u,x) = \hat{\mathfrak{b}}(u,x) + \big(\hat{m}^{(1)}(u,x),\; h_{\mathrm{TL}}\hat{\nabla}_0 m^{(1)}(u,x),\; \dots,\; h_{\mathrm{TL}}\hat{\nabla}_d m^{(1)}(u,x)\big)^{\top},
\end{equation}
so that the target regression estimator is  $\hat{m}^{\mathrm{TL}}(u,x)=\hat{b}(u,x)+\hat{m}^{(1)}(u,x)$.

\subsubsection{Convergence Rate Analysis}
\begin{theorem}\label{Thm4.2}
Assume that conditions \textup{(C1)}–\textup{(C5)} and \textup{(K1)}–\textup{(K3)} hold for $(X_{t,T}, R_{t,T}) \in \{(X^{(0)}_{t_0,T_0},1), (X^{(0)}_{t_0,T_0},\varepsilon^{(0)}_{t_0,T_0}), (X^{(1)}_{t_1,T_1},1), (X^{(1)}_{t_1,T_1},\varepsilon^{(1)}_{t_1,T_1})\}$.
Let $r=\min(\rho_1,\rho_2,1)$ and suppose the bandwidths $h_0,h_{\mathrm{TL}}$ satisfy \eqref{assh}.
Assume that $m^{(0)}$, $m^{(1)}$, and the bias function $b(u,x)=m^{(0)}(u,x)-m^{(1)}(u,x)$ are three times continuously differentiable with bounded derivatives, that $\inf_{u,x}f(u,x)>0$, $\inf_{u,x}g(u,x)>0$, and that $S\subset\mathbb{R}^d$ is compact and satisfies~\eqref{CS}. Let $T_1\gg T_0$. Then, the transfer learning locally linear estimator satisfies
\begin{equation}\label{eqthm4.2} 
    \begin{aligned} 
        &\sup_{u\in[0,1],\,x\in S}\Big|\hat{m}^{\tl}(u,x)-m^{(0)}(u,x)\Big| \\ &= O_p\!\left(\sqrt{\frac{\log T_0}{T_0 h_{\tl}^{\,d+1}}} +\frac{\eta_{2,b}}{T_0^{r}h_{\tl}^{\,d-1}} +\eta_{2,b}h_{\tl}^{2}+ \sqrt{\frac{\log T_1}{T_1 h_1^{\,d+1}}}+\frac{1}{T_1^{r}h_1^{\,d-1}}+h_1^{2} \right).
    \end{aligned} 
\end{equation} 
\end{theorem}

Theorem~\ref{Thm4.2} sharpens the convergence analysis of Theorem~\ref{Thm4.1} by showing how the local linear structure improves performance in the presence of temporal nonstationarity. The stochastic fluctuation remains of the same order as for the Nadaraya--Watson estimator, but the deterministic remainder generated by local nonstationarity improves from \(O(T_0^{-r}h_{\mathrm{TL}}^{-d})\) to \(O(T_0^{-r}h_{\mathrm{TL}}^{-(d-1)})\). This gain reflects the ability of first--order local fitting to accommodate smooth temporal drift within the locally stationary regime.

The improvement stems from the locally linear expansion, which corrects temporal curvature and stabilizes the effective design near boundaries, yielding a more accurate approximation to the evolving data--generating mechanism. By incorporating derivatives in both rescaled time and covariates, the estimator aligns the local stationary approximation with the direction of smooth structural change, thereby reducing the accumulation of nonstationary bias. Theorem~\ref{Thm4.2} thus demonstrates that the locally linear formulation enhances robustness to temporal evolution and provides a sharper characterization of the attainable rate under locally stationary dependence.

\section{Simulation Studies}\label{sec-simulation}
\subsection{Simulation Design}
To illustrate the theoretical results developed in Section~4, we conduct a simulation study examining how the proposed transfer-learning estimators behave in finite samples under controlled domain discrepancies.  
We consider two types of samples: the \emph{given sample}, which represents the primary data available for estimation, and the \emph{auxiliary sample}, which provides additional information from related environments.  
For clarity, we subsequently refer to the given sample as the \emph{target domain} (indexed by $\ell = 0$) and to the auxiliary sample as the \emph{source domain} (indexed by $\ell = 1$).
For each domain we observe
\[
    Y^{(\ell)}_{t_\ell} = m^{(\ell)}\!\big(u^{(\ell)}_{t_\ell},x^{(\ell)}_{t_\ell}\big) + \eta^{(\ell)}_{t_\ell}, \qquad t_\ell=1,\ldots,T_\ell,
\]
where $u^{(\ell)}_{t_\ell}=t_\ell/T_\ell\in[0,1]$ is the rescaled time index and $\eta^{(\ell)}_{t_\ell}\stackrel{\text{i.i.d.}}{\sim}N(0,\tau_\ell^2)$ are independent observational errors. 
The regression surface in the target domain is specified as
\[
    m^{(0)}(u,x) = 2\sin(\pi x)\,(0.5u+2) + u(1-x) + 2,
\]
a smooth bivariate function with uniformly bounded second derivatives. 
This choice guarantees that the bias expansions derived in Section~4 apply directly, while preserving nontrivial interaction between $u$ and $x$.

The covariates $\{x^{(\ell)}_{t_\ell}\}$ are generated from a time-varying autoregressive process of order two,
\[
    x^{(\ell)}_{t_\ell} = a_1(u^{(\ell)}_{t_\ell})\,x^{(\ell)}_{t_\ell-1} + a_2(u^{(\ell)}_{t_\ell})\,x^{(\ell)}_{t_\ell-2} + s(u^{(\ell)}_{t_\ell})\,\zeta^{(\ell)}_{t_\ell}, \qquad  \zeta^{(\ell)}_{t_\ell}\stackrel{\text{i.i.d.}}{\sim}N(0,1),
\]
with $x^{(\ell)}_{-1}=x^{(\ell)}_0=0$. The coefficient functions are given by
\[
    a_1(u)=0.15\cos(\pi/3)+0.3\bigl((u-0.5)^2-\tfrac{1}{12}\bigr),\qquad
    a_2(u)=0.15\cos(2\pi/3)+0.3\bigl((u-0.5)^2-\tfrac{1}{12}\bigr),
\]
and the volatility curve is defined as
\[
    s(u)=\max\{\,0.10+0.10(0.5+0.5\sin(0.5\pi u)),\,10^{-3}\,\}.
\]
These functions vary smoothly in $u$ and satisfy $\sup_u(|a_1(u)|+|a_2(u)|)<1$, ensuring that the roots of $1-a_1(u)z-a_2(u)z^2$ remain outside the unit circle uniformly in $u$. 
Consequently, the process $\{x^{(\ell)}_{t_\ell}\}$ is stable and locally stationary in the sense of \cite{Dahlhaus.Richter.Wu_2019}. 
Each covariate series is subsequently linearly rescaled to $[0,1]$ so that the target and source domains have comparable support.

The target domain contains $T_0=2000$ observations generated as
\[
    Y^{(0)}_{t_0} = m^{(0)}(u^{(0)}_{t_0},x^{(0)}_{t_0}) + \eta^{(0)}_{t_0}, \qquad \eta^{(0)}_{t_0}\sim N(0,0.1^2),
\]
whereas the source domain consists of $T_1=20000$ observations drawn from a related but biased regression surface,
\[
    m^{(1)}(u,x) = m^{(0)}(u,x) + b_\gamma(u,x), \qquad Y^{(1)}_{t_1} = m^{(1)}(u^{(1)}_{t_1},x^{(1)}_{t_1}) + \eta^{(1)}_{t_1}, \quad \eta^{(1)}_{t_1}\sim N(0,0.1^2).
\]
The bias function $b_\gamma$ controls the structural difference between the two domains. We consider three smooth bias families of comparable curvature:
\[
    b_\gamma^{\text{quad}}(u,x) = \gamma\frac{u^2+x^2}{2},\qquad 
    b_\gamma^{\text{cubic}}(u,x) = \gamma\frac{u^3+x^3}{6},\qquad
    b_\gamma^{\text{exp}}(u,x) = \gamma\frac{e^u+e^x}{e}.
\]
Each satisfies $\sup_{u,x}\|\nabla^2 b_\gamma(u,x)\|_\infty = |\gamma|$, so $\gamma$ directly determines the curvature and intensity of the discrepancy while maintaining an identical smoothness order across bias types. 
Varying $\gamma$ over $\{-10,-8,\ldots,10\}$ produces a smooth continuum from nearly identical domains to strongly misspecified ones, enabling a clear visualization of the bias–variance transition predicted by the theory.

Four kernel estimators are compared. 
The first pair are the Nadaraya–Watson and locally linear estimators fitted solely on the target sample, denoted by $\widehat{m}_{\mathrm{NW\text{-}T}}$ and $\widehat{m}_{\mathrm{LL\text{-}T}}$. 
The second pair, $\widehat{m}_{\mathrm{NW\text{-}TL}}$ and $\widehat{m}_{\mathrm{LL\text{-}TL}}$, are their transfer-learning analogues, which incorporate information from the source domain through residual-based bias correction:
\[
    \widehat{m}_{\mathrm{TL}}(u,x) = \widehat{m}_{\mathrm{T}}(u,x) + \widehat{b}(u,x), \qquad \widehat{b}(u,x) = \text{smooth of } \big\{(u^{(1)}_{t_1},x^{(1)}_{t_1},Y^{(1)}_{t_1} -\widehat{m}_{\mathrm{T}}(u^{(1)}_{t_1},x^{(1)}_{t_1}))\big\}.
\]
All estimators use the product Epanechnikov kernel
\[
    K(u,x) = \frac{3}{4}(1-u^2)\mathbf{1}_{(|u|\le1)}\, \frac{3}{4}(1-x^2)\mathbf{1}_{(|x|\le1)},
\]
and bandwidths $(h_u,h_x)$ (and $(h_u^{\mathrm{TL}},h_x^{\mathrm{TL}})$ for bias smoothing) are chosen by ten-fold least-squares cross-validation within each domain. 
All models are fitted on the entire target dataset, and the bias-correction step relies only on target residuals and source predictions, preventing information leakage.

Each estimated surface is bilinearly interpolated onto a common $N\times N$ grid of $[0,1]^2$ with $N=50$. 
At grid points $(u_i,x_j)$, we compute the pointwise squared error
\[
    E^{(m)}_{ij} = [\widehat{m}(u_i,x_j)-m^{(0)}(u_i,x_j)]^2,
\]
and summarize it by the median across all grid locations,
\[
    \mathrm{L2}^{(m)}_{\text{iteration}} = \operatorname{median}_{i,j} E^{(m)}_{ij}.
\]
Using the median reduces the influence of boundary regions where kernel weights vanish or numerical instability may occur. 
Each configuration $(\gamma,\text{bias type})$ is replicated fifty times, and the resulting median errors are averaged across replications. 
Line plots of the mean median error against $\gamma$, together with boxplots displaying the empirical variation, summarize the performance of the four estimators.

The simulation design is chosen to mirror the theoretical setting of Section~4 as closely as possible. The locally stationary covariate process ensures that the uniform bias and variance expansions established in the theory remain valid over the entire domain, while the family of discrepancy surfaces \(\{b_\gamma\}\) introduces a controlled curvature scale that corresponds naturally to the bandwidth order \(O(h_u^{2}+h_x^{2})\). When \(|\gamma|\) is small, the discrepancy varies smoothly within each local neighborhood, allowing the transfer learning estimators to draw effectively on the larger source sample and achieve substantial variance reduction. As \(|\gamma|\) grows, the curvature of \(b_\gamma\) eventually exceeds the resolution of the local smoother, reducing the accuracy of the bias correction and causing the performance of the transfer estimators to approach that of the target-only estimators. This design therefore provides a clear finite-sample illustration of the transition between beneficial and negligible transfer anticipated by the asymptotic theory.

\subsection{Simulation Results}

This section reports the finite-sample performance of the four estimators described in Section~5.1 under the three bias families: quadratic, cubic, and exponential. 
Each plot shows the mean of the median grid errors across fifty replications as a function of the domain discrepancy parameter $\gamma$. 
The corresponding boxplots that illustrate the distribution of the errors are provided in Appendix~\ref{app:simulation_boxplots}.

\begin{figure}[!ht]
    \centering
    \includegraphics[width=0.85\textwidth]{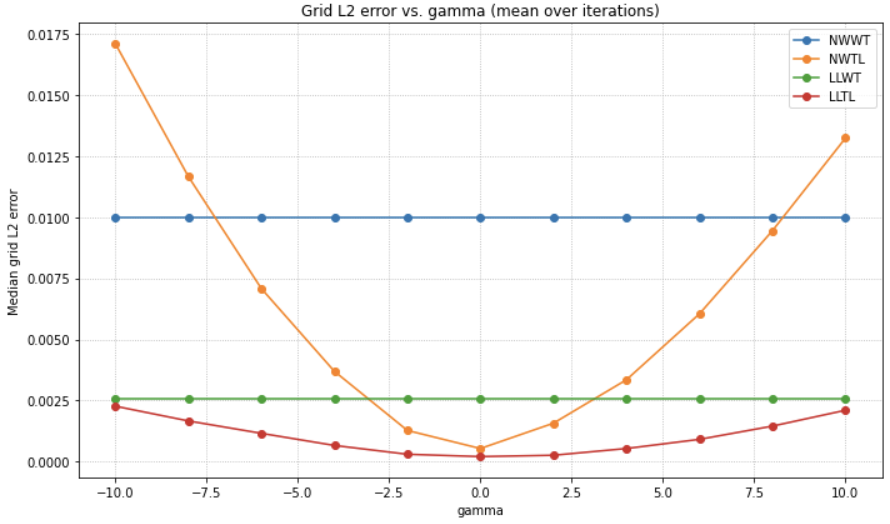}
    \caption{Mean of median grid errors versus $\gamma$ for the quadratic bias family. 
    Each curve represents the average over fifty replications for the Nadaraya–Watson and locally linear estimators trained on the target data and their transfer-learning counterparts.}
    \label{fig:quad_medianerror}
\end{figure}

\textit{Quadratic bias.}
Figure~\ref{fig:quad_medianerror} displays the results for the quadratic bias family. 
The two estimators fitted only on the target data, namely the Nadaraya–Watson estimator and the locally linear estimator, form almost horizontal curves since their performance does not depend on $\gamma$. 
The locally linear estimator remains consistently below the Nadaraya–Watson estimator, which is consistent with the smaller asymptotic bias shown in Section~3. 
Both transfer-learning estimators exhibit a clear ``V'' shape centered at $\gamma=0$. 
When $\gamma$ is small, the discrepancy between source and target domains is minor, and transfer learning effectively reduces variance by borrowing information from the large source sample. 
As the absolute value of $\gamma$ increases, the discrepancy surface becomes more curved, leading to greater bias and higher estimation error. 
At extreme values of $\gamma$, both transfer-learning curves approach the target-only benchmarks, illustrating the theoretical prediction that transfer benefits vanish when the domains diverge substantially. 
In the quadratic case, the locally linear transfer-learning estimator achieves the lowest overall error and remains below all other curves throughout the investigated range of $\gamma$.

\begin{figure}[!ht]
    \centering
    \includegraphics[width=0.85\textwidth]{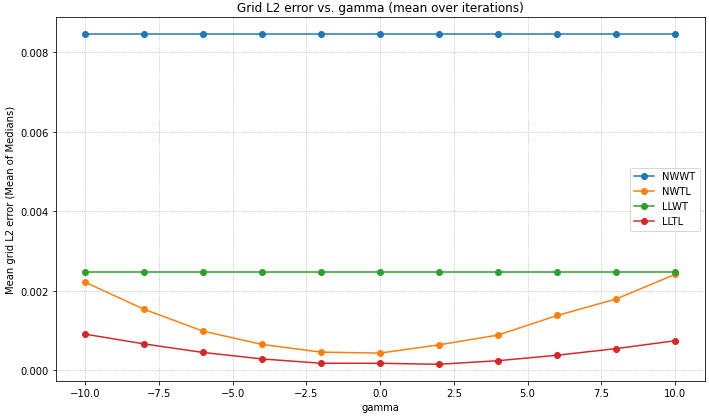}
    \caption{Mean of median grid errors versus $\gamma$ for the cubic bias family. 
    The curves illustrate the bias–variance transition predicted by theory: transfer-learning estimators achieve lower error for small $|\gamma|$ but converge to the target-only estimators as the discrepancy curvature increases.}
    \label{fig:cubic_medianerror}
\end{figure}

\textit{Cubic bias.}
Figure~\ref{fig:cubic_medianerror} reports the results for the cubic bias family. 
The same overall pattern is observed as in the quadratic case. 
The Nadaraya–Watson and locally linear estimators trained only on the target domain remain horizontal, with the locally linear version consistently exhibiting lower error. 
The two transfer-learning estimators again display a pronounced ``V'' shape, confirming the symmetric deterioration in performance as the absolute value of $\gamma$ increases. 
For moderate discrepancies, both transfer-learning estimators outperform their target-only counterparts, verifying that the theoretical smoothness conditions of Theorem~4.1 yield tangible finite-sample benefits. 
At large curvature levels, the estimation errors of the transfer-learning methods approach those of the target-only estimators, demonstrating the predicted phase transition where the bias term dominates. 
The locally linear transfer-learning estimator maintains the lowest error throughout, underscoring its robustness to curvature misspecification.

\begin{figure}[!ht]
    \centering
    \includegraphics[width=0.85\textwidth]{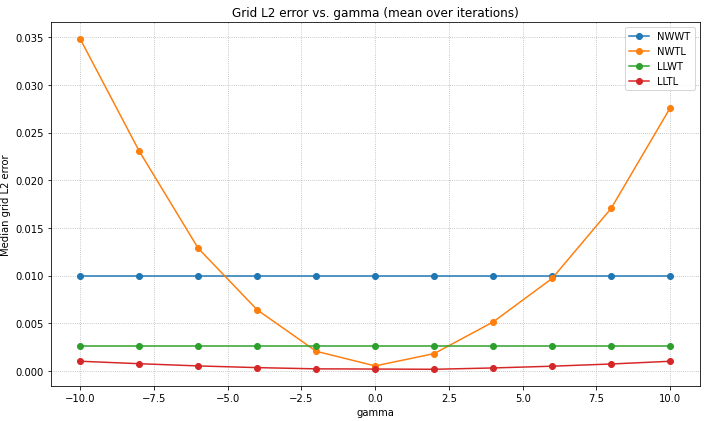}
    \caption{Mean of median grid errors versus $\gamma$ for the exponential bias family. 
    The transfer-learning estimators display the same ``V''-shaped pattern observed in the polynomial cases, and the locally linear transfer-learning estimator remains superior across the entire range of $\gamma$.}
    \label{fig:exp_medianerror}
\end{figure}

\textit{Exponential bias.}
Figure~\ref{fig:exp_medianerror} presents the results for the exponential bias family, which produces slightly sharper curvature in $b_\gamma(u,x)$ compared with the polynomial cases. 
The general behavior remains consistent with theoretical expectations. 
The two target-only estimators remain nearly flat across $\gamma$, with the locally linear version yielding smaller errors than the Nadaraya–Watson version. 
The transfer-learning estimators again form ``V''-shaped curves, symmetric around $\gamma=0$, and their error increases with the magnitude of $\gamma$. 
The Nadaraya–Watson transfer-learning estimator converges rapidly to the target-only level as the discrepancy grows, whereas the locally linear transfer-learning estimator remains below both target-only curves over the entire range $\gamma\in[-10,10]$. 
This persistence of improvement confirms that locally linear bias correction continues to be effective even under moderately large curvature, illustrating the advantage of the local linear adjustment in maintaining bias control.

Taken together, these findings demonstrate clear agreement with the theoretical framework in Section~4. 
The “V’’ shapes observed across all bias families reflect the rate decomposition in Theorem~\ref{Thm4.2}. 
For small $|\gamma|$, the curvature of $b_\gamma(u,x)$ is mild and the transfer variance term of order $(T_1 h_{\tl}^{d+1})^{-1/2}$ dominates the curvature bias, leading to strictly lower error than the target-only estimators. 
As $|\gamma|$ increases, the curvature-driven component $\eta_{2,b}/(T_0^{r}h_{\tl}^{\,d-1})$ grows and eventually offsets the variance gain, producing the symmetric increase in estimation error. 
For large curvature, the discrepancy bias dominates and the transfer-learning estimators converge to the target-only benchmarks, precisely matching the theoretical phase transition.
The consistently superior performance of the locally linear transfer-learning estimator illustrates the sharper bias control predicted by the expansions in Section~3 and Theorem~\ref{Thm4.2}. 
Because the locally linear correction eliminates the leading bias in $u$ and $x$, the locally stationary remainder enters at a smaller order than for the Nadaraya-Watson version, allowing it to remain effective even under moderate curvature. 
These finite-sample patterns therefore provide direct empirical confirmation of the theoretical bias-variance tradeoff underlying the transfer-learning procedure.

\section{Empirical Application}\label{sec-empirical}
\subsection{Application Setting}

We demonstrate the proposed methodology through an empirical study of international fuel prices and crude oil prices. The target domain, indexed by the superscript $0$, corresponds to the United States, and the source domain, indexed by the superscript $1$, corresponds to Korea. In each specification, $(U^{(\ell)}, X^{(\ell)}, Y^{(\ell)})$ for $\ell \in \{0,1\}$ denotes the rescaled time, covariate, and response variables, respectively. The response $Y^{(\ell)}$ represents the retail fuel price, and the covariate $X^{(\ell)}$ represents either the crude oil price or the price of a related fuel product. The primary goal is to estimate the time-varying regression function
\[
m^{(\ell)}(u,x) = \mathbb{E}\bigl[Y^{(\ell)} \mid U^{(\ell)}=u,\, X^{(\ell)}=x\bigr],
\]
and to compare three estimators on the target domain: the target-only estimator $\hat m_0$, the transfer learning estimator $\hat m^{\mathrm{TL}}$, and the pooled estimator $\hat m^{\mathrm{Pool}}$.

The sample period spans 2010–2024. Daily Korean retail fuel prices for regular gasoline, diesel, premium gasoline, and kerosene are obtained from the Opinet database, while weekly U.S.\ retail prices for regular gasoline, diesel, and premium gasoline are collected from the U.S.\ Energy Information Administration (EIA). Daily spot prices for West Texas Intermediate (WTI) and Brent crude oil are included as covariates. All series are cleaned by removing nonnumeric entries and filling short missing segments using the local mean of neighboring observations. Each series is then converted into a triple $(U^{(\ell)}, X^{(\ell)}, Y^{(\ell)})$ defined on a normalized time scale $[0,1]$.

Because the Korean data are daily and the U.S.\ data are weekly, the two time indices are first harmonized. The source time index $U^{(1)}$ is constructed by ordering all daily observations by date and assigning $U_{t_1,T_1}^{(1)} = t_1/T_1$ for $t_1 = 1,\dots,T_1$, where $T_1$ is the total number of source observations. The target time index $U^{(0)}$ is then obtained by interpolating the source time map onto the weekly U.S.\ observation dates. This mapping embeds both domains into the same rescaled interval $[0,1]$ and ensures that comparable calendar periods correspond to similar rescaled time points.

The construction of the covariate $X^{(\ell)}$ depends on the specification. When both domains share a common exogenous variable, such as crude oil prices, the lagged daily crude oil price is used as $X_{t_1,T_1}^{(1)}$ in the source domain to ensure temporal precedence. For the target domain, $X_{t_0,T_0}^{(0)}$ is defined as the average of daily crude prices over the weekly interval $[\tau_{t_0-1}, \tau_{t_0})$, where $\tau_{t_0}$ is the $t_0$th weekly observation date. This yields a coherent pair $(X^{(0)}, X^{(1)})$ representing the same underlying process at different sampling frequencies. This specification is applied to diesel prices with WTI and Brent crude oil as covariates.

In alternative specifications, $X^{(\ell)}$ represents a related fuel price within each country. For example, diesel prices $(Y^{(0)}, Y^{(1)})$ are regressed on regular gasoline prices $(X^{(0)}, X^{(1)})$, or regular gasoline prices are regressed on diesel prices. Another case considers premium gasoline prices as responses with regular gasoline prices as covariates. In all such designs, the covariate precedes the response by one period within each domain, preserving a causal ordering.

For comparability, each covariate $X^{(\ell)}$ is linearly rescaled to approximately $[0,1]$ within its domain. For the transfer learning estimator, the source response $Y^{(1)}$ remains on its original scale. For the pooled estimator, $Y^{(1)}$ is standardized to match the mean and variance of the target response $Y^{(0)}$, which aligns the global level while maintaining the relative dynamics of the source domain.

The empirical analysis is based on the constructed triples $(U^{(\ell)}, X^{(\ell)}, Y^{(\ell)})$. The target sample is partitioned deterministically into a training and test set by allocating every fourth observation to the test set, ensuring an even temporal spread of evaluation points. The source sample is fully used for the transfer learning estimator and, after removing time points with duplicated $U$ values, is combined with the target training sample for the pooled estimator.

The three estimators considered in the analysis differ in how they incorporate information from the source domain. The target-only estimator $\hat m_0$ uses only the target training data and serves as a benchmark. The transfer learning estimator $\hat m^{\mathrm{TL}}$ first estimates the regression function $\hat m^{(1)}$ in the source domain and then corrects its systematic bias on the target domain via a nonparametric bias function $\hat b(u,x)$, yielding
\[
\hat m^{\mathrm{TL}}(u,x) = \hat m^{(1)}(u,x) + \hat b(u,x).
\]
This procedure leverages the structural similarity between the two domains while allowing local adjustments for domain-specific differences. The pooled estimator $\hat m^{\mathrm{Pool}}$ combines the target and filtered source samples directly into a unified dataset, thereby exploiting cross-domain information through joint kernel smoothing. The three estimators thus represent distinct strategies for transferring information: no transfer, bias-corrected transfer, and full pooling.

For each specification, the regression surface $m(u,x)$ is estimated on a uniform grid over $[0,1]^2$. Both the Nadaraya–Watson estimator $\hat m_{\mathrm{NW}}(u,x)$ and the locally linear estimator $\hat m_{\mathrm{LL}}(u,x)$ are implemented, with bandwidths $(h_u,h_x)$ selected by $k$-fold least-squares cross validation. The estimators $\hat m_0$, $\hat m^{\mathrm{TL}}$, and $\hat m^{\mathrm{Pool}}$ are then evaluated on the held-out target test set. Predictive performance is assessed by empirical $L_2$ and $L_\infty$ losses, and the fitted surfaces are visualized in both two and three dimensions.

\subsection{Results}

We report empirical findings from four transfer learning applications linking U.S.\ and Korean fuel prices. Each experiment compares the target--only estimator \(\hat m_0\), the transfer learning estimator \(\hat m^{\mathrm{TL}}\), and the pooled estimator \(\hat m^{\mathrm{Pool}}\), implemented using both the Nadaraya--Watson (NW) and locally linear (LL) methods. The settings include diesel prices with WTI and Brent crude as covariates, diesel prices regressed on regular gasoline, and gasoline or premium gasoline prices regressed on related fuels. Throughout the analysis, the United States is treated as the target domain and Korea as the source domain.

The results align closely with the theoretical predictions of Section~4. The transfer learning estimator \(\hat m^{\mathrm{TL}}\) consistently outperforms the target--only estimator \(\hat m_0\), whereas the pooled estimator \(\hat m^{\mathrm{Pool}}\) often performs worse. This pattern reflects the bias decomposition in Theorem~4.1, in which the source--target discrepancy is captured by a smooth bias function \(b(u,x)\) that varies gradually over time. When this bias is sufficiently smooth—its derivatives uniformly bounded and its variation smaller than the local bandwidth—the correction effectively removes cross--domain differences without inflating variance. The pooled estimator, by implicitly assuming \(b(u,x)\equiv 0\), suffers from negative transfer whenever structural heterogeneity is present. The LL estimator further strengthens performance by reducing boundary bias and adapting to local curvature, which leads to smaller \(L^2\) and \(L^\infty\) errors relative to NW across all settings.

\paragraph*{Diesel with WTI crude as covariate.}
Table~\ref{tab:diesel_WTI} summarizes the results when diesel prices are regressed on WTI crude oil prices. The transfer learning estimator reduces the \(L^2\) error by approximately 23\% under NW and 21\% under LL relative to the target--only benchmark. The pooled estimator performs markedly worse, indicating that the diesel--WTI relationship differs across the two economies and cannot be treated as homogeneous. Figures~\ref{fig:diesel_WTI_nw} and \ref{fig:diesel_WTI_ll} show that \(\hat m^{\mathrm{TL}}\) effectively adjusts the Korean structure toward the U.S.\ trend, consistent with the behavior predicted by Theorem~4.2.

\begin{table}[H]
\centering
\caption{Diesel fuel prices with WTI crude as covariate. Target: U.S.\ weekly diesel. Source: Korean daily diesel.}
\label{tab:diesel_WTI}
\begin{tabular}{lcccc}
\toprule
& \multicolumn{2}{c}{Nadaraya-Watson (NW)} & \multicolumn{2}{c}{Locally Linear (LL)} \\
\cmidrule(lr){2-3} \cmidrule(lr){4-5}
& $L^2$ & $L^\infty$ & $L^2$ & $L^\infty$ \\
\midrule
Baseline & 0.180 & 0.727 & 0.167 & 0.851\\
Transfer Learning & 0.139 & 0.533 & 0.132 & 0.476\\
Pooled & 0.282 & 0.915 & 0.229 & 1.194\\
\bottomrule
\end{tabular}
\end{table}

\paragraph*{Diesel with Brent crude as covariate.}
Table~\ref{tab:diesel_brent} shows similar patterns when Brent crude is used as the covariate. The transfer learning estimator again achieves the lowest prediction errors, reducing the \(L^2\) metric from 0.169 to 0.135 under NW and from 0.124 to 0.103 under LL. The improvement suggests that the cross-domain bias in the diesel-Brent relation is smoother than in the diesel-WTI case, consistent with the theoretical prediction that smoother bias functions yield sharper gains.

\begin{table}[H]
\centering
\caption{Diesel fuel prices with Brent crude as covariate. Target: U.S.\ weekly diesel. Source: Korean daily diesel.}
\label{tab:diesel_brent}
\begin{tabular}{lcccc}
\toprule
& \multicolumn{2}{c}{Nadaraya-Watson (NW)} & \multicolumn{2}{c}{Locally Linear (LL)} \\
\cmidrule(lr){2-3} \cmidrule(lr){4-5}
& $L^2$ & $L^\infty$ & $L^2$ & $L^\infty$ \\
\midrule
Baseline & 0.169 & 0.722 & 0.124 & 0.685 \\
Transfer Learning & 0.135 & 0.522 & 0.103 & 0.516 \\
Pooled & 0.306 & 1.192 & 0.262 & 0.922\\
\bottomrule
\end{tabular}
\end{table}

\paragraph*{Diesel with gasoline as covariate.}
Table~\ref{tab:diesel_gasoline} reports the results when diesel prices are regressed on regular gasoline. The transfer learning estimator continues to outperform the target-only baseline, indicating that the structural comovement of related fuel types can be effectively utilized through smooth bias correction. The pooled estimator performs worse due to temporal and scale mismatches between the two domains, illustrating the negative transfer that arises when \(b(u,x)\) varies more rapidly than the local bandwidth.

\begin{table}[H]
\centering
\caption{Diesel fuel prices with gasoline as covariate. Target: U.S.\ weekly diesel. Source: Korean daily diesel.}
\label{tab:diesel_gasoline}
\begin{tabular}{lcccc}
\toprule
& \multicolumn{2}{c}{Nadaraya-Watson (NW)} & \multicolumn{2}{c}{Locally Linear (LL)} \\
\cmidrule(lr){2-3} \cmidrule(lr){4-5}
& $L^2$ & $L^\infty$ & $L^2$ & $L^\infty$ \\
\midrule
Baseline & 0.176 & 0.699 & 0.197 & 0.894 \\
Transfer Learning & 0.132 & 0.449 & 0.145 & 0.539 \\
Pooled & 0.247 & 0.723 & 0.173 & 0.679\\
\bottomrule
\end{tabular}
\end{table}

\paragraph*{Gasoline with diesel as covariate.}
The results in Table~\ref{tab:gasoline_diesel} show a smaller gap among the estimators. The improvement from transfer learning remains positive but is modest, suggesting that the bias between U.S. and Korean gasoline--diesel relations is less smooth or exhibits partial asymmetry. This is consistent with the theoretical insight that higher local variability in the bias function limits the potential efficiency gains, though \(\hat m^{\mathrm{TL}}\) still improves upon \(\hat m_0\).

\begin{table}[H]
\centering
\caption{Regular gasoline prices with diesel as covariate. Target: U.S.\ weekly gasoline. Source: Korean daily gasoline.}
\label{tab:gasoline_diesel}
\begin{tabular}{lcccc}
\toprule
& \multicolumn{2}{c}{Nadaraya-Watson (NW)} & \multicolumn{2}{c}{Locally Linear (LL)} \\
\cmidrule(lr){2-3} \cmidrule(lr){4-5}
& $L^2$ & $L^\infty$ & $L^2$ & $L^\infty$ \\
\midrule
Baseline & 0.188 & 0.662 & 0.168 & 0.637 \\
Transfer Learning & 0.174 & 0.547 & 0.167 & 0.606 \\
Pooled & 0.208 & 0.618 & 0.184 & 0.608\\
\bottomrule
\end{tabular}
\end{table}

\paragraph*{Premium gasoline with regular as covariate.}
Table~\ref{tab:premium_regular} presents the strongest transfer effect. The LL transfer estimator achieves the smallest prediction errors, reducing the \(L^2\) error by nearly half relative to the target--only approach. This setting features an almost constant bias function, as the premium--regular spread evolves smoothly across countries. As a result, the theoretical efficiency bound for \(\hat m^{\mathrm{TL}}\) is nearly achieved, while the pooled estimator suffers from over-smoothing induced by redundant heterogeneity.

\begin{table}[H]
\centering
\caption{Premium gasoline prices with regular gasoline as covariate. Target: U.S.\ weekly premium. Source: Korean daily premium.}
\label{tab:premium_regular}
\begin{tabular}{lcccc}
\toprule
& \multicolumn{2}{c}{Nadaraya-Watson (NW)} & \multicolumn{2}{c}{Locally Linear (LL)} \\
\cmidrule(lr){2-3} \cmidrule(lr){4-5}
& $L^2$ & $L^\infty$ & $L^2$ & $L^\infty$ \\
\midrule
Baseline & 0.178 & 0.565 & 0.050 & 0.203 \\
Transfer Learning & 0.102 & 0.434 & 0.052 & 0.186 \\
Pooled & 0.309 & 0.785 & 0.301 & 0.789\\
\bottomrule
\end{tabular}
\end{table}

Overall, the empirical evidence reinforces the theoretical conclusions of Section~4. Transfer learning yields consistent gains when the inter--domain bias \(b(u,x)\) is locally smooth, whereas naive pooling increases bias in the presence of persistent heterogeneity. The locally linear transfer estimator \(\hat m^{\mathrm{TL}}_{\mathrm{LL}}\) performs most reliably across settings, validating the theoretical rate advantage and illustrating the practical value of smooth bias correction under local stationarity.

\subsection{Economic Relevance}
The empirical patterns have clear economic implications that complement their statistical interpretation. Retail fuel prices reflect the pass-through of global crude oil costs together with refining margins, distribution, and taxes. According to the U.S.\ Energy Information Administration (EIA), crude oil accounts for roughly one half of the retail gasoline price and more than forty percent of the retail diesel price, with the remaining share determined by refining, distribution, and taxation \cite{EIA_FuelBreakdown}. This composition explains why crude benchmarks, particularly Brent, play a central role in cross--country fuel price dynamics. The strong improvement observed when Brent is used as a covariate for U.S.\ diesel prices reflects the dominant contribution of global crude costs to retail outcomes. The bias correction in the transfer learning estimator aligns market-specific deviations in price pass--through, enabling global information to enhance domestic prediction accuracy.

The experiments linking gasoline and diesel prices underscore the economic interdependence among refined products. Because gasoline and diesel are joint outputs of the same refining process, their prices move together through refinery yield constraints and crack--spread relationships. The finding that U.S.\ diesel prices benefit substantially from Korean gasoline information, while the reverse linkage yields smaller gains, mirrors the asymmetric structure of global refinery operations. Diesel demand typically governs refinery utilization, with gasoline adjusting accordingly. The transfer learning estimator reproduces this asymmetry by adapting cross--domain information to the U.S.\ context, illustrating how bias correction can recover directional relationships between jointly produced goods.

The application to premium and regular gasoline highlights the role of substitution on both the supply and demand sides. The two grades share similar input costs and closely aligned consumption patterns. The sizable reduction in prediction errors under transfer learning, particularly with the locally linear estimator, indicates that auxiliary information from one grade nearly captures the dynamics of the other. This demonstrates how transfer learning can exploit information from near substitutes when one product is observed more frequently or with greater precision.

Taken together, the results show that transfer learning is both statistically effective and economically meaningful. The gains arise precisely where economic reasoning suggests cross--market information should be informative: when global crude price transmission dominates, when products share production technologies, and when substitution patterns generate strong co--movements. By formalizing how local bias can be corrected through smooth cross--market adaptation, the framework provides a data--driven tool for integrating heterogeneous international information, improving forecasting and market monitoring in settings where cross--country data disparities are substantial.

\section{Discussion}\label{sec-diss}
\subsection{Discussions on Multivariate Locally Linear Estimation}
The multivariate locally linear estimator provides a significant theoretical and practical advancement over the local-constant Nadaraya-Watson approach for non-stationary time series. By expressing the estimation through a local first-order Taylor expansion, the procedure moves beyond simple averaging. The estimator is derived as the solution to a locally weighted least-squares problem that simultaneously fits the regression surface $m(u,x)$ and its scaled partial derivatives with respect to rescaled time, $\nabla_0 m$, and the covariates, $\nabla_j m$. This mechanism is crucial for its improved performance, particularly in mitigating boundary effects. The Nadaraya-Watson estimator suffers from first-order boundary bias, arising from the asymmetry of the kernel support near the domain edges. In contrast, the locally linear estimator's explicit fitting of local slopes corrects for this asymmetry, yielding a classical smoothing bias term of order $O_p(h^2)$ that holds uniformly across the entire compact domain, $[0,1] \times S$.

This improvement is most evident in the uniform error decomposition established in the theoretical results. The convergence rate for the complete estimator vector, $\hat{\mathfrak{m}}(u,x) - \mathfrak{m}(u,x)$, is shown to be composed of three distinct parts:
$$
O_p\left(\sqrt{\frac{\log T}{Th^{d+1}}} + \frac{1}{T^r h^{d-1}} + h^2\right)
$$
The first component, the stochastic fluctuation term of order $O_p\left(\sqrt{\frac{\log T}{Th^{d+1}}}\right)$, represents the standard optimal rate for a nonparametric problem of this dimensionality and is shared by the Nadaraya-Watson estimator. The second component is the aforementioned uniform second-order smoothing bias, $O_p(h^2)$.

The main theoretical contribution is the refinement of the second component, namely the remainder term that arises from the local-stationarity approximation. The analysis shows this nonstationarity remainder to be of order $O_p\left(\frac{1}{T^r h^{d-1}}\right)$, where $r = \min\{\rho, 1\}$ quantifies the quality of the local stationary approximation. This rate is a direct and notable improvement over the corresponding $O_p\left(\frac{1}{T^r h^d}\right)$ remainder associated with the local-constant framework. The reduction in the remainder's dependence on dimensionality, from $h^{-d}$ to $h^{-(d-1)}$, demonstrates that the locally linear estimator's capacity to model local drift makes it inherently less sensitive to violations of strict stationarity. This structural advantage stabilizes the local design matrix and provides a more robust and accurate approximation of the evolving regression surface, especially in higher dimensions or near the temporal boundaries where traditional kernel smoothers are least reliable.

\subsection{Discussions on Transfer Learning Framework} \label{subsec:TL-discussion}
\label{sec:oracle-neglig-dominance}

We now derive the oracle rate of the locally linear transfer estimator by optimizing its bias–variance–remainder trade-off, establish conditions under which the auxiliary (source) component becomes asymptotically negligible, and then compare the resulting rate with the target-only benchmarks to determine when transfer learning provides a strict improvement.

\vspace{0.5em}
\noindent
Starting from \eqref{eqthm4.2}, we define the oracle target rate as
\begin{equation}\label{eq:LL-TL-oracle-again}
    A^\star(T_0;\eta_{2,b}) :=\inf_{h_{\tl}>0} \left\{ \sqrt{\frac{\log T_0}{T_0 h_{\tl}^{\,d+1}}} +\frac{\eta_{2,b}}{T_0^{r}h_{\tl}^{\,d-1}} +\eta_{2,b}h_{\tl}^{2} \right\},
\end{equation}

and determine its asymptotic behavior by balancing the three terms.

\emph{Case 1 ($r \ge \frac{d+1}{d+5}$).}
When the process is strongly locally stationary, the nonstationary remainder is dominated by the bias–variance trade-off.  
Balancing the first and third terms in~\eqref{eq:LL-TL-oracle-again} gives
\begin{equation*}
    h_{\tl}\asymp \eta_{2,b}^{-\frac{2}{d+5}}(\log T_0)^{1/(d+5)}T_0^{-1/(d+5)},
\end{equation*}
which leads to
\begin{equation*}
    \sup_{u\in[0,1],\,x\in S}\Big|\hat{m}^{\tl}(u,x)-m^{(0)}(u,x)\Big|
    = O_P\!\Big(\eta_{2,b}^{\frac{d+1}{d+5}}(\log T_0)^{\frac{2}{d+5}}T_0^{-\frac{2}{d+5}}\Big).
\end{equation*}
This rate matches the stationary-optimal locally linear order, while curvature alignment through $\eta_{2,b}<1$ provides a multiplicative gain that reduces both variance and bias constants.

\emph{Case 2 ($r < \tfrac{d+1}{d+5}$ and $(\log T_0)^{1/2}T_0^{\frac{r(d+5)-(d+1)}{2(d+1)}}<\eta_{2,b}\le 1$).}
When nonstationarity dominates, the reduced remainder term $T_0^{-r}h_{\tl}^{-(d-1)}$ determines the balance.  
Equating the remainder and bias in~\eqref{eq:LL-TL-oracle-again} gives
\begin{equation*}
    h_{\tl}\asymp \eta_{2,b}^{1/(d+2)}T_0^{-\,r/(d+2)},
\end{equation*}
and hence
\begin{equation*}
    \sup_{u\in[0,1],\,x\in S}\Big|\hat{m}^{\tl}(u,x)-m^{(0)}(u,x)\Big|
    =O_P\!\big(\eta_{2,b}\,T_0^{-\,2r/(d+2)}\big).
\end{equation*}
This order is strictly faster than the target-only locally linear rate $O_P(T_0^{-2r/(d+1)})$ whenever $\eta_{2,b}<1$, as the transferred curvature information partially corrects the target’s nonstationary remainder.

\emph{Case 3 ($r<\tfrac{d+1}{d+5}$ and $\eta_{2,b}<(\log T_0)^{1/2}T_0^{\frac{r(d+5)-(d+1)}{2(d+1)}}$).}
When curvature alignment is sufficiently strong, the remainder term becomes negligible, and the variance–bias trade-off once again governs the rate.  
Balancing the first and third terms yields
\begin{equation*}
    h_{\tl}\asymp \eta_{2,b}^{-\frac{2}{d+5}}(\log T_0)^{1/(d+5)}T_0^{-1/(d+5)},
\end{equation*}
which implies
\begin{equation*}
    \sup_{u\in[0,1],\,x\in S}\Big|\hat{m}^{\tl}(u,x)-m^{(0)}(u,x)\Big|
    =O_P\!\Big(\eta_{2,b}^{\frac{d+1}{d+5}}(\log T_0)^{\frac{2}{d+5}}T_0^{-\frac{2}{d+5}}\Big).
\end{equation*}
This again achieves the stationary-optimal order with an improved constant determined by $\eta_{2,b}$.

\vspace{0.5em}
\noindent
With the oracle rate characterized, we now establish when the auxiliary (source) term in the full estimator becomes negligible.  
The source component
\[
    \sqrt{\frac{\log T_1}{T_1 h_1^{\,d+1}}} +\frac{1}{T_1^{r}h_1^{\,d-1}} +h_1^{2}
\]
is asymptotically dominated by the target-side rate if
\begin{equation}\label{eq:source-negl-final}
    \frac{T_1 h_1^{d+1}}{\log T_1}\gg (A^\star)^{-2}, \qquad
    T_1^{r}h_1^{d-1}\gg (A^\star)^{-1}, \qquad
    h_1\ll (A^\star)^{1/2},\quad A^\star:=A^\star(T_0;\eta_{2,b}),
\end{equation}
where $A^\star(T_0;\eta_{2,b})$ is the rate determined in the three cases above.  
The feasible range for $h_1$ that satisfies these conditions is
\begin{equation*}
    \max\!\left\{ \Big(\tfrac{\log T_1}{(A^\star)^2 T_1}\Big)^{\!\frac{1}{d+1}}, \ \Big(\tfrac{1}{A^\star T_1^{r}}\Big)^{\!\frac{1}{d-1}} \right\} \ll h_1 \ll (A^\star)^{1/2},
\end{equation*}
and the interval is nonempty provided that
\begin{equation*}
    \frac{T_1}{\log T_1}\gg (A^\star)^{-\frac{d+5}{2}}, \qquad T_1\gg (A^\star)^{-\frac{d+1}{2r}}.
\end{equation*}
Under these conditions, the source term is asymptotically negligible, and the convergence rate of the transfer estimator is fully determined by the target-side rate $A^\star(T_0;\eta_{2,b})$.

\emph{Comparison with target-only estimators.}
The target-only locally linear estimator satisfies
\[
\sup_{u,x}\big|\hat m^{(0)}(u,x)-m^{(0)}(u,x)\big|=
\begin{cases}
O_P\!\big(T_0^{-\,2r/(d+1)}\big), & r<\tfrac{d+1}{d+5},\\[2pt]
O_P\!\big((\log T_0)^{2/(d+5)}T_0^{-\,2/(d+5)}\big), & r\ge\tfrac{d+1}{d+5},
\end{cases}
\]
and the Nadaraya-Watson estimator attains $O_P((\log T_0/T_0)^{2/(d+3)})$.  
In Case~1, the transfer estimator achieves the same order as the stationary-optimal locally linear estimator but with a smaller constant $\eta_{2,b}^{(d+1)/(d+5)}$.  
In Case~2, it is strictly faster than the target-only locally linear estimator  with a smaller constant $\eta_{2,b}$.  
In Case~3, the transfer estimator achieves the same order as the stationary-optimal locally linear estimator because $- \frac{2}{d+5}<-\,\frac{2r}{d+1}$ when $r<\tfrac{d+1}{d+5}$.
If $\eta_{2,b}$ are sufficiently small, the rate becomes substantially faster than that of the standard locally linear estimator, establishing when and how transfer learning achieves superior asymptotic performance. 
If both $\eta_{1,b}$ and $\eta_{2,b}$ are sufficiently small so that the discrepancy varies smoothly within the effective bandwidth neighborhood, then an analogous phenomenon appears in the Nadaraya Watson setting as well. In that case the additional curvature induced by local stationarity remains dominated by the smoothing scale, and the bias correction retains its accuracy, which enables the transfer learned estimator to outperform the target only Nadaraya Watson estimator in asymptotic order. This identifies a parallel regime in which even the simpler local constant methods benefit from information borrowing and achieve strictly improved convergence behavior.

\subsection{Concluding Remarks}

This study develops a framework that connects local stationarity, nonparametric kernel regression, and transfer learning.  
By extending the classical Nadaraya–Watson estimator to a multivariate locally linear form, we establish a nonparametric approach that remains valid under temporal dependence and smoothly time-varying structures. The proposed bias-corrected transfer procedure enables principled information sharing across related domains, improving estimation in data-scarce environments. Together, these developments provide both theoretical guarantees and empirical evidence that transfer learning can be rigorously justified for dependent, nonstationary time series when the discrepancy between domains is sufficiently smooth.  
The results delineate clear boundaries of effectiveness: as the discrepancy exceeds the smoothing scale, transfer offers no additional benefit, thereby identifying when knowledge transfer is both safe and efficient.

Despite these contributions, several limitations remain.  
The current analysis assumes comparable design densities between source and target domains, excluding explicit covariate-shift or heteroskedastic structures.  
Relaxing these assumptions would allow the framework to accommodate more general forms of distributional heterogeneity.  
The bandwidth selection procedure relies on empirical cross-validation, which is practical but lacks a fully theoretical foundation under dependence.  
Developing adaptive or plug-in criteria for bandwidth and local likelihood selection could strengthen both theoretical rigor and empirical robustness.  
Moreover, the present framework focuses on mean-function transfer; extending it to conditional variance, quantile, or spectral transfer would broaden its applicability to richer forms of dynamic dependence.  

Future work may refine the multivariate locally linear design by incorporating higher-order local polynomials or semiparametric components that better capture nonlinear temporal curvature.
In high-dimensional or sparse settings, integrating local stationarity with regularization techniques such as functional Lasso or additive decomposition could yield scalable alternatives.  
Extending the framework to multi-source or hierarchical transfer scenarios, or to frequency-domain representations, would further enhance its capacity to integrate heterogeneous time-varying information.  
These directions represent natural continuations of the present study and hold promise for building a more comprehensive statistical theory of transfer learning under dependence and nonstationarity.  
Ultimately, this work underscores a broader principle for studying evolving systems: effective adaptation across domains depends on balancing local approximation accuracy with global structural alignment. This insight provides a foundation for future developments in nonparametric inference for complex, heterogeneous, and time-dependent data environments.

\bibliography{ref}

\newpage
\appendix
\section{Appendix A}

\subsection{Proof of \hyperref[Thm3.1]{Theorem 3.1}}\label{prf3.1}
\begin{proof}

We write, 
\begin{equation*}
    \begin{aligned}
        \hat{\Psi}(u,x) &= \frac 1 {T} \mbf{D}^{\top} \mbf{W} \mbf{R}_{T} \\:&= 
        \begin{pmatrix}
            \frac 1 {T} \sum_{t=1}^T K_h \left(u - \frac{t}{T}\right)  \prod_{j=1}^{d} K_h \left(x^j - X_{t,T}^j\right) \RtT \\[0.5em]
            \frac 1 {T}  \sum_{t=1}^T \left( \frac{\tT - u}{h} \right)  K_h \left(u - \frac{t}{T}\right)   \prod_{j=1}^{d} K_h \left(x^j - X_{t,T}^j\right) \RtT \\[0.5em]
            \frac 1 {T}  \sum_{t=1}^T \left( \frac{\XtT^1 - x^1}{h} \right) K_h \left(u - \frac{t}{T}\right)  \prod_{j=1}^{d} K_h \left(x^j - X_{t,T}^j\right) \RtT \\ \vdots \\ 
            \frac 1 {T}  \sum_{t=1}^T \left( \frac{\XtT^d - x^d}{h} \right) K_h \left(u - \frac{t}{T}\right)  \prod_{j=1}^{d} K_h \left(x^j - X_{t,T}^j\right)\RtT
        \end{pmatrix} = 
        \begin{pmatrix}
            \hat{\Psi}_0(u,x) \\ \hat{\Psi}_1(u,x) \\ \hat{\Psi}_2(u,x) \\ \vdots \\ \hat{\Psi}_{d+1}(u,x)
        \end{pmatrix}.
    \end{aligned}
\end{equation*}
Using the shorthand $\tilde K(v):=vK(v)$ and $\tilde K_h(v):=\frac{1}{h}\tilde K(v/h)$,
we can rewrite the slope–weighted kernel factors as
\[
\Big(\frac{t/T-u}{h}\Big)K_h\!\Big(u-\tfrac{t}{T}\Big)
= -\,\tilde K_h\!\Big(u-\tfrac{t}{T}\Big),\qquad
\Big(\frac{X_{t,T}^k-x^k}{h}\Big)K_h(x^k-X_{t,T}^k)
= -\,\tilde K_h(x^k-X_{t,T}^k).
\]
Thus, each component $\hat{\Psi}_j(u,x)$ is a kernel average of the form 
\begin{equation*}
    \begin{aligned}
\hat\Psi_1(u,x)
&= \frac1T\sum_{t=1}^T\!\Big(\frac{t/T-u}{h}\Big)
K_h\!\Big(u-\tfrac{t}{T}\Big)\prod_{j=1}^d K_h(x^j-X_{t,T}^j)\,R_{t,T}\\
&= -\,\frac1T\sum_{t=1}^T \tilde K_h\!\Big(u-\tfrac{t}{T}\Big)\prod_{j=1}^d K_h(x^j-X_{t,T}^j)\,R_{t,T},
\\[0.4em]
\hat\Psi_{k+1}(u,x)
&= \frac1T\sum_{t=1}^T\!\Big(\frac{X_{t,T}^k-x^k}{h}\Big)
K_h\!\Big(u-\tfrac{t}{T}\Big)\prod_{j=1}^d K_h(x^j-X_{t,T}^j)\,R_{t,T}\\
&= -\,\frac1T\sum_{t=1}^T K_h\!\Big(u-\tfrac{t}{T}\Big)
\tilde K_h(x^k-X_{t,T}^k)\!\!\prod_{j\ne k}\! K_h(x^j-X_{t,T}^j)\,R_{t,T}.
\end{aligned}
\end{equation*}
As we assumed $K$ is bounded, Lipschitz, and compactly supported; therefore $\tilde K(v)=vK(v)$
is also bounded and Lipschitz on the same compact support, and $\tilde K_h$ inherits these
properties. Consequently, each component $\hat\Psi_j(u,x)$, $j=1,\dots,d+1$, has the same
stochastic order and admits the same maximal inequality as $\hat\Psi_0(u,x)$ after replacing
one of the kernels by $\tilde K$. In particular, all arguments used for $\hat\Psi_0$ apply
verbatim to $\hat\Psi_j$.

It suffices to show 
\begin{align*}
    \sup_{u\in[0,1],\,x\in S} \Abs{\hat{\Psi}_0(u,x)-\E[\hat{\Psi}_0(u,x)]}{}
= O_P\Big(\sqrt{\frac{\log T}{T\,h^{d+1}}}\Big), 
\end{align*}
and it is exactly same as the Theorem 4.1 of \cite{Vogt_2012}. 
\end{proof}

\subsection{Proof of \hyperref[Thm3.2]{Theorem 3.2}}\label{prf3.2}

\begin{proof}
We write, 
\begin{equation*}
    \begin{aligned}
        \hat{\frakm}(u,\,x) - \frakm(u,\,x) &= \hat{\frakm}(u,\,x) - \E[\hat{\frakm}(u,\,x) \vert \mbf{X}_T ] + \E[\hat{\frakm}(u,\,x) \vert \mbf{X}_T ] - \frakm(u,\,x)\\
        &= (\mbf{D}^{\top} \mbf{W} \mbf{D})^{-1} \mbf{D}^{\top} \mbf{W} \Big(\mbf{Y}_T\, - \, \mbf{D} \frakm(u,x) \Big)\\
        &=(\mbf{D}^{\top} \mbf{W} \mbf{D})^{-1} \mbf{D}^{\top} \mbf{W} \Big(\mbf{Y}_T\, - \, \mbf{m} \Big) \\
        &\qquad+(\mbf{D}^{\top} \mbf{W} \mbf{D})^{-1} \mbf{D}^{\top} \mbf{W} \Big(\mbf{m}\, - \, \mbf{D} \frakm(u,x) \Big)
    \end{aligned}
\end{equation*}
So, 
\begin{equation*}
    \begin{aligned}
        & \hat{\frakm}(u,x) - \E[{\hat{\frakm}(u,x)}\mid \mbf{X}_{T}] = (\mbf{D}^{\top} \mbf{W} \mbf{D})^{-1} \mbf{D}^{\top} \mbf{W} \Big(\mbf{Y}_T\, - \, \mbf{m} \Big) \\
        & \E[{\hat{\frakm}(u,x)}\mid \mbf{X}_{T}] - \frakm(u,x) = (\mbf{D}^{\top} \mbf{W} \mbf{D})^{-1} \mbf{D}^{\top} \mbf{W} \Big(\mbf{m}\, - \, \mbf{D} \frakm(u,x) \Big)
    \end{aligned}
\end{equation*}

We can prove the theorem by proving the following intermediate results. 

\noindent (i) By Theorem 3.1, with $R_{t,T} =  \varepsilon_{t,T}, \; \mbf{R}_{T} = \mbf{Y}_T - \mbf{m}$, 
\begin{align*}
    \sup_{u \in [0,1] \, x \in S} \Norm{\frac 1 {T} \mbf{D}^{\top} \mbf{W} (\mbf{Y}_T - \mbf{m})}{2}= O_p \Big( \sqrt{\frac{\log T}{Th^{d+1}}} \Big).
\end{align*}
\noindent (ii) Applying the arguments for Theorem 3.1 to $\mbf{R}_{T} = \mbf{m}\, - \, \mbf{D} \frakm(u,x)$ yields
\begin{equation*}
    \begin{aligned}
        \sup_{u \in [0,1] \, x \in S} & \Bigg\| \frac 1 {T} \mbf{D}^{\top} \mbf{W} (\mbf{m}\, - \, \mbf{D} \frakm(u,x)) - \E\left[{\frac 1 {T} \mbf{D}^{\top} \mbf{W} (\mbf{m}\, - \, \mbf{D} \frakm(u,x))}\right] \Bigg\|_{2} \\ &= O_p \Big( \sqrt{\frac{\log T}{Th^{d+1}}} \Big).
    \end{aligned}
\end{equation*}
\noindent (iii) It holds that 
\begin{equation*}
    \begin{aligned}
        \sup_{u \in [0,1], x \in S} &\Bigg\|\E\left[{\frac 1 {T} \mbf{D}^{\top} \mbf{W} (\mbf{m}\, - \, \mbf{D} \frakm(u,x))} \right]\Bigg\|_{2} = O \left( \frac{1}{T^r h^{d-1}} + h^2\right)
    \end{aligned}
\end{equation*}

\noindent (iv) We have that 
\begin{equation*}
    \begin{aligned}
        \sup_{u \in [0,1], x \in S} \Norm{\frac{1}{T}\mbf{D}^{\top} \mbf{W} \mbf{D} - \mbf{M}(u,x) f(u,x)}{F} = o_p(1).
    \end{aligned}
\end{equation*}
where  $\|\cdot\|_{F}$ is frobenius norm defined for $(d+2) \times (d+2)$ matrices and $(d+2)\times(d+2)$ matrix $\mathbf{M(u,x)}$ is entrywise defined by, 
\begin{equation*}
    \begin{aligned}
        &[ \mbf{M} (u,x)]_{0,0} = \mu_{0,0}(u,x) = \frac{1}{T} \sum_{t=1}^{T} K_h\left(u-\tT\right)\int_{S}  \prod_{j=1}^{d} K_h(x^j -z^j )  dz \\
        &[\mbf{M} (u,x)]_{1,1} = \mu_{1,1}(u,x) = \frac{1}{T} \sum_{t=1}^{T}\left( \frac{\tT -u}{h} \right)^2 K_h\left(u-\tT\right)\int_{S}    \prod_{j=1}^{d} K_h(x^j -z^j )  dz \\
        &[\mbf{M} (u,x)]_{k+1,k+1} = \mu_{k+1,k+1}(u,x) = \frac{1}{T} \sum_{t=1}^{T} K_h\left(u-\tT\right)\int_{S} \left( \frac{z^k-x^k}{h}\right)^2   \prod_{j=1}^{d} K_h(x^j -z^j )  dz \\
        &[ \mbf{M} (u,x)]_{0,1} =[\mbf{M} (u,x)]_{1,0} = \mu_{0,1}(u,x) =\mu_{1,0}(u,x) \\& \quad= \frac{1}{T} \sum_{t=1}^{T} \left( \frac{\tT -u}{h} \right) K_h\left(u-\tT\right)\int_{S}     \prod_{j = 1}^d K_h(x^j-z^j)  dz \\
        &[ \mbf{M} (u,x)]_{0,k+1} =[\mbf{M} (u,x)]_{k+1,0} = \mu_{0,k+1}(u,x) =\mu_{k+1,0}(u,x) \\& \quad= \frac{1}{T} \sum_{t=1}^{T} K_h\left(u-\tT\right)\int_{S}   \left( \frac{z^k-x^k}{h}\right)  \prod_{j = 1}^d K_h(x^j-z^j)  dz \\
        &[ \mbf{M} (u,x)]_{1,k+1} =[\mbf{M} (u,x)]_{k+1,1} = \mu_{1,k+1}(u,x) =\mu_{k+1,1}(u,x) \\& \quad= \frac{1}{T} \sum_{t=1}^{T} \left( \frac{\tT -u}{h} \right) K_h\left(u-\tT\right)\int_{S}   \left( \frac{z^k-x^k}{h}\right)  \prod_{j = 1}^d K_h(x^j-z^j)  dz\\
        &[ \mbf{M} (u,x)]_{k+1,k'+1} =[\mbf{M} (u,x)]_{k'+1,k+1} = \mu_{k+1,k'+1}(u,x) =\mu_{k'+1,k+1}(u,x) \\& \quad= \frac{1}{T} \sum_{t=1}^{T}  K_h\left(u-\tT\right)\int_{S}   \left( \frac{z^k-x^k}{h}\right)\left( \frac{z^{k'}-x^{k'}}{h}\right)  \prod_{j = 1}^d K_h(x^j-z^j)  dz 
    \end{aligned}
\end{equation*}
where $z= (z_1, \dots, z_d)^{\top}$, $k, k' =1, \dots, d$, $k \ne k'$.

\noindent (v) There exists $\varepsilon > 0$ such that, with probability tending to one, 
\begin{equation*}
    \hat{\lambda}_{\inf} := \inf_{u \in [0,1], x \in S} \lambda_{\min} \left( \frac{1}{T} \mbf{D}^{\top} \mbf{W} \mbf{D}\right) > \varepsilon
\end{equation*}

Combining the intermediate results (i)-(v) we arrive at 
\begin{equation*}
    \begin{aligned}
    \begin{split}
        \sup_{u \in [0,1]. x \in S} &\left\| \hat{\frakm} (u,x) - \frakm(u,x) \right\|_{2} \\
        & = \sup_{u \in [0,1]. x \in S} \left\|(\mbf{D}^{\top} \mbf{W} \mbf{D})^{-1} \left\{\mbf{D} \mbf{W} \left(\mbf{Y}\, - \, \mbf{m} \right) + \mbf{D} \mbf{W} \left(\mbf{m}\, - \, \mbf{D} \frakm(u,x) \right) \right\} \right\|_{2} \\
        & = \sup_{u \in [0,1]. x \in S}\Bigg[  \Bigg\| \left( \frac{1}{T} \mbf{D}^{\top}\mbf{W}\mbf{D} \right)^{-1} \\ 
        & \qquad \qquad \times   \left\{ \frac 1 {T} \mbf{D}^{\top} \mbf{W} (\mbf{m}\, - \, \mbf{D} \frakm(u,x)) - \E\left[{\frac 1 {T} \mbf{D}^{\top} \mbf{W} (\mbf{m}\, - \, \mbf{D} \frakm(u,x))} \right]  \right.\\
        & \qquad \qquad \qquad \qquad \qquad \left. \left. + \E\left[{\frac 1 {T} \mbf{D}^{\top} \mbf{W} (\mbf{m}\, - \, \mbf{D} \frakm(u,x))} \right] \right\} \Bigg\|_{2}\right]
    \end{split}
    \end{aligned}
\end{equation*}

Then, we can write 
\begin{equation*}
    \begin{aligned}
    \begin{split}
        \sup_{u \in [0,1]. x \in S} &\left\| \hat{\frakm} (u,x) - \frakm(u,x) \right\|_{2} \\
        & \le \sup_{u \in [0,1]. x \in S} \left[ \left\| \frac 1 {T} \mbf{D}^{\top} \mbf{W} (\mbf{m}\, - \, \mbf{D} \frakm(u,x)) - \E\left({\frac 1 {T} \mbf{D}^{\top} \mbf{W} (\mbf{m}\, - \, \mbf{D} \frakm(u,x))} \right)\right\|_{2} \right.\\
        &\qquad \left. + \left\| \E\left({\frac 1 {T} \mbf{D}^{\top} \mbf{W} (\mbf{m}\, - \, \mbf{D} \frakm(u,x))} \right) \right\|_{2} \right] \cdot \hat{\lambda}_{\inf}^{-1} \\
        & = \hat{\lambda}_{\inf}^{-1} \;\cdot\; O_p \left( \sqrt{\frac{\log T}{Th^{d+1}}} + \frac{1}{T^r h^{d-1}} + h^2\right)
    \end{split}
    \end{aligned}
\end{equation*}
with $r = \min\{\rho, 1 \}$. Moreover, (iv), (v) immediately imply that $\hat{\lambda}_{\inf}^{-1} = O_p(1)$. This completes the proof.

\vspace*{1cm}

\subsubsection{Proof of (iii)}
\noindent (iii) It holds that 
\begin{equation*}
    \begin{aligned}
        \sup_{u \in [0,1], x \in S} &\Norm{\E\left[{\frac 1 {T} \mbf{D}^{\top} \mbf{W} (\mbf{m}\, - \, \mbf{D} \frakm(u,x))} \right]}{2} \\& = O \left( \frac{1}{T^rh^{d-1}} + h^2\right)
    \end{aligned}
\end{equation*}

Define $\ZtT = (\tT, \, \XtT)^{\top}$, $Z_t \left( \tT \right) = (\tT , \, X_t (\tT))^{\top}$, $z = (u,x)^{\top}$, $X_t (\tT) = \left(X_t^1 (\tT), X_t^2 (\tT), \dots, X_t^d (\tT)\right)^{\top}$
\begin{equation*}
    \begin{aligned}
        &\frac {1} {T} \mbf{D}^{\top} \mbf{W} (\mbf{m}\, - \, \mbf{D} \frakm(u,x)) \\& 
        = \begin{pmatrix}
            \frac{1}{T} \sum_{t=1}^T K_h \left(u - \frac{t}{T}\right) \prod_{j=1}^{d}K_h \left(x^j - X_{t,T}^j\right) \left( m(\ZtT) - m(z) - \nabla m(z)^{\top} \left( \ZtT - z\right)  \right) \\[0.5em]
            \frac{1}{T} \sum_{t=1}^T \left( \frac{\tT-u}{h} \right) K_h \left(u - \frac{t}{T}\right) \prod_{j=1}^{d}K_h \left(x^j - X_{t,T}^j\right) \left( m(\ZtT) - m(u,x) - \nabla m(z)^{\top} \left( \ZtT - z\right)  \right)\\[0.5em]
            \frac{1}{T} \sum_{t=1}^T \left( \frac{X_{t,T}^1-x^1}{h} \right)  K_h \left(u - \frac{t}{T}\right) \prod_{j=1}^{d}K_h \left(x^j - X_{t,T}^j\right) \left( m(\ZtT) - m(z) - \nabla m(z)^{\top} \left( \ZtT - z\right)  \right) \\ \vdots \\
            \frac{1}{T} \sum_{t=1}^T \left( \frac{X_{t,T}^d-x^d}{h} \right)  K_h \left(u - \frac{t}{T}\right) \prod_{j=1}^{d}K_h \left(x^j - X_{t,T}^j\right) \left( m(\ZtT) - m(z) - \nabla m(z)^{\top} \left( \ZtT - z\right)  \right)
        \end{pmatrix} 
    \end{aligned}
\end{equation*}

Similarly with the proof of Theorem 3.1, applying $\tilde{K}(v) = v \cdot K(v)$, it suffices to show 
\begin{equation*}
    \begin{aligned}
        \E  \left[ \frac{1}{T} \sum_{t=1}^T K_h \left(u - \frac{t}{T}\right) \prod_{j=1}^{d}K_h \left(x^j - X_{t,T}^j\right) \left( m(Z_{t,T}) - m(z) - \nabla m(z)^{\top} \left( Z_{t,T} - z\right)  \right)\right] \\= O \left( \frac{1}{T^r h^{d-1}} + h^2\right).
    \end{aligned}
\end{equation*}

We define $\bar{K} : \bbR \rightarrow \bbR$ be a Lipschitz continuous function with support $[-qC_1, qC_1]$ for some $q > 1$. Assume that $\bar{K}(x) = 1$ for all $x \in [-C_1, C_1]$ and write $\bar{K}_h(x) = \bar{K}( x /h)$. Then, 
\begin{equation*}
    \begin{aligned}
        \E & \left[ \frac{1}{T} \sum_{t=1}^T K_h \left(u - \frac{t}{T}\right) \prod_{j=1}^{d}K_h \left(x^j - X_{t,T}^j\right) \left( m(Z_{t,T}) - m(z) - \nabla m(z)^{\top} \left( Z_{t,T} - z\right)  \right)\right] \\
        &= Q_1(u,x) + Q_2(u,x) + Q_3(u,x) + Q_4(u,x)
    \end{aligned}
\end{equation*}
with 
\begin{align*}
    Q_i(u,x) = \frac{1}{T} \sum_{t=1}^{T} K_h \left( u - \tT \right) q_i(u,x)
\end{align*}
and
\begin{align*}
    q_{1}(z) \;=\; \E \Bigg[ \prod_{j=1}^d \bar{K}_h (x^j - \XtT^j)\, \Bigg\{ \prod_{j=1}^{d}K_h \left(x^j - X_{t,T}^j\right)\;-\;  \prod_{j=1}^{d}K_h\Big(x^j - X_t^j\Big(\frac{t}{T}\Big)\Big) \Bigg\} \\
    \;\times\; \Big\{ m\left(\ZtT\right) \;-\; m(z) - \nabla m(z)^{\top} (\ZtT - z) \Big\} \Bigg],
\end{align*}

\begin{align*}
    q_{2}(z) \;=\; \E &\Bigg[ \prod_{j=1}^d \bar{K}_h (x^j - \XtT^j)\, \prod_{j=1}^{d} K_h\Big(x^j - X_t^j \Big(\tT \Big)\Big) \\
    &\times \Bigg\{ m (\ZtT) - m \left(Z_t \left(\frac t T\right) \right) - \nabla m\left(z \right)^{\top} \left(\ZtT - Z_t \left( \tT \right)  \right) \Bigg\}\Bigg],
\end{align*}

\begin{align*}
    q_{3}(z) &\;=\; \E \Bigg[ \Bigg\{ \prod_{j=1}^d \bar{K}_h(x^j - \XtT^j ) - \prod_{j=1}^d\bar{K}_h \Big(x^j - X_t^j \Big( \tT\Big) \Big) \Bigg\}\prod_{j=1}^{d} K_h \Big( x^j - X_t^j \Big( \tT\Big)  \Big)  \\
    &\times \Bigg\{ m\left(Z_t \left( \frac tT \right) \right) - m(z) - \nabla m(z)^{\top} \left(Z_t \left( \frac tT \right) - z\right)\Bigg\}\Bigg],
\end{align*}

\begin{align*}
    q_{4}(z) \;=\; \E \Bigg[ \prod_{j=1}^{d}K_h \Big(x^j - X_t^j \Big( \tT\Big) \Big) \Bigg\{ m\left(Z_t \left( \frac tT \right) \right) - m(z) - \nabla m(z)^{\top} \left(Z_t \left( \frac tT \right) - z\right)\Bigg\}\Bigg].
\end{align*}
We first consider $Q_1(u,x)$. As the kernel $K$ is bounded and Lipschitz, we can find a constant $C_{K, Lip} < \infty$ and sufficiently large $C$ which varies by line, with 
\begin{align*}
     \Big|\prod_{j=1}^{d} K_h\Big(x^j - \XtT^j \Big) - \prod_{j=1}^{d} K_h\Big(x^j - X_t^j \Big(\frac{t} {T}\Big)\Big) \Big| &= \frac 1 {h^d} \left| \prod_{j=1}^{d} K \left( \frac{x^j - \XtT^j}{h} \right) - \prod_{j=1}^{d} K \left( \frac{x^j - X_t^j\Big(\frac{t} {T}\Big)}{h} \right)\right|\\
     &\le \frac C {h^d} \sum_{j=1}^d \left| K \left( \frac{x^j - \XtT^j}{h} \right) -  K \left( \frac{x^j - X_t^j\Big(\frac{t} {T}\Big)}{h} \right)\right|\\
     &\le \frac C {h^d} \sum_{j=1}^d \left| K \left( \frac{x^j- \XtT^j}{h} \right) - K \left( \frac{x^j - X_t^j\Big(\frac{t} {T}\Big)}{h} \right)\right|^r \\ 
     & \le \frac{C_{K, Lip}}{h^d} \sum_{j=1}^d \left| \frac{\XtT^j - X_t^j(t/T)}{h}\right|^r
\end{align*}
for $r = \min\{\rho, 1\}$. 
Using the above inequality, we obtain
\begin{align*}
    &|Q_1 (u,x) |\\
    &\le \frac{C_{K, Lip}}{T} \sum_{t=1}^T K_h \Big( u - \tT \Big) \\
    &\times \E \left[ \frac{1}{h^d} \sum_{j=1}^d \left|  \frac{\XtT^j - X_t^j(t/T)}{h}\right|^r \times \prod_{j=1}^d \bar{K}_h (x^j - \XtT^j) \Big| m(\ZtT) \;-\; m(z) - \nabla m(z)^{\top} (\ZtT - z)\Big|\right]
\end{align*}
with $r = \min\{\rho, 1\}$. $m$ is Lipschitz and $\bar{K}_h$ is bounded by 1 and, $\prod_{j=1}^d \bar{K}_h (x^j - \XtT^j) \Big| m(\ZtT) \;-\; m(x) - \nabla m(z)^{\top} (\ZtT - z) \Big| \le C_{m,1}h^2 $ for some $C_{m,1} < \infty$. Since $K$ is Lipschitz ($\because$(C5)), $\Big| \XtT^j - X_t^j\Big(\tT\Big) \Big| \le \frac{1}{T} U_{t,T}(\frac{t}{T})$ and the variables $U_{t,T} (\frac{t}{T})$ have finite rth moment, we can infer that 
\begin{align*}
    |Q_1 (u,x) | &\le \frac{C_{K, Lip} C_{m,1} }{Th^{d-2}} \sum_{t=1}^T K_h \Big( u - \tT \Big) \E\left[ \sum_{j=1}^d\left| \frac{\XtT^j - X_t^j(t/T)}{h}\right|^r \right]\\
    &\le \frac{C_{K, Lip} C_{m,1} }{Th^{d-2}} \sum_{t=1}^T K_h \Big( u - \tT \Big) \E\Bigg[\sum_{j=1}^d \Big|\frac 1 {Th}  U_{t,T}\Big(\frac{t}{T}\Big) \Big|^r \Bigg]\\
    &\le\frac{C_K C_{K, Lip} C_{m,1} C_U }{T^r h^{d+r-2}}
\end{align*}
uniformly in $u$ and $x$.

Now, we consider 
\begin{align*}
    Q_2&(u,x) = \frac{1}{T} \sum_{t=1}^{T} K_h \Big( u - \tT \Big) \\ &\qquad \times \E \Bigg[ \prod_{j=1}^d \bar{K}_h (x^j - \XtT^j)\, \prod_{j=1}^dK_h\Big(x^j - X_t^j\Big(\tT \Big)\Big)\\ &\qquad  \times \Bigg\{ m (\ZtT ) - m \Big( Z_t (t/T)\Big) - \nabla m(z)^{\top} \left(\ZtT -Z_t (t/T)\right)  \Bigg\}\Bigg].
\end{align*}
For the formula in the curly bracket, 
\begin{equation*}
    \begin{aligned}
        m (\ZtT ) &- m \Big( Z_t (t/T)\Big) - \nabla m(z)^{\top} \left(\ZtT -Z_t (t/T)\right) \\
        &= m (\ZtT ) - m \Big( Z_t (t/T)\Big) - \nabla m\Big( Z_t (t/T)\Big)^{\top} \Big(\ZtT -Z_t (t/T)\Big) \\
        &\qquad \qquad+  \nabla m\Big( Z_t (t/T)\Big)^{\top} \Big(\ZtT -Z_t (t/T)\Big) - \nabla m(z)^{\top} \left(\ZtT -Z_t (t/T)\right). 
    \end{aligned}
\end{equation*}
Notice that $\ZtT -Z_t (t/T) = \Big(0, \XtT - X_t (t/T)\Big)^{\top}$. Recall that putting $u = t/T$ in the definition of local stationary series, $\Big| \XtT^j - X_t^j\Big(\tT\Big) \Big| \le \frac{1}{T} U_{t,T}(\frac{t}{T})$. 

As, $m$ is twice continuously differentiable, for sufficiently large constant $C_{m,2}$, 
\begin{equation*}
    \begin{aligned}
        \Bigg|m (\ZtT ) - m \Big( Z_t (t/T)\Big) - \nabla m(z)^{\top} \left(\ZtT -Z_t (t/T)\right)  \Bigg| & \le 2C_{m,2}\,C_1 \cdot \left( \frac {h} {T^r}\right)
    \end{aligned}
\end{equation*}

Similarly, for the sufficiently large constant $C$ which may vary line by line, using the boundedness of $K, \;\bar{K}_h$, 
\begin{align*}
    \sup_{u,x }\Big| Q_2(u,x) \Big| &\le \frac{C}{T} \sum_{t=1}^T K_h \Big( u - \tT \Big) \cdot \left( \frac{1}{h^d}\right)\left( \frac {h} {T^r}\right) \\
    &\le \frac{C} {T^r h^{d-1}}.
\end{align*}

For $Q_3 (u,x)$,
\begin{equation*}
    \begin{aligned}
        Q_3(u,x) = &\frac{1}{T} \sum_{t=1}^{T} K_h \Big( u - \tT \Big) \\
        &\E \Bigg[ \Bigg\{ \prod_{j=1}^d \bar{K}_h(x^j - \XtT^j ) - \prod_{j=1}^d\bar{K}_h \Big( x^j - X_t^j \Big( \tT\Big) \Big) \Bigg\} \\
        &\times \prod_{j=1}^dK_h \Big( x^j - X_t^j \Big( \tT\Big)  \Big) \Bigg\{ m \Big(Z_t (t/T) \Big) - m(z) - \nabla m(z)^{\top} \Big(Z_t (t/T)  - z\Big)\Bigg\}\Bigg].
    \end{aligned}
\end{equation*}
Using similar arguments for sufficiently large constant $C < \infty$, 
\begin{equation*}
    \begin{aligned}
        \sup_{u,\,x} |Q_3(u,x)| \le \frac{C}{T^r h^{d+r-2}}
    \end{aligned}
\end{equation*}

Finally, consider
\begin{equation*}
    \begin{aligned}
        Q_4(u,x) &= \frac{1}{T} \sum_{t=1}^{T} K_h \Big( u - \tT \Big) \\ 
        &\times \E \Bigg[ \prod_{j=1}^d K_h \Big(x^j - X_t^j \Big( \tT\Big) \Big) \Bigg\{ m\left(Z_t \left( \frac tT \right) \right) - m(z) - \nabla m(z)^{\top} \left(Z_t \left( \frac tT \right) - z\right)\Bigg\}\Bigg] \\
        & = \frac{1}{T} \sum_{t=1}^{T} K_h \Big( u - \tT \Big) \int_S \prod_{j=1}^d K_h (x^j-y^j) \\& \qquad \qquad  \left\{m\left(\tT, y\right) - m(u,x) - \nabla m(u,x)^{\top} \begin{pmatrix}
            \tT - u \\[0.1em] y^1-x^1 \\ \vdots \\ y^d - x^d
        \end{pmatrix} \right\} f_{X_t(\tT)}(y) dy.
    \end{aligned}
\end{equation*}
for $y = (y_1, \dots , y_d)^{\top}$. Also define 
\begin{equation*}
    \begin{aligned}
        & G(\varphi,y) =  \left\{m\left(\varphi, y\right) - m(u,x) - \nabla m(u,x)^{\top} \begin{pmatrix}
            \varphi - u \\ y^1-x^1 \\ \vdots \\ y^d - x^d
        \end{pmatrix}  \right\} f_{X_t(\varphi)}(y) \;\;\text{ for } \varphi \in [0,1]\\
        &H\left(\varphi, x\right) = \int_S \prod_{j=1}^{d} K_h (x^j-y^j) G \left(\varphi, y\right) dy \\
        & g(\varphi) = K_h ( u - \varphi) H \left(\varphi, x \right).
    \end{aligned}
\end{equation*}

Then, 
\begin{equation*}
    \begin{aligned}
        & \left| \frac{1}{T} \sum_{t=1}^{T} K_h \Big( u - \tT \Big) H \left(\tT, x \right) \right| \\
        & = \left| \frac{1}{T} \sum_{t=1}^{T} K_h \Big( u - \tT \Big)  H \left(\tT, x \right) - \int_{[0,1]} K_h (u-\varphi) H(\varphi, x) d\varphi \right| + \left|  \int_{[0,1]} K_h (u-\varphi) H(\varphi, x) d\varphi \right|
    \end{aligned}
\end{equation*}
We first bound the first absolute term. 
\begin{equation*}
    \begin{aligned}
        \left| \frac{1}{T} \sum_{t=1}^{T} K_h \Big( u - \tT \Big)  H \left(\tT, x \right) - \int_{[0,1]} K_h (u-\varphi) H(\varphi, x) d\varphi \right| &= \left|\frac{1}{T} \sum_{t=1}^{T} g \left( \tT \right) - \int_{0}^1 g(\varphi) d \varphi \right| \\
        & \le \frac{1}{T} \sup_{\varphi \in [0,1]} |g'(\varphi)| = O\left( \frac{1}{T} \right)
    \end{aligned}
\end{equation*}
The last equality holds because 
\begin{equation*}
    \begin{aligned}
        g'(\varphi) &= \frac{d} {d\varphi} g(\varphi) = \frac{d}{d\varphi} \left[ \frac1 h K(\frac{u-\varphi}{h})H(\varphi, x) \right]  \\
        &= - \frac{1}{h^2} K'\left( \frac{u-\varphi}{h}\right) H(\varphi,x) + \frac{1}{h} K \left( \frac{u - \varphi}{h}\right) \partial_{\varphi}H(\varphi, x)
    \end{aligned}
\end{equation*}
For the first term of $g'(\varphi)$, we can substitute $v= \frac{u-\varphi}{h}$ the formula into
\begin{equation*}
    \begin{aligned}
        - \frac{1}{h^2} K'\left( \frac{u-\varphi}{h}\right) H(\varphi,x) = -\frac{1}{h^2} \left( K'(v)H(u-hv, x) \right). 
    \end{aligned}
\end{equation*}
Since the kernel is bounded and Lipschitz and has compact support, 
\begin{equation*}
    \begin{aligned}
        H(u-hv, x) &= \int_S \prod_{j=1}^d K_h (x^j-y^j) \left\{m\left(u-hv, y\right) - m(u,x) - \nabla m(u,x)^{\top} \begin{pmatrix}
            -hv \\ y^1-x^1 \\ \vdots \\ y^d - x^d
        \end{pmatrix} \right\} f_{X_t(u-hv)}(y) dy \\
        &= O(h^2).
    \end{aligned}
\end{equation*}
So, 
\begin{equation*}
    - \frac{1}{h^2} K'\left( \frac{u-\varphi}{h}\right) H(\varphi,x) = O(1).     
\end{equation*}

For the second term $\frac{1}{h} K \left( \frac{u - \varphi}{h}\right) \partial_{\varphi}H(\varphi, x)$,
\begin{equation*}
    \begin{aligned}
        \partial_{\varphi}H(\varphi, x) = \int_S \prod_{j=1}^d K_h(x^j - y^j) \, \frac{\partial G}{\partial \varphi}(\varphi, y) \, dy.
    \end{aligned}
\end{equation*}
And since $K_h (u-\varphi), \;K_h(x^j - y^j) \ne 0$, 
\begin{equation*}
    \begin{aligned}
        \frac{\partial G}{\partial \varphi}(\varphi, y) = (\partial_0 m(\varphi, y) - \partial_0 m(u, x))\cdot f_{X_t(v)}(y) = O(h).
    \end{aligned}
\end{equation*}
As a result, $\sup_{\varphi \in [0,1]}|g'(\varphi) | = O(1) $. 
For the term $ \left|  \int_{[0,1]} K_h (u-\varphi) H(\varphi, x) d\varphi \right|$, 
\begin{equation*}
    \begin{aligned}
        \int_{[0,1]} & K_h (u-\varphi) H(\varphi, x) d\varphi    \\ 
        & = \int_{[0,1] } \int_{S } K_h (u-\varphi) \prod_{j=1}^d K_h(x^j - y^j) \\ & \qquad \qquad \qquad \qquad  \left\{m\left(\varphi, y\right) - m(u,x) - \nabla m(u,x)^{\top} \begin{pmatrix}
            \varphi - u \\ y^1-x^1 \\ \vdots \\ y^d - x^d
        \end{pmatrix} \right\} f(\varphi, y) \; dy d \varphi \\
        & = \int_{[0,1] \times S} K(w_0) \prod_{j=1}^d K(w_j)\left\{ \frac{h^2}{2} (w^{\top} \nabla^2 m(z) w) +o(h^2) \right\} \\ & \qquad \qquad \qquad \qquad \{ f(u,x) + \nabla f(u,x)^{\top} \cdot hw+o(h) \} dw \\
        & \left(\because \text{integration by substitution}, w_0 = \frac{u - \varphi}{h}, \; w_j = \frac{x^j-y^j}{h}, \; w = (w_0, w_1, \dots, w_d)^{\top}\right)\\
        & = O(h^2) 
    \end{aligned}
\end{equation*}
This completes the proof of 
\begin{align}
    \sup_{u \in [0,1], x \in S} \Abs{Q_4(u,x)}{} = O\left(\frac{1}{T} + h^2\right)
\end{align}

Combining the results on $Q_1(u,x), \, Q_2(u,x), \, Q_3(u,x), \, Q_4(u,x)$ yields
\begin{equation*}
    \begin{aligned}
        \sup_{u \in [0,1], x \in S} &\left\|\E\left[{\frac 1 {T} \mbf{D}^{\top} \mbf{W} (\mbf{m}\, - \, \mbf{D} \frakm(u,x))}\right]\right\|_{2} \\
        & = O \left( \frac{1}{T^rh^{d-1}} + h^2\right)
    \end{aligned}
\end{equation*}

and it completes the proof of (iii). 

\vspace*{1cm}

\subsubsection{Proof of (iv)}
\noindent (iv) We write, 
\begin{equation*}
    \begin{aligned}
        \sup_{u \in [0,1], x \in S} \Big\|{\frac{1}{T}\mbf{D}^{\top} \mbf{W} \mbf{D} - \mbf{M}(u,x) f(u,x)}\Big\|_{F} = o_p(1)
    \end{aligned}
\end{equation*}

For the proof, we define $B := \{(u,x) \in \bbR^{d+1} : u \in [0,1], x \in S\}$, \(a_T = \sqrt{\frac{\log T}{T h^{d+1}}}\).

For each element of $\frac{1}{T}\mbf{D}^{\top} \mbf{W} \mbf{D}$, we can write as 
\begin{equation*}
    \frac{1}{T}  \sum_{t=1}^T \left( \frac{\tT - u}{h} \right)^{l_0} \left( \frac{\XtT^k - x^k}{h} \right)^{l} \left( \frac{\XtT^{k'} - x^{k'}}{h} \right)^{l'} K_h \left( u - \tT \right) \prod_{j=1}^d K_h \left( x^j - \XtT^j \right)
\end{equation*} for $k, k'=1, 2, \dots , d$, $k \ne k'$, $l_0, l, l' = 0, 1, 2$.
For each element of $\mbf{M}(u,x)$, we can write as
\begin{equation*}
    \frac{1}{T} \sum_{t=1}^T \left( \frac{\tT - u}{h} \right)^{l_0} K_h \left(u -\tT \right)\int_{S} \left( \frac{z^{k}-x^{k}}{h}\right)^{l} \left( \frac{z^{k'}-x^{k'}}{h}\right)^{l'} \prod_{j=1}^d K_h(x^j-z^j) dz
\end{equation*}for $k=1, 2, \dots , d$, $l = 0, 1, 2$.

Now for the convenience of the proof we define 
\begin{equation*}
    \begin{aligned}
         &\hat{\psi}(u,x) := \frac{1}{T}  \sum_{t=1}^T \left( \frac{\tT - u}{h} \right)^{l_0} \left( \frac{\XtT^{k} - x^{k}}{h} \right)^{l} \left( \frac{\XtT^{k'} - x^{k'}}{h} \right)^{l'} K_h \left( u - \tT \right) \prod_{j=1}^d K_h \left( x^j - \XtT^j \right) \\
         & \mu(u,x) := \frac{1}{T} \sum_{t=1}^T \left( \frac{\tT - u}{h} \right)^{l_0} K_h \left(u -\tT \right)\int_{S} \left( \frac{z^{k}-x^{k}}{h}\right)^{l} \left( \frac{z^{k'}-x^{k'}}{h}\right)^{l'}  \prod_{j=1}^d K_h(x^j-z^j) dz \\
         & \pi(u,x) := \mu(u,x)f(u,x).
    \end{aligned}
\end{equation*}

Cover the region \( B \subset \mathbb{R}^{d+1} \) with
\(N = O( {(a_T h)}^{-{d+1}})\)
balls of the form
\[
B_n = \left\{ (u, x) \in \mathbb{R}^{d+1}: \|(u, x) - (u_n, x_n)\|_{2} \le a_T h \right\},
\]

where \( (u_n, x_n) \) denotes the center of the ball \( B_n \).

Define the kernel envelope function \( K^*: \mathbb{R} \to \mathbb{R} \) by
\[
K^*(v) = C^*\, \mathbb{I}(|v| \le 2C_1), \; K_h^*(v) = \frac{1}{h} K^*\left( \frac{v}{h} \right).
\]
where \( C_1 < \infty \) is a constant such that the kernel function \( K(v) = 0 \) for all \( |v| > 2C_1 \), and \( C^* > 0 \) is chosen large enough to ensure the following inequality holds.

For all \( (u, x) \in B_n \) and sufficiently large \( T \), the following holds:
\[
\begin{aligned}
\Bigg| &\left( \frac{\tT - u}{h} \right)^{l_0} \left( \frac{\XtT^{k} - x^{k}}{h} \right)^{l} \left( \frac{\XtT^{k'} - x^{k'}}{h} \right)^{l'} K_h\left(u - \frac{t}{T} \right) \prod_{j=1}^d K_h(x^j - \XtT^j) \\
&\quad - \left( \frac{\tT - u_n}{h} \right)^{l_0} \left( \frac{\XtT^{k} - x_n^{k}}{h} \right)^{l} \left( \frac{\XtT^{k'} - x_n^{k'}}{h} \right)^{l'} K_h(u_n - \tT) \prod_{j=1}^d K_h(x_n^j - \XtT^j) \Bigg| \\
&\le a_T \cdot K_h^*(u_n - \tT )\; \prod_{j=1}^d K_h^*(x_n^j - \XtT^j).
\end{aligned}
\]

Define 
\begin{equation*}
    \begin{aligned}
        &\tilde{\psi} (u, x) = \frac 1 {T} \sum_{t=1}^T K_h^*(u_n - \tT )\; \prod_{j=1}^d K_h^*(x_n^j - \XtT^j). \\ 
    \end{aligned}    
\end{equation*}
Then, 
\begin{equation*}
    \begin{aligned}
        \E\left[ \tilde{\psi} (u, x)\right] &= \frac{1}{T} \sum_{t=1}^{T} K_h^* \left( u- \tT \right) \int_S \prod_{j=1}^d K_h^* (x^j - z^j) f_{\XtT}(z) dz \\
        & = \frac{1}{Th} \sum_{t=1}^{T} I \left( \left|u -\tT \right| \le 2C_1 h\right) \int_S  \prod_{j=1}^d I \left( |\omega_j| \le 2C_1 \right) f_{\XtT}(x - h \omega) d\omega \\&\left(\omega_j = \frac{x^j - z^j}{h}, \;\; {\omega}^{\top} = (\omega_1, \omega_2, \dots, \omega_d)\right) \\& \le M.
    \end{aligned}
\end{equation*}
This is because at most $O(Th)$ of $u$ are $I \left( \left|u -\tT \right| \le 2C_1 h\right) \ne 0$. Which leads to $\sup_{u \in [0,1], x \in S}\E | \tilde{\psi}(u, x)|\le M < \infty$ for some sufficiently large $M$. We obtain 
\begin{align*}
    \sup_{(u,x) \in B_n}& \Big| \hat{\psi}(u,x) - \pi(u,x) \Big|
    \\& \le  \Big|\hat{\psi}(u_n, x_n ) - \pi(u_n, x_n ) \Big| + \sup_{(u,x) \in B_n} \left\{ \Big|\hat{\psi}(u, x ) - \hat{\psi}(u_n, x_n ) \Big| + \Big|{\pi}(u_n, x_n ) - \pi(u, x ) \Big| \right\} \\
    &\le \Big|\hat{\psi}(u_n, x_n ) - \E \hat{\psi}(u_n, x_n ) \Big| +\Big|\E \hat{\psi}(u_n, x_n ) - \pi(u_n, x_n ) \Big| + a_T  \Big|\tilde{\psi}(u_n, x_n )\Big| + C a_T h \\
    & \Big(\because \Big| \hat{\psi}(u,x) - \hat{\psi}(u_n, x_n )\Big| \le a_T \Big| \tilde{\psi}(u_n, x_n ) \Big|, \hat{\psi} \text{ and } \pi \text{ are Lipschitz continuous}\Big) \\
    &\le \Big|\hat{\psi}(u_n, x_n ) - \E \hat{\psi}(u_n, x_n ) \Big| + a_T \Big( \Big|\tilde{\psi}(u_n, x_n )\Big| - \E\Big| \tilde{\psi}(u_n, x_n ) \Big| \Big) + Ma_T + C a_T h + o_p(1)\\
    & \le \Big|\hat{\psi}(u_n, x_n ) - \E \hat{\psi}(u_n, x_n ) \Big| + \Big| \tilde{\psi}(u_n, x_n) - \E \tilde{\psi}(u_n, x_n) \Big| + Ma_T + C a_T h+  o_p(1).
\end{align*}
The third inequality is because 
\begin{equation*}
    \begin{aligned}
        &\Big|\E [\hat{\psi}(u_n, x_n ) ]- \pi(u_n, x_n )\Big| = o_p(1),
    \end{aligned}
\end{equation*}
so for all $\varepsilon > 0 \quad\lim_{T \rightarrow \infty} \Prob \left(\Big|\E [\hat{\psi}(u_n, x_n )]- \pi(u_n, x_n ) \Big| > \varepsilon\right) = 0$. We will prove this later.
As a consequence, 
\begin{align*}
    &\Prob \Big(\sup_{(u,x) \in B} \Big| \hat{\psi}(u,x) - \pi(u,x) \Big| > 4M a_T + \varepsilon\Big) \\
    &\quad \le N \max_{1 \le n \le N} \Prob \Big( \sup_{(u,x) \in B_n} \Big| \hat{\psi}(u,x) - \E[\hat{\psi}(u,x)  ]\Big| > 4M a_T + \varepsilon \Big) \le \hat{Q}_T + \tilde{Q}_T
\end{align*}
with
\begin{align*}
    & \hat{Q}_T =  N \max_{1 \le n \le N} \Prob \Big(  \Big| \hat{\psi}(u_n,x_n) - \E[\hat{\psi}(u_n,x_n) ] \Big| > M a_T \Big) \\
    &\tilde{Q}_T = N \max_{1 \le n \le N} \Prob \Big( \Big| \tilde{\psi}(u_n,x_n) - \E[\tilde{\psi}(u_n,x_n)]\Big| > M {a_T} \Big)
\end{align*}
\\
As $\hat{Q_T}$ and $\tilde{Q}_T$ can be analyzed in the same way, we restrict attention to $\hat{Q}_T$ in what follows. To bound $\hat{Q}_T$ we write 
\begin{align}
    \Prob \Big( |\hat{\psi}(u,x) -\E[\hat{\psi}(u,x)  ]  | > Ma_T \Big) \\
    = \Prob \Big(\Big| \sum_{t=1}^T Z_{t,T}(u,x)  \Big| > M a_T T  \Big) \nonumber
\end{align}
with 
\begin{equation*}
    \begin{aligned}
        Z_{t,T}(u,x) = \left( \frac{\tT - u}{h} \right)^{l_0} K_h\Big(u - \tT \Big) \Bigg\{ \left( \frac{\XtT^{k} - x^{k}}{h} \right)^{l} \left( \frac{\XtT^{k'} - x^{k'}}{h} \right)^{l'} \prod_{j=1}^d K_h(x^j - \XtT^j)  \\ -\E \Big[\left( \frac{\XtT^{k} - x^{k}}{h} \right)^{l} \left( \frac{\XtT^{k'} - x^{k'}}{h} \right)^{l'} \prod_{j=1}^d K_h(x^j - \XtT^j)    \Big] \Bigg\}.
    \end{aligned}
\end{equation*}
Note that $Z_{t,T}(u,x)$ is a function of $X_{t,T}$. Thus, the array ${Z_{t,T}(u,x)}$ inherits the mixing properties of ${X_{t,T}}$. The assumptions for the theorems in the paper state that the array involving $\{\XtT\}$ is $\alpha$-mixing with coefficients denoted by $\alpha(k)$. This implies that the process $\{\ZtT\}$ is itself $\alpha$-mixing. 

Since,  $Z_{t,T}(u,x)$ is a measurable function of the random variable $\XtT$, $Z_{t,T}(u,x)$ inherits the mixing property with mixing coefficients $\alpha_T^Z$ satisfying $\alpha_T^Z(k) \le \alpha(k)$. This allows us to apply the following lemma. 

\begin{lemma} [Theorem 2.1 in \cite{Liebscher1996}] \label{lemmaA.1}
We define $Z_{t,T}$ be a zero-mean triangular array such that 
$|Z_{t,T}| \le b_T$ with strong mixing coefficients~$\alpha(k)$.
Then for any $\varepsilon > 0$ and $S_T \le T$, with $\varepsilon > 4S_T b_T$
\begin{align*}
    \mathbb{P}\Big(\,\Big|\sum_{t=1}^T Z_{t,T}\Big| > \varepsilon\Big) \; \le \;4 \,\exp\!\Big(-\,\frac{\varepsilon^2}{64 \sigma_{S_t, T}^2\,T^{\frac{T}{S_T}}\;+\; \frac 8 3\varepsilon b_T S_T} \Big) \;+\; 4\frac T {S_T}\,\alpha(S_T).
\end{align*}
where $\sigma_{S_T, T}^2 = \sup_{0 \le j \le T-1} \E [(\sum_{t = j+1}^{\min\{j + S_T, T\}}Z_{t,T} )^2]$.
\end{lemma}
We apply this exponential inequality with $\varepsilon = Ma_T T $, $b_T = C / h^{d+1} $ for some sufficiently large $C$, and $S_T = a_T^{-1} $. Moreover, an extension of Theorem 1 in \cite{Hansen_2008} shows that $\sigma_{S_T , T}^2 \le \Theta S_T / h^{d+1}$ with a constant $\Theta$ independent of $(u,x)$. We now get,
\begin{align*}
    &\Prob \Big( \Big| \sum_{t=1}^T Z_{t,T} (u,x) \Big| >  Ma_T T\Big) \\
    & \le 4 \exp \Big( - \frac{\varepsilon^2}{64\Theta S_T  \frac{T}{S_Th^{d+1}} + \frac{8}{3} \varepsilon S_T b_T}  \Big) + 4 \frac{T}{S_T} \alpha(S_T)\\
    &\le 4 \exp \Big( - \frac{M^2 \log T}{64 \Theta + \frac{8}{3}CM} \Big) + 4 \frac{T}{S_T} A S_T^{-\beta} \qquad ( \because \text{first term is just calculation, second term is from (K2)} )\\
    &\le 4 \exp \Big( - \frac{M \log T}{64  + \frac{8}{3}C} \Big) + 4 AT S_T^{-1-\beta} \qquad ( \because \text{we choose} \; M > \Theta) \\
    &= 4T^{-\frac{M}{64 + 3C}} + 4 AT S_T^{-1-\beta}
\end{align*}
Since, $ N \le C h^{-(d+1)} a_T^{-(d+1)}$, it follows that 
\begin{align*}
    \hat{Q}_T \le O(R_{1T}) + O(R_{2T})
\end{align*}
with 
\begin{align*}
    &R_{1T} = h^{-(d+1)} a_T^{-(d+1)} T^{- \frac{M}{64+3C}}\\
    &R_{2T} = h^{-(d+1)} a_T^{-(d+1)} T S_T^{-1-\beta}
\end{align*}
Choosing $M$ sufficiently large, we obtain that $R_{1T} \le T^{-\eta}$ for some small $\eta >0$. As $\frac{ \log \log T \,\log T}{T^{\theta}\,h^{d+1}} = o(1) $ by assumption, we further get that 
\begin{align*}
    R_{2T} &= h^{-(d+1)} a_T^{-(d+1)} T{(a_T S_T)}^{1+\beta}\\
    &= \Big(\frac{1}{h^{d+1}}\Big)^{1+\frac{\beta - d}{2}} T^{1-\frac{\beta -d}{2}} \\
    &= o (T^{\theta(1+\frac{\beta - d}{2}) + 1-\frac{\beta -d}{2}}) = o(T^{-\frac{1+\beta}{s}}) = o(1).
\end{align*}
This yield the result
\begin{equation*}
    \begin{aligned}
        \Prob \Big(\sup_{(u,x) \in B_n} \Big| \hat{\psi}(u,x) - \E[\hat{\psi}(u,x)  ]\Big| > 4M a_T + \varepsilon \Big)  \longrightarrow 0.
    \end{aligned}
\end{equation*}
In particular, since each entry of $T^{-1}\mathbf D^{\top}(u,x)\mathbf W(u,x)\mathbf D(u,x)$ is a kernel average of the same type as $\hat{\psi}(u,x)$, the preceding display implies the uniform convergence of the local design matrix:
\begin{equation*}
    \begin{aligned}
        \sup_{u \in [0,1], x \in S} \Bigg\|\frac{1}{T}\mbf{D}^{\top} \mbf{W} \mbf{D}- \E\Big[\frac{1}{T}\mbf{D}^{\top} \mbf{W} \mbf{D}    \Big]  \Bigg\|_F = o_p(1)
    \end{aligned}
\end{equation*}
where $\|\cdot\|_{F}$ denotes the Frobenius norm.

Now, we prove the part that we missed. 
\begin{equation*}
    \begin{aligned}
        &\Big|\E [\hat{\psi}(u, x ) ] - \pi(u, x ) \Big| \\
        &=\Bigg|\E  \Bigg[ \frac{1}{T}  \sum_{t=1}^T \left( \frac{\tT - u}{h} \right)^{l_0} \left( \frac{\XtT^{k} - x^{k}}{h} \right)^{l} \left( \frac{\XtT^{k'} - x^{k'}}{h} \right)^{l'} K_h \left( u - \tT \right) \prod_{j=1}^d K_h \left( x^j - \XtT^j \right)    \Bigg] \\ & \qquad \qquad \qquad \qquad - \mu(u,x) f(u,x)\Bigg| = o_p(1)
    \end{aligned}
\end{equation*} 
for pointwise $(u,x)$ can be analyzed. 

To be specific 
\begin{equation*}
    \begin{aligned}
        \E & \Bigg[ \! \frac{1}{T}  \sum_{t=1}^T \left( \frac{\tT - u}{h} \right)^{l_0} \left( \frac{\XtT^{k} - x^{k}}{h} \right)^{l} \left( \frac{\XtT^{k'} - x^{k'}}{h} \right)^{l'} \! K_h \left( u - \tT \right) \prod_{j=1}^d K_h \left(x^j - \XtT^j \right) \! \Bigg] - \mu(u,x) f(u,x)\\
        & = \E \Bigg[ \frac{1}{T}  \sum_{t=1}^T \left( \frac{\tT - u}{h} \right)^{l_0} \left( \frac{\XtT^{k} - x^{k}}{h} \right)^{l} \left( \frac{\XtT^{k'} - x^{k'}}{h} \right)^{l'} K_h \left( u - \tT \right) \prod_{j=1}^d K_h\left(x^j - \XtT^j \right)    \Bigg] \\
        & \qquad- \frac{1}{T} \sum_{t=1}^{T} \left( \frac{\tT - u}{h} \right)^{l_0} \frac{1}{h}K_h(u - \tT) \int_{S}  \left( \frac{z^{k}-x^{k}}{h}\right)^{l} \left( \frac{z^{k'}-x^{k'}}{h}\right)^{l'} \prod_{j=1}^d K_h(x^j-z^j) dz \; \cdot f(u,x) \\
        & = \E\Bigg[\frac{1}{T} \sum_{t=1}^T \left( \frac{\tT - u}{h} \right)^{l_0} \left( \frac{\XtT^{k} - x^{k}}{h} \right)^{l} \left( \frac{\XtT^{k'} - x^{k'}}{h} \right)^{l'} K_h \left( u - \tT \right) \prod_{j=1}^d K_h \left(x^j - \XtT^j \right)    \Bigg] \\
        & \qquad- \frac{1}{T} \sum_{t=1}^{T} \left( \frac{\tT - u}{h} \right)^{l_0} K_h(u - \tT) \int_{S}  \left( \frac{z^{k}-x^{k}}{h}\right)^{l} \left( \frac{z^{k'}-x^{k'}}{h}\right)^{l'} \prod_{j=1}^d K_h(x^j-z^j) f(\tT,z) dz\; +\;o(1) 
    \end{aligned}
\end{equation*}

The last equality holds because since $f$ is continuously differentiable, by the taylor expansion of $f$, $f(\tT,z) = f(u,x) + O(h)$ and $O(h)$ is $o(1)$. 

Since, $X_t(u)$ is strictly stationary process with density $f(u,x) := f_{X_t(u)} (x)$, 
\begin{equation*}
    \begin{aligned}
        \E &\Bigg[ \frac{1}{T}  \sum_{t=1}^T \left( \frac{\tT - u}{h} \right)^{l_0} \left( \frac{X_t^{k}\left(\tT\right) - x^{k}}{h} \right)^{l} \left( \frac{X_t^{k'}\left(\tT\right) - x^{k'}}{h} \right)^{l'} K_h \left( u - \tT \right) \prod_{j=1}^d K_h \left(x^j - X_t^j \left(\tT\right) \right)    \Bigg] \\
        &= \frac{1}{T}  \sum_{t=1}^T  \left( \frac{\tT - u}{h} \right)^{l_0}  K_h \left( u - \tT \right) \Bigg[ \int_{S} \left( \frac{z^k-x^k}{h}\right)^{l} \left( \frac{z^{k'}-x^{k'}}{h}\right)^{l'} \prod_{j=1}^d K_h \left(x^j - z^j \right) f_{X_t(\tT)}(z) dz \Bigg].
    \end{aligned}
\end{equation*}

We define $\bar{K} : \bbR \rightarrow \bbR$ be a Lipschitz continuous function with support $[-qC_1, qC_1]$ for some $q > 1$. Assume that $\bar{K}(x) = 1$ for all $x \in [-C_1, C_1]$ and write $\bar{K}_h(x) = \bar{K}( \frac x h)$. Then, 
\begin{equation*}
    \begin{aligned}
        &\E \Bigg[ \frac{1}{T}  \sum_{t=1}^T \left( \frac{\tT - u}{h} \right)^{l_0} \left( \frac{\XtT^{k} - x^{k}}{h} \right)^{l} \left( \frac{\XtT^{k'} - x^{k'}}{h} \right)^{l'} K_h \left( u - \tT \right) \prod_{j=1}^d K_h\left(x^j - \XtT^j \right) \\ 
        & \quad-\frac{1}{T}  \sum_{t=1}^T \left( \frac{\tT - u}{h} \right)^{l_0} \left( \frac{X_t^{k}\left(\tT\right) - x^{k}}{h} \right)^{l} \left( \frac{X_t^{k'}\left(\tT\right) - x^{k'}}{h} \right)^{l'} K_h \left( u - \tT \right) \prod_{j=1}^d K_h \left(x^j - X_t^j \left(\tT\right) \! \right) \! \Bigg] \\
        & = Q_1 (u,x) + Q_2(u,x) +Q_3(u,x)
    \end{aligned}
\end{equation*}
with 
\begin{equation*}
    \begin{aligned}
         Q_i(u,x) = \frac{1}{T} \sum_{t=1}^T \left( \frac{\tT - u}{h} \right)^{l_0}  K_h \left( u - \tT \right) q_i(u,x)
    \end{aligned}
\end{equation*}
and
\begin{align*}
    q_{1}(u,x) \;=\; \E \Bigg[\prod_{j=1}^d \bar{K}_h (x^j - \XtT^j)\, \left\{\prod_{j=1}^d K_h\Big(x^j - \XtT^j \Big)\;-\;  \prod_{j=1}^d K_h\Big(x^j - X_t^j\Big(\frac{t}{T}\Big)\Big) \right\} 
    \\ \;\times\; \left( \frac{\XtT^{k} - x^{k}}{h} \right)^{l} \left( \frac{\XtT^{k'} - x^{k'}}{h} \right)^{l'}\; \Bigg],
\end{align*}

\begin{align*}
    q_{2}(u,x) \;&=\; \E \Bigg[\prod_{j=1}^d \bar{K}_h (x^j - \XtT^j)\, \prod_{j=1}^d K_h\Big(x^j - X_t^j\Big(\tT \Big)\Big) \qquad \qquad \\ &\times \Bigg\{\left( \frac{\XtT^{k} - x^{k}}{h} \right)^{l} \left( \frac{\XtT^{k'} - x^{k'}}{h} \right)^{l'} - \left( \frac{X_t^{k}\left(\tT\right) - x^{k}}{h} \right)^{l} \left( \frac{X_t^{k'}\left(\tT\right) - x^{k'}}{h} \right)^{l'} \Bigg\}\, \Bigg],
\end{align*}

\begin{align*}
    q_{3}(u,x) \;&=\; \E \Bigg[ \Bigg\{\prod_{j=1}^d \bar{K}_h(x^j - \XtT^j ) - \prod_{j=1}^d \bar{K}_h \Big(x^j - X_t^j \Big( \tT\Big) \Big) \Bigg\} \qquad \qquad \\ &\times \prod_{j=1}^d K_h \Big( x^j - X_t^j \Big( \tT\Big)  \Big) \left( \frac{X_t^{k}\left(\tT\right) - x^{k}}{h} \right)^{l} \left( \frac{X_t^{k'}\left(\tT\right) - x^{k'}}{h} \right)^{l'} \,\Bigg],
\end{align*}

We first consider $Q_1(u,x)$. As the kernel $K$ is bounded and Lipschitz, we can find a constant $C_{K, Lip} < \infty$ and sufficiently large $C$ which varies by line, with 
\begin{align*}
     \Big|\prod_{j=1}^{d} K_h\Big(x^j - \XtT^j \Big) - \prod_{j=1}^{d} K_h\Big(x^j - X_t^j \Big(\frac{t} {T}\Big)\Big) \Big| &= \frac 1 {h^d} \left| \prod_{j=1}^{d} K \left( \frac{x^j - \XtT^j}{h} \right) - \prod_{j=1}^{d} K \left( \frac{x^j - X_t^j\Big(\frac{t} {T}\Big)}{h} \right)\right|\\
     &\le \frac C {h^d} \sum_{j=1}^d \left| K \left( \frac{x^j - \XtT^j}{h} \right) -  K \left( \frac{x^j - X_t^j\Big(\frac{t} {T}\Big)}{h} \right)\right|\\
     &\le \frac C {h^d} \sum_{j=1}^d \left| K \left( \frac{x^j- \XtT^j}{h} \right) - K \left( \frac{x^j - X_t^j\Big(\frac{t} {T}\Big)}{h} \right)\right|^r \\ 
     & \le \frac{C_{K, Lip}}{h^d} \sum_{j=1}^d \left| \frac{\XtT^j - X_t^j(t/T)}{h}\right|^r
\end{align*}
for $r = \min\{\rho, 1\}$. 

Using above inequality, we obtain
\begin{align*}
    |Q_1 &(u,x) |\\
    &\le \frac{C_{K, Lip}}{T} \sum_{t=1}^T  \left( \frac{\tT - u}{h} \right)^{l_0} K_h \Big( u - \tT \Big) \\ & \qquad \times \E\left[\frac{1}{h^d}\sum_{j=1}^d\left| \frac{\XtT^j - X_t^j(t/T)}{h}\right|^r \times \prod_{j=1}^d \bar{K}_h (x^j - \XtT^j) \left( \frac{\XtT^{k} - x^{k}}{h} \right)^{l} \left( \frac{\XtT^{k'} - x^{k'}}{h} \right)^{l'} \right]
\end{align*}
with $r = \min\{ \rho, 1 \}$. The term $\prod_{j=1}^d \bar{K}_h (x^j - \XtT^j) \left( \frac{\XtT^{k} - x^{k}}{h} \right)^{l} \left( \frac{\XtT^{k'} - x^{k'}}{h} \right)^{l'}$ in the above expression can be bounded. Since, $K$ is Lipschitz, $|\XtT^j - X_t^j (\tT)| \le \frac{1}{T} U_{t,T}(\tT)$  and the variables $U_{t,T} (\frac{t}{T})$ have finite rth moment, we can infer that for sufficiently large constant $C$ which may vary line by line, 

\begin{align*}
    |Q_1 &(u,x) |\\
    &\le \frac{C_{K, Lip}}{Th^d} \sum_{t=1}^T \left( \frac{\tT - u}{h} \right)^{l_0}   K_h \Big( u - \tT \Big) \E\Bigg[ \sum_{j=1}^d \Big|\frac 1 {Th}  U_{t,T}\Big(\frac{t}{T}\Big) \Big|^r \Bigg]\\
    &\le\; \frac{C_K C_{K, Lip} C_U}{T^r h^{r+d}}
\end{align*}

Similarly, 
\begin{align*}
    |Q_2 &(u,x) |\\
    &\le \frac{C}{T} \sum_{t=1}^T \left( \frac{\tT - u}{h} \right)^{l_0} K_h \Big( u - \tT \Big) \\ 
    & \qquad \E\Bigg[\prod_{j=1}^d\bar{K}_h (x^j - \XtT^j)\, \prod_{j=1}^d K_h\Big(x^j - X_t^j\Big(\tT \Big)\Big) \times \sum_{j=1}^d \Big|\frac 1 {T}  U_{t,T}\Big(\frac{t}{T}\Big) \Big|^r \Bigg] \\
    & \le \frac{C }{T^r h^{d}} 
\end{align*}
and, 
\begin{align*}
    |Q_3 &(u,x) |\\
    &\le \frac{C}{T} \sum_{t=1}^T \left( \frac{\tT - u}{h} \right)^{l_0} K_h \Big( u - \tT \Big) \\ 
    & \qquad \E\Bigg[ \Bigg\{ \prod_{j=1}^d \bar{K}_h(x^j - \XtT^j) - \prod_{j=1}^d \bar{K}_h \Big(x^j - X_t^j \Big( \tT\Big) \Big) \Bigg\} \\ & \qquad \qquad \times \prod_{j=1}^d K_h \Big( x^j - X_t^j \Big( \tT\Big)  \Big) \left( \frac{X_t^{k}\left(\tT\right) - x^{k}}{h} \right)^{l} \left( \frac{X_t^{k'}\left(\tT\right) - x^{k'}}{h} \right)^{l'} \Bigg] \\
    & \le\frac{C}{T^r h^{r+d}}
\end{align*}
Combining the results on $Q_1(u,x), \; Q_2(u,x), \; Q_3(u,x)$ yields 
\begin{align*}
    \sup_{u \in [0,1], \; x \in S} \E \Bigg[ \frac{1}{T}  \sum_{t=1}^T \left( \frac{\tT - u}{h} \right)^{l_0} \left( \frac{\XtT^{k} - x^{k}}{h} \right)^{l} \left( \frac{\XtT^{k'} - x^{k'}}{h} \right)^{l'} K_h \left( u - \tT \right) \prod_{j=1}^d K_h \left( x^j - \XtT^j \right) \Bigg] \qquad \\ - \mu(u,x) f(u,x) = o_p(1)
\end{align*}

\medskip

\subsubsection{Proof of (v)}
\noindent (v) There exists $\varepsilon > 0$ such that, with probability tending to one, 
\begin{equation*}
    \hat{\lambda}_{\inf} := \inf_{u \in [0,1], x \in S} \lambda_{\min} \left( \frac{1}{T} \mbf{D}^{\top} \mbf{W} \mbf{D}\right) > \varepsilon
\end{equation*}

First prove that 
\begin{equation*}
    \inf_{u \in [0,1], x \in S } \lambda_{\min} (\mbf{M} (u,x)) >0.
\end{equation*}
We begin by observing a quadratic-form identity. For 
$\omega=(\alpha_0,\alpha_1,\beta^\top)^\top \in \mathbb R^{d+2}$, with 
$\alpha_0,\alpha_1 \in \mathbb R$ and $\beta\in\mathbb R^d$, define 
\[
\ell_\omega(s,w):=\alpha_0+\alpha_1s+\beta^\top w,\qquad (s,w)\in\mathbb R\times\mathbb R^d.
\]
A direct regrouping of terms shows that for all $(u,x)\in[0,1]\times S$,
\begin{equation}\label{eq:qf}
\omega^\top \mathbf M(u,x)\,\omega
=\frac{1}{T}\sum_{t=1}^{T} K_h\!\Big(u-\tfrac{t}{T}\Big)\int_{S}
\Big[\ell_\omega\!\Big(\tfrac{\tfrac{t}{T}-u}{h},\,\tfrac{z-x}{h}\Big)\Big]^2
\prod_{j=1}^{d} K_h(x^j-z^j)\,dz.
\end{equation}

We next establish the continuity of $\mathbf M(u,x)$. Each entry is a weighted average of functions in $u$ and $x$, and both components can be treated separately. For the temporal part, note that for $r=0,1,2$ the map
\[
u \mapsto K_h\!\Big(u-\tfrac{t}{T}\Big)\Big(\tfrac{\tfrac{t}{T}-u}{h}\Big)^r
\]
is continuous in $u$ for each $t$, and finite sums preserve continuity.

Turning to the covariate part, consider the functions
\[
\begin{aligned}
f_0(w)&:=\prod_{j=1}^{d} K_h(w^j),\qquad
f_{1,k}(w):=-\tfrac{w^k}{h}\prod_{j=1}^{d} K_h(w^j),\\
f_{2,k}(w)&:=\Big(\tfrac{w^k}{h}\Big)^2\prod_{j=1}^{d} K_h(w^j),\qquad
f_{k,k'}(w):=\tfrac{w^k w^{k'}}{h^2}\prod_{j=1}^{d} K_h(w^j).
\end{aligned}
\]
Since $K$ is continuous with compact support, each $K_h$ is continuous with compact support, and hence so is every $f_\bullet$. In particular, each $f_\bullet$ is uniformly continuous on $\mathbb R^d$. For $\rho\ge 0$, define the modulus of continuity
\[
\omega_{f_\bullet}(\rho):=\sup_{\|\eta\|\le \rho}\sup_{w\in\mathbb R^d}|f_\bullet(w+\eta)-f_\bullet(w)|.
\]
Uniform continuity implies that $\omega_{f_\bullet}(\rho)\to 0$ as $\rho\to 0$. For any $x,y\in\mathbb R^d$, a change of variable $w=x-z$ yields
\[
\begin{aligned}
&\Big|\int_S \big(f_\bullet(x-z)-f_\bullet(y-z)\big)\,dz\Big|
=\Big|\int_S \big(f_\bullet(w)-f_\bullet(w+(y-x))\big)\,dz\Big|\\
& \qquad \le \int_S \omega_{f_\bullet}(\|x-y\|)\,dz
=\mathrm{Leb}(S)\,\omega_{f_\bullet}(\|x-y\|).
\end{aligned}
\]
Because $\mathrm{Leb}(S)<\infty$ by compactness, this upper bound vanishes as $y\to x$. Thus each $\int_S f_\bullet(x-z)\, dz$ is uniformly continuous in $x$, and no additional regularity of $S$ is required.

Now we know that each entry of the $(d+2) \times (d+2)$ matrix $\mathbf{M}(u,x)$ is a continuous function on the compact domain $[0,1] \times S$. The continuity of the matrix's eigenvalues-a crucial property for the subsequent argument-is a direct consequence of this fact.

This property stems from the relationship between a matrix, its characteristic polynomial, and its eigenvalues. The eigenvalues, $\lambda$, of $\mathbf{M}(u,x)$ are defined as the roots of the characteristic polynomial, given by $p(\lambda; u, x) = \det(\mathbf{M}(u,x) - \lambda \mathbf{I})$. The coefficients of this polynomial in $\lambda$ are constructed from sums and products of the entries of $\mathbf{M}(u,x)$. Since each entry is a continuous function of $(u,x)$, it follows that the coefficients of the characteristic polynomial are also continuous functions of $(u,x)$.

A well-established result from analysis states that the roots of a polynomial are continuous functions of its coefficients. Because $\mathbf{M}(u,x)$ is symmetric, its eigenvalues are real for all $(u,x)$. Consequently, each eigenvalue, and in particular the smallest eigenvalue $\lambda_{\min}(\mathbf{M}(u,x))$, is a real-valued continuous function on the compact set $[0,1] \times S$.

Fix $(u,x)\in[0,1]\times S$ and let $\omega=(\alpha_0,\alpha_1,\beta^\top)^\top\neq 0$. 
From \eqref{eq:qf}, we have
\[
\omega^\top \mathbf M(u,x)\,\omega
=\frac{1}{T}\sum_{t=1}^{T} K_h\!\Big(u-\tfrac{t}{T}\Big)\int_{S}
\Big[\alpha_0+\alpha_1\,\tfrac{\tfrac{t}{T}-u}{h}+\beta^\top\tfrac{z-x}{h}\Big]^2
\prod_{j=1}^{d} K_h(x^j-z^j)\,dz.
\]

Consider first the case $\beta\neq 0$.  
Select $t$ with $K_h(u-t/T)>0$, which is possible for every $u$ since $Th^{d+1}\to\infty$. 
The map
\[
z\ \longmapsto\ \alpha_0+\alpha_1\Big(\tfrac{t/T-u}{h}\Big)+\beta^\top\Big(\tfrac{z-x}{h}\Big)
\]
is a nonconstant affine function of $z$, and its zero set is therefore an affine hyperplane of Lebesgue measure zero. 
By condition \eqref{CS}, the set $S\cap(x+(-C_1h,C_1h)^d)$ has positive measure, and within this set 
$K_h(x^j-z^j)>0$ whenever $|x^j-z^j|<C_1h$. 
Hence $\prod_{j=1}^{d}K_h(x^j-z^j)>0$ for almost every $z$ in $S\cap(x+(-C_1h,C_1h)^d)$. 
The inner integral is therefore strictly positive, which implies $\omega^\top \mathbf M(u,x)\,\omega>0$.  

Now suppose $\beta=0$ and $(\alpha_0,\alpha_1)\neq(0,0)$.  
In this case the integrand does not depend on $z$, so that
\[
\omega^\top \mathbf M(u,x)\,\omega
=\Bigg[\frac{1}{T}\sum_{t=1}^{T} K_h\!\Big(u-\tfrac{t}{T}\Big)\Big(\alpha_0+\alpha_1\,\tfrac{t/T-u}{h}\Big)^{\!2}\Bigg]
\cdot
\int_{S}\ \prod_{j=1}^{d} K_h(x^j-z^j)\,dz.
\]
The integral factor is strictly positive by the condition of the compact set $S$.  
Moreover, since $Th^{d+1}\to\infty$, there exist distinct $t_1,t_2$ with $|u-t_i/T|<C_1h$. 
The function $s\mapsto a+bs$ with $(a,b)\neq(0,0)$ cannot vanish at two distinct points, so the time sum in brackets is strictly positive. 
Thus $\omega^\top \mathbf M(u,x)\,\omega>0$ also in this case.

We conclude that for every $(u,x)\in[0,1]\times S$ and every $\omega\neq 0$, the quadratic form satisfies $\omega^\top \mathbf M(u,x)\,\omega>0$. 
Hence $\mathbf M(u,x)$ is positive definite pointwise.

Combining the temporal and spatial parts, we conclude that every entry of $\mathbf M(u,x)$ is continuous on $[0,1]\times S$. As $\mathbf M(u,x)$ is symmetric with continuous entries, its eigenvalues vary continuously with $(u,x)$. In particular, $\lambda_{\min}(\mathbf M(u,x))$ is a continuous function of $(u,x)$.
 We now establish a uniform strictly positive lower bound for the smallest eigenvalue of the design matrix. 
Since $\lambda_{\min}(\mathbf M(u,x))$ is a continuous function of $(u,x)$ and $[0,1]\times S$ is compact, the extreme value theorem implies that $\lambda_{\min}(\mathbf M(u,x))$ attains a positive minimum on this set. Hence
\begin{equation}\label{eq:eigen}
    \inf_{u\in[0,1],\,x\in S}\lambda_{\min}\!\big(\mathbf M(u,x)\big) > 0.
\end{equation}
By Assumption (C2), the density $f(u,x)$ is bounded away from zero on $[0,1]\times S$, so that
\[
\inf_{u\in[0,1],\,x\in S}\lambda_{\min}\!\big(\mathbf M(u,x)f(u,x)\big) 
= \Big(\inf_{u,x} f(u,x)\Big)\cdot\Big(\inf_{u,x}\lambda_{\min}(\mathbf M(u,x))\Big) > 0.
\]

Define 
\begin{equation*}\label{eq:AB-def}
A(u,x):=\tfrac{1}{T}\mathbf D^\top \mathbf W\,\mathbf D,
\qquad
B(u,x):=\mathbf M(u,x)\,f(u,x),
\end{equation*}
which are symmetric for every $(u,x)\in[0,1]\times S$. 
By the Hoffman-Wielandt inequality, which bounds eigenvalue deviations by the Frobenius norm,
\begin{equation*}
\lambda_{\min}\big(A(u,x)\big)\;\ge\;\lambda_{\min}\big(B(u,x)\big)\;-\;\big\|A(u,x)-B(u,x)\big\|_{F}.
\end{equation*}
Taking the infimum of the left-hand side and combining the infimum and supremum of the terms on the right yields
\begin{equation}\label{eq:inf-sup}
\hat\lambda_{\inf}:=\inf_{u,x}\lambda_{\min}\big(A(u,x)\big)
\;\ge\;\inf_{u,x}\lambda_{\min}\big(B(u,x)\big)
-\sup_{u,x}\|A(u,x)-B(u,x)\|_{F}.
\end{equation}
The first term on the right-hand side is strictly positive by \eqref{eq:eigen}, while part~(iv) ensures that 
\[
\sup_{u,x}\|A(u,x)-B(u,x)\|_{F}=o_p(1).
\]
Consequently, for any $\delta \in (0,m)$, 
\[
\Pr\!\Big(\hat\lambda_{\inf}\;\ge\;m-\delta\Big)\;\longrightarrow\;1,
\]
where
\begin{equation}\label{eq:eps-def}
m:=\Big(\inf_{u,x} f(u,x)\Big)\cdot
\Big(\inf_{u,x}\lambda_{\min}\big(\mathbf M(u,x)\big)\Big)\;>\;0.
\end{equation}
Defining 
\[
\varepsilon := \tfrac{1}{2} m > 0,
\]
we obtain
\[
\hat\lambda_{\inf}\;\ge\;2\varepsilon - o_p(1),
\qquad\text{so that}\qquad
\Pr\!\big(\hat\lambda_{\inf}>\varepsilon\big)\;\longrightarrow\;1.
\]

This establishes the desired uniform positive definiteness of $\tfrac{1}{T}\mathbf D^\top \mathbf W \mathbf D$, completing the proof of part (v).

\end{proof}
\bigskip
\section{Appendix B}
\subsection{Proof of \hyperref[Thm4.1]{Theorem 4.1}}\label{prf4.1}
\begin{proof}

We will first, derive the convergence rate of $\sup_{u\in I_h,\,x\in S}\Big|\hat{m}^{\tl}(u,x)-m^{(0)}(u,x)\Big|$. We write 
\begin{equation}
    \begin{aligned}
        \sup_{u\in I_h,\,x\in S}|\hat m^{\tl}(u,x)-m^{(0)}(u,x)| \le \sup_{u\in I_h,\,x\in S}|\hat{b}(u,x) - b(u,x)| + \sup_{u\in I_h,\,x\in S}|\hat{m}^{(1)}(u,x) - {m}^{(1)}(u,x)|
    \end{aligned}
\end{equation}
By the assumptions, we know that,
\begin{equation}
    \begin{aligned}
        \sup_{u\in I_h,\,x\in S}\Big|\hat{m}^{(1)}(u,x)-m^{(1)}(u,x)\Big| = O_P\Big(\sqrt{\frac{\log T_1}{T_1\,h_1^{d+1}}} + \frac{1}{T_1^r\,h_1^d} + h_1^2\Big).
    \end{aligned}
\end{equation}

For the first term $\sup_{u\in I_h,\,x\in S}|\hat{b}(u,x) - b(u,x)|$, we can elaborate as
\begin{equation*}
    \begin{aligned}
        &\Big|\widehat{m^{(0)}(u, x) - m^{(1)}(u, x) } - \Big( m^{(0)}(u, x) - m^{(1)}(u, x) \Big) \Big| \\
        &= \Big|\frac{\displaystyle \sum_{t_0=1}^{T_0}K_{h_\tl}\Big(u-\tTz\Big) \prod_{j=1}^d K_{h_\tl}\Big(x^j -X_{t_0,T_0}^{j,(0)} \Big) \Big(Y_{t_0,T_0}^{(0)}-\hat{m}^{(1)}(\tTz,X_{{t_0,T_0}}^{(0)})\Big)}{\displaystyle \sum_{t=1}^{T_0}K_{h_\tl}\Big(u-\tTz\Big) \prod_{j=1}^d K_{h_\tl}\Big(x^j - X_{t_0,T_0}^{j, (0)} \Big)} -  b(u,\,  x)\Big| \\ 
        &=\Big| \frac{\displaystyle \sum_{t_0=1}^{T_0}K_{h_\tl} \Big(u-\tTz\Big) \prod_{j=1}^dK_{h_\tl}\Big(x^j - X_{t_0,T_0}^{j, (0)} \Big)\Big(Y_{t_0,T_0}^{(0)}-\hat{m}^{(1)}(\tTz,X_{{t_0,T_0}}^{(0)}) -  b(u,\, x)\Big)}{\displaystyle \sum_{t_0=1}^{T_0}K_{h_\tl}\Big(u-\tTz\Big) \prod_{j=1}^d K_{h_\tl}\Big(x^j - X_{t_0,T_0}^{j,(0)} \Big)}\Big|
    \end{aligned}
\end{equation*}
Now we put, 
\begin{align*}
    &Y_{t_0,T_0}^{(0)} = m^{(0)} \Big( \tTz, X_{t_0,T_0}^{(0)} \Big) + \varepsilon_{t_0, T_0}^{(0)} = m^{(1)} \left( \tTz,  X_{t_0,T_0}^{(0)}\right) + b\left(\tTz, X_{t_0,T_0}^{(0)}\right) + \varepsilon_{t_0, T_0}^{(0)}
\end{align*}

Then, 
\begin{align*}
    Y_{t_0,T_0}^{(0)}-\hat{m}^{(1)}(\tTz,X_{{t_0,T_0}}^{(0)}) -  b(u, x) = m^{(1)} \Big( \tTz, X_{t_0,T_0}^{(0)} \Big)- \hat{m}^{(1)}(\tTz ,X_{{t_0,T_0}}^{(0)}) + b\left(\tTz, X_{t_0,T_0}^{(0)}\right) -b(u, x)+ \varepsilon_{t_0, T_0}^{(0)} 
\end{align*}
By, triangle inequality, 
\begin{equation*}
    \begin{aligned}
        &\Big|\widehat{m^{(0)}(u, x) - m^{(1)}(u, x) } - \Big( m^{(0)}(u, x) - m^{(1)}(u, x) \Big) \Big| \\
        &=\Bigg| \frac{\displaystyle \sum_{t_0=1}^{T_0}K_{h_\tl}\Big(u-\tTz\Big) \prod_{j=1}^d K_{h_\tl}\Big(x^j-X_{t_0,T_0}^{j, (0)} \Big)\Big(Y_{t_0,T_0}^{(0)}-\hat{m}^{(1)}(\tTz, X_{{t_0,T_0}}^{(0)}) -  b(u, x)\Big)}{\displaystyle \sum_{t_0=1}^{T_0}K_{h_\tl}\Big(u-\tTz\Big) \prod_{j=1}^d K_{h_\tl}\Big(x^j - X_{t_0,T_0}^{j,(0)}\Big)}\Bigg| \\
        &\le \Bigg| \frac{\displaystyle \sum_{t_0=1}^{T_0}K_{h_\tl}\Big(u_0-\tTz\Big) \prod_{j=1}^d K_{h_\tl}\Big(x^j -X_{t_0,T_0}^{j, (0)}\Big)\Big(m^{(1)} \Big( \tTz, X_{t_0,T_0}^{(0)} \Big)- \hat{m}^{(1)}(\tTz, X_{{t_0,T_0}}^{(0)}) \Big)}{\displaystyle \sum_{t=1}^{T_0}K_{h_\tl}\Big(u-\tTz\Big) \prod_{j=1}^d K_{h_\tl}\Big(x^j -X_{t_0,T_0}^{j, (0)}\Big)}\Bigg| \\ 
        &+\Bigg| \frac{\displaystyle \sum_{t_0=1}^{T_0}K_{h_\tl}\Big(u-\tTz\Big) \prod_{j=1}^d K_{h_\tl}\Big(x^j-X_{t_0,T_0}^{j, (0)}\Big)\Big( b(\tTz, X_{t_0,T_0}^{(0)}) -b(u, x) \Big)}{\displaystyle \sum_{t_0=1}^{T_0}K_{h_\tl}\Big(u-\tTz\Big)\prod_{j=1}^d K_{h_\tl}\Big(x^j-X_{t_0,T_0}^{j, (0)}\Big)}\Bigg|\\
        &+\Bigg| \frac{\displaystyle \sum_{t_0=1}^{T_0}K_{h_\tl}\Big(u-\tTz\Big) \prod_{j=1}^d K_{h_\tl}\Big(x^j -X_{t_0,T_0}^{j,(0)}\Big)\Big(\varepsilon_{t_0, T_0}^{(0)} \Big)}{\displaystyle \sum_{t_0=1}^{T_0}K_{h_\tl}\Big(u-\tTz\Big)\prod_{j=1}^d K_{h_\tl}\Big(x^j-X_{t_0,T_0}^{j, (0)}\Big)}\Bigg|
    \end{aligned}
\end{equation*}

For the first term, since we know that $\Big|m^{(1)} \Big( \tTz, X_{t_0,T_0}^{(0)} \Big)- \hat{m}^{(1)}(\tTz, X_{{t_0,T_0}}^{(0)})\Big| = O_P\Big(\sqrt{\frac{\log T_1}{T_1\,h_1^{d+1}}} + \frac{1}{T_1^r\,h_1^d} + h_1^2\Big)$, 
\begin{equation}
    \begin{aligned}
        \Bigg| \frac{\displaystyle \sum_{t_0=1}^{T_0}K_{h_\tl}\Big(u_0-\tTz\Big) \prod_{j=1}^d K_{h_\tl}\Big(x^j -X_{t_0,T_0}^{j, (0)}\Big)\Big(m^{(1)} \Big( \tTz, X_{t_0,T_0}^{(0)} \Big)- \hat{m}^{(1)}(\tTz, X_{{t_0,T_0}}^{(0)}) \Big)}{\displaystyle \sum_{t=1}^{T_0}K_{h_\tl}\Big(u-\tTz\Big) \prod_{j=1}^d K_{h_\tl}\Big(x^j -X_{t_0,T_0}^{j, (0)}\Big)}\Bigg| \\ = O_P\left(\sqrt{\frac{\log T_0}{T_0\,h_{\tl}^{d+1}}}\Big(\sqrt{\frac{\log T_1}{T_1\,h_1^{d+1}}} + \frac{1}{T_1^r\,h_1^d} + h_1^2\Big)\right).
    \end{aligned}
\end{equation}

We now move on to the second term. Since we assumed $m^{(0)}$ and $m^{(1)}$ is twice continuously partially differentiable w.r.t. for both $u$ and $x$, $b(u, x)$ is twice continuously differentiable. We define $Z_{t_0, T_0}^{(0)} = (t_0/T_0, X_{t_0,T_0}^{(0)}), \; z = (u,x)$, 
\begin{align*}
    b\left(\tTz, X_{t_0,T_0}^{(0)}\right) -b(u, x) = \nabla_0 b(u,x) \left(\tTz -u\right) \;+\; \nabla_1 b(u,x) \left( X_{t_0,T_0}^{1, (0)} - x^1 \right)  \\ +  \cdots + \nabla_d b(u,x) \left( X_{t_0,T_0}^{d, (0)} - x^d \right) + O(\norm{Z_{t_0, T_0}^{(0)} - z}{}^2)
\end{align*}

We also defined $\sup_{(u,x)} \|\nabla b(u,x)\|_{2} = \eta_{1,b}$.
It is central to bounding the bias function $b(u,x)$. By the Mean Value Theorem, for any point $Z_{t_0,T_0}^{(0)} = (\tTz, X_{t_0,T_0}^{(0)})$ in the neighborhood of $z=(u,x)$, there exists a $\tilde{z}$ on the line segment between $z$ and $Z_{t_0,T_0}^{(0)}$ such that:
\begin{equation*}
    b(Z_{t_0,T_0}^{(0)}) - b(z) = \nabla b(z)^{\top} (Z_{t_0,T_0}^{(0)} - z).
\end{equation*}
The kernel function $K_{h_\tl}$ has compact support, meaning it is non-zero only when $\|\frac{Z_{t_0,T_0}^{(0)} - z}{h_{\tl}}\|_{2} \le C_1$ for some constant $C_1$. This implies that within the support of the kernel, $\|Z_{t_0,T_0}^{(0)} - z\| = O(h_{\tl})$. Combining these facts, we can bound the magnitude of the difference:
\begin{align*}
    |b(Z_{t_0,T_0}^{(0)}) - b(z)| &\le \|\nabla b(z)\|_{2} \cdot \|Z_{t_0,T_0}^{(0)} - z\| \\
    &= (\eta_{1,b}) \cdot O(h_{\tl}) = O(\eta_{1,b} h_{\tl}).
\end{align*}
This demonstrates how the parameter $\eta_{1,b}$ controls the local variation of the bias function.

The error in estimating the bias term can be decomposed into its stochastic part and its expectation:
\begin{equation*}
    \begin{aligned}
        \Bigg| \displaystyle \sum_{t_0=1}^{T_0}&K_{h_\tl}\Big(u-\tTz\Big) \prod_{j=1}^d K_{h_\tl}\Big(x^j -X_{t_0,T_0}^{j, (0)}\Big)\Big( b(\tTz, X_{t_0,T_0}^{(0)}) -b(u, x) \Big)\Bigg| \\
        &\le \Bigg| \displaystyle \sum_{t_0=1}^{T_0}K_{h_\tl}\Big(u-\tTz\Big)\prod_{j=1}^d K_{h_\tl}\Big(x^j -X_{t_0,T_0}^{j, (0)}\Big)\Big( b(\tTz, X_{t_0,T_0}^{(0)}) -b(u, x) \Big) \\ &\qquad - \E \Big[ \displaystyle \sum_{t_0=1}^{T_0}K_{h_\tl}\Big(u-\tTz\Big)\prod_{j=1}^d K_{h_\tl}\Big(x^j -X_{t_0,T_0}^{j, (0)}\Big)\Big( b(\tTz, X_{t_0,T_0}^{(0)}) -b(u, x) \Big) \Big]\Bigg| \\
        & \qquad +\Bigg| \E \Big[\displaystyle \sum_{t_0=1}^{T_0}K_{h_\tl}\Big(u-\tTz\Big) \prod_{j=1}^d K_{h_\tl}\Big(x^j -X_{t_0,T_0}^{j, (0)}\Big)\Big( b(\tTz, X_{t_0,T_0}^{(0)}) -b(u, x) \Big) \Big] \Bigg|
    \end{aligned}
\end{equation*}
We now analyze the first term in the inequality above. 

The uniform convergence bounds in works like \cite{Hansen_2008}, \cite{Vogt_2012} depend critically on such moment bounds of the averaged random variables. The presence of the $\eta_{1, b}$ factor in the moments of $R_{t, T}$ introduces a multiplicative factor of $\eta_{1, b}$ into the final uniform rate for the stochastic term.
This leads to the following result:
\begin{equation*}
    \begin{aligned}
        \sup_{u\in [0,1], x \in S}\Bigg| \displaystyle \sum_{t_0=1}^{T_0}&K_{h_\tl}\Big(u-\tTz\Big)\prod_{j=1}^d K_{h_\tl}\Big(x^j -X_{t_0,T_0}^{j, (0)}\Big)\Big( b(\tTz, X_{t_0,T_0}^{(0)}) -b(u, x) \Big) \\ 
        &\qquad - \E \Big[ \displaystyle \sum_{t_0=1}^{T_0}K_{h_\tl}\Big(u-\tTz\Big)\prod_{j=1}^d K_{h_\tl}\Big(x^j -X_{t_0,T_0}^{j, (0)}\Big)\Big( b(\tTz, X_{t_0,T_0}^{(0)}) -b(u, x) \Big) \Big]\Bigg| \\
        & = O_p \left(\eta_{1,b} h_{TL} \sqrt{\frac{\log T_0 }{T_0 h_{\tl}^{d+1}}} \right).
    \end{aligned}
\end{equation*}

Also, applying $\sup_{u,x}\|\nabla b(u,x)\|_{2} = \eta_{1, b}$ and similar arguments used in the proof of Theorem 4.2 of \cite{Vogt_2012}, 
\begin{equation*}
    \begin{aligned}
        \sup_{u\in I_h, x \in S}\Bigg| \E \Big[  \sum_{t_0=1}^{T_0} &K_{h_\tl}\Big(u-\tTz\Big)\prod_{j=1}^d K_{h_\tl}\Big(x^j -X_{t_0,T_0}^{j, (0)}\Big)\Big( b(\tTz, X_{t_0,T_0}^{(0)}) -b(u, x) \Big) \Big]\Bigg| \\
        & \qquad = O_p \left(\frac{ \eta_{1,b}}{T_0^{r} h_{\tl}^{d}}+\eta_{2, b} h_{\tl}^2 \right)
    \end{aligned}
\end{equation*}

Further elaborate, define $\bar{K} : \mathbb{R} \rightarrow \mathbb{R}$ be a Lipschitz continuous function with support $[-qC_1, qC_1]$ for some $q > 1$. Assume that $\bar{K}(x) = 1$ for all $x \in [-C_1, C_1]$ and write $\bar{K}_h(x) = \bar{K}( x /h)$. Then, the expectation of the bias estimator's numerator is decomposed as:
\begin{equation*}
    \begin{aligned}
        \E & \left[ \frac{1}{T_0} \sum_{t_0=1}^{T_0} K_{h_{\tl}} \left(u - \frac{t_0}{T_0}\right) \prod_{j=1}^d K_{h_{\tl}} \left(x^j - X_{t_0, T_0}^{j, (0)}\right) \left( b(\frac{t_0}{T_0},X_{t_0, T_0}^{(0)} ) -b(u,x)\right)  \right] \\
        &= Q_1(u,x) + Q_2(u,x) + Q_3(u,x) + Q_4(u,x)
    \end{aligned}
\end{equation*}
with $Q_i(u,x) = \frac{1}{T_0} \sum_{t_0=1}^{T_0} K_{h_{\tl}} ( u - \frac{t_0}{T_0} ) q_i(u,x)$. We analyze the uniform bound of each term.

and
\begin{align*}
    q_{1}(z) \;=\; \E \Bigg[ \prod_{j=1}^{d}\bar{K}_{h_{\tl}} (x^j - X_{t_0, T_0}^{j, (0)})\, \Bigg\{ \prod_{j=1}^{d} K_{h_{\tl}}\Big(x^j - X_{t_0, T_0}^{j, (0)} \Big)\;-\;  \prod_{j=1}^{d} K_{h_{\tl}} \Big(x^j - X_{t_0}^{j,(0)}\Big(\frac{t_0}{T_0}\Big)\Big) \Bigg\} \\
    \;\times\; \Big\{ b(\frac{t_0}{T_0},X_{t_0, T_0}^{(0)} ) -b(u,x) \Big\} \Bigg],
\end{align*}

\begin{align*}
    q_{2}(z) \;=\; \E &\Bigg[  \prod_{j=1}^{d} \bar{K}_{h_{\tl}} (x^j - X_{t_0, T_0}^{j, (0)})\,\prod_{j=1}^{d} K_{h_{\tl}} \Big(x^j - X_{t_0}^{j, (0)}\Big(\frac{t_0}{T_0}\Big)\Big) \\
    &\times \Bigg\{ b(\frac{t_0}{T_0},X_{t_0, T_0}^{(0)} ) -b\Big(\frac{t_0}{T_0},X_{t_0}^{(0)} \Big(\frac{t_0}{T_0}\Big) \Big) \Bigg\} \Bigg],
\end{align*}

\begin{align*}
    q_{3}(z) &\;=\; \E \Bigg[ \Bigg\{  \prod_{j=1}^{d}\bar{K}_{h_{\tl}} (x^j - X_{t_0, T_0}^{j, (0)} ) - \prod_{j=1}^{d}\bar{K}_{h_{\tl}} \Big(x^j - X_{t_0}^{j, (0)}\Big(\frac{t_0}{T_0}\Big) \Big) \Bigg\} \\
    &\times \prod_{j=1}^{d} K_{h_{\tl}} \Big( x^j - X_{t_0}^{j,(0)}\Big(\frac{t_0}{T_0}\Big) \Big) \Bigg\{ b\Big(\frac{t_0}{T_0},X_{t_0}^{(0)} \Big(\frac{t_0}{T_0}\Big) \Big) -b(u,x) \Bigg\} \Bigg],
\end{align*}

\begin{align*}
    q_{4}(z) \;=\; \E \Bigg[ \prod_{j=1}^{d}K_{h_{\tl}} \Big(x^j - X_{t_0}^{j, (0)}\Big(\frac{t_0}{T_0}\Big) \Big) \Bigg\{ b\Big(\frac{t_0}{T_0},X_{t_0}^{(0)} \Big(\frac{t_0}{T_0}\Big) \Big)-b(u,x) \Bigg\} \Bigg].
\end{align*}
\textbf{Analysis of $Q_1(u,x)$ and $Q_3(u,x)$}
These terms capture the interaction between the local stationarity approximation error in the kernel weights and the deviation of the bias function $b$. Consider $q_1(z)$:
From the local stationarity assumption (C1) and Lipschitz continuity of the kernel (C5), $\prod_{j=1}^{d} K_{h_{\tl}}\Big(x^j - X_{t_0, T_0}^{j, (0)} \Big)\;-\;  \prod_{j=1}^{d} K_{h_{\tl}} \Big(x^j - X_{t_0}^{j,(0)}\Big(\frac{t_0}{T_0}\Big)\Big)$ is bounded by $O(T_0^{-r}h_{\tl}^{-d-r})$. For $b(\frac{t_0}{T_0},X_{t_0, T_0}^{(0)} ) -b(u,x)$, using the assumption $\sup \|\nabla b\|_2 = \eta_{1, b}$ and the fact that we are inside the kernel support ($\|Z_{t_0,T_0}^{(0)} - z\|_2 = O(h_{\tl})$), we have $|b(Z_{t_0,T_0}^{(0)}) - b(u,x)| = O(\eta_{1,b}h_{\tl})$.
Combining these bounds, $|q_1(z)| = O(T_0^{-r}h_{\tl}^{-d-r}) \times O(\eta_{1, b}h_{\tl}) = O(\eta_{1,b} T_0^{-r} h_{\tl}^{-d-r+1})$.
The analysis for $q_3(z)$ is analogous. Summing over $t_0$ and noting that the sum contains $O(T_0 h_{\tl})$ effective terms due to the time kernel, we obtain:
$$ \sup_{u, x}|Q_1(u,x)| \le \frac{C \eta_{1, b}}{T_0^r h_{\tl}^{d+r-1}} \quad \text{and} \quad \sup_{u, x}|Q_3(u,x)| \le \frac{C \eta_{1,b}}{T_0^r h_{\tl}^{d+r-1}}. $$

\textbf{Analysis of $Q_2(u,x)$}
This term captures the error from approximating the process at time $t_0/T_0$.
The local stationarity assumption (C1) implies $\|Z_{t_0, T_0}^{(0)} - Z_{t_0}^{(0)}(t_0/T_0)\|_{2} \le \frac{1}{T_0} U_{t_0, T_0}$. Applying the Lipschitz property of $b$, which stems directly from the bound on its gradient:
$$ \Bigg| b(\frac{t_0}{T_0},X_{t_0, T_0}^{(0)} ) -b\Big(\frac{t_0}{T_0},X_{t_0}^{(0)} \Big(\frac{t_0}{T_0}\Big) \Big) \Bigg| \le \|\nabla b\|_2 \cdot \|Z_{t_0, T_0}^{(0)} - Z_{t_0}^{(0)}(t_0/T_0)\|_{2} \le \eta_{1,b} \frac{1}{T_0} U_{t_0, T_0} = O(\eta_{1,b} \, T_0^{-r}). $$
The expectation of the kernel products is $O(h_{\tl}^{-d})$. Therefore, after summing, we get:
$$ \sup_{u, x}|Q_2(u,x)| \le \frac{C \eta_{1, b}}{T_0^r h_{\tl}^{d}}. $$

\textbf{Analysis of $Q_4(u,x)$}
By performing a Taylor expansion of $b$ around $z=(u,x)$ and using the symmetry of the kernel $K$, the first-order terms vanish, and the leading bias is determined by the second derivatives of $b$:
$$\frac{h_{\tl}^2}{2} \int w^\top \nabla^2 b(u,x) w \prod K(w_j) dw_j. $$
The magnitude of this term is $O(h_{\tl}^2)$, and it is scaled by the second-order derivatives of $b$. Since the second derivative of $b$ is governed by $\eta_{2,b}$, the bias term's magnitude is proportional to $\eta_{2,b}$. Thus, we have:
$$ \sup_{u, x}|Q_4(u,x)| \le C\eta_{2, b} h_{\tl}^2. $$

\textbf{Conclusion on Bias}
The preceding analysis established the bound for the deterministic bias of the estimator, which arises from the expectation of the term involving $b(\tTz, X_{t_0,T_0}^{(0)}) - b(u,x)$. This part of the proof is where $\|\nabla b(u,x)\|_2 = \eta_{1, b}$, is crucial.
Combining the rates for $Q_1, Q_2, Q_3, Q_4$, the dominant terms come from $Q_2$ and $Q_4$. Summing these dominant terms yields the final bound for the expectation:
\begin{equation*}
    \begin{aligned}
        \sup_{u\in I_h, x \in S}\Bigg| \E \Big[  \sum_{t_0=1}^{T_0} &K_{h_\tl}\Big(u-\frac{t_0}{T_0}\Big)\prod_{j=1}^d K_{h_\tl}\Big(x^j -X_{t_0,T_0}^{j, (0)}\Big)\Big( b(\frac{t_0}{T_0}, X_{t_0,T_0}^{(0)}) -b(u, x) \Big) \Big]\Bigg| \\
        & = O \left( \frac{\eta_{1,b}}{T_0^{r} h_{\tl}^{d}} + \eta_{2,b} h_{\tl}^2 \right).
    \end{aligned}
\end{equation*}

Now, we recall the total error for the bias estimator $\hat{b}(u,x)$. The error $\hat{b}(u,x) - b(u,x)$ is a result of three distinct sources:
\begin{equation*}
    \begin{aligned}
        &\sup_{u\in I_h, x \in S}\Big|\hat{b}(u, x)  - b(u, x) \Big| \\
        &\le \sup_{u\in I_h, x \in S} \Bigg| \frac{\displaystyle \sum_{t_0=1}^{T_0}K_{h_\tl}\Big(u_0-\tTz\Big) \prod_{j=1}^d K_{h_\tl}\Big(x^j -X_{t_0,T_0}^{j, (0)}\Big)\Big(m^{(1)} \Big( \tTz, X_{t_0,T_0}^{(0)} \Big)- \hat{m}^{(1)}(\tTz, X_{{t_0,T_0}}^{(0)}) \Big)}{\displaystyle \sum_{t=1}^{T_0}K_{h_\tl}\Big(u-\tTz\Big) \prod_{j=1}^d K_{h_\tl}\Big(x^j -X_{t_0,T_0}^{j, (0)}\Big)}\Bigg| \\ 
        &+\sup_{u\in I_h, x \in S} \Bigg| \frac{\displaystyle \sum_{t_0=1}^{T_0}K_{h_\tl}\Big(u-\tTz\Big) \prod_{j=1}^d K_{h_\tl}\Big(x^j-X_{t_0,T_0}^{j, (0)}\Big)\Big( b(\tTz, X_{t_0,T_0}^{(0)}) -b(u, x) \Big)}{\displaystyle \sum_{t_0=1}^{T_0}K_{h_\tl}\Big(u-\tTz\Big)\prod_{j=1}^d K_{h_\tl}\Big(x^j-X_{t_0,T_0}^{j, (0)}\Big)}\Bigg|\\
        &+\sup_{u\in I_h, x \in S} \Bigg| \frac{\displaystyle \sum_{t_0=1}^{T_0}K_{h_\tl}\Big(u-\tTz\Big) \prod_{j=1}^d K_{h_\tl}\Big(x^j -X_{t_0,T_0}^{j,(0)}\Big)\Big(\varepsilon_{t_0, T_0}^{(0)} \Big)}{\displaystyle \sum_{t_0=1}^{T_0}K_{h_\tl}\Big(u-\tTz\Big)\prod_{j=1}^d K_{h_\tl}\Big(x^j-X_{t_0,T_0}^{j, (0)}\Big)}\Bigg|
    \end{aligned}
\end{equation*}

The rates for these components are as follows:
The first term is bounded by $\sup_{u,x}|\hat{m}^{(1)} - m^{(1)}| = O_P\left(\sqrt{\frac{\log T_0}{T_0 h_{\tl}^{d+1}}}\left\{\sqrt{\frac{\log T_1}{T_1 h_1^{d+1}}} + \frac{1}{T_1^r h_1^d} + h_1^2\right\}\right)$. The second term converges in $O(\eta_{1,b}h_{\tl}\sqrt{\frac{\log T_0}{T_0 h_{\tl}^{d+1}}} + \frac{\eta_{1,b}}{T_0^r h_{\tl}^d} +\eta_{2,b} h_{\tl}^2)$ rate. This is the only component where the properties of $b$, and thus $\eta_{1,b},\, \eta_{2,b}$, play a role. The third part comes from the target data's regression noise, $\varepsilon_{t_0, T_0}^{(0)}$. This term is independent of the bias function $b$. It follows the standard rate for Nadaraya-Watson estimators, which is $O_P(\sqrt{\frac{\log T_0}{T_0 h_{\tl}^{d+1}}})$. 

Combining these, the total error for the bias estimator is:
\[
\sup_{u,x}|\hat{b}(u,x) - b(u,x)| = O_P\left(\sqrt{\frac{\log T_0}{T_0 h_{\tl}^{d+1}}} + \frac{\eta_{1, b}}{T_0^r h_{\tl}^d} +\eta_{2,b} h_{\tl}^2 + \sqrt{\frac{\log T_0}{T_0 h_{\tl}^{d+1}}}\left\{\sqrt{\frac{\log T_1}{T_1 h_1^{d+1}}} + \frac{1}{T_1^r h_1^d} + h_1^2\right\}\right).
\]

Finally, we use the decomposition of the total error for the transfer learning estimator:
\[
\sup_{u \in I_h,\,x \in S}\Big|\hat{m}^{\tl}(u,x) - m^{(0)}(u,x)\Big|
\;\le\; 
\sup_{u,x}\Big|\hat{b}(u,x)-b(u,x)\Big| 
+ \sup_{u,x}\Big|\hat{m}^{(1)}(u,x)-m^{(1)}(u,x)\Big|.
\]
This yields the final convergence rate:
\[
\sup_{u \in I_h,\,x \in S}\Big|\hat{m}^{\tl}(u,x) - m^{(0)}(u,x)\Big|
= O_P\left(\sqrt{\frac{\log T_0}{T_0 h_{\tl}^{d+1}}} + \frac{\eta_{1, b}}{T_0^r h_{\tl}^d} +\eta_{2,b} h_{\tl}^2 + \sqrt{\frac{\log T_1}{T_1 h_1^{d+1}}} + \frac{1}{T_1^r h_1^d} + h_1^2 \right).
\]
This result concludes the rate of the transfer learning estimator. 
\end{proof}
\subsection{Proof of \hyperref[Thm4.2]{Theorem 4.2}}\label{prf4.2}
\begin{proof}
By the assumptions, we know that,
\begin{equation}\label{eqthm4.2normhatm1-m1}
    \begin{aligned}
        \sup_{u\in [0,1],\,x\in S}\Norm{\hat{\frakm}^{(1)}(u,x)-\frakm^{(1)}(u,x)}{2}
        = O_P\Big(\sqrt{\frac{\log T_1}{T_1\,h_1^{d+1}}} + \frac{1}{T_1^r h_1^{d-1}} + h_1^2\Big).
    \end{aligned}
\end{equation}

Then,
\begin{equation}\label{eqthm4.2hatmtl-m0}
    \begin{aligned}
        \hat{\mathfrak{m}}^{\mathrm{TL}}(u,x) - \mathfrak{m}^{(0)}(u,x) &= \hat{\mathfrak{b}}(u,x) + \hat{\mathfrak{m}}^{(1),\mathrm{TL}}(u,x) - \mathfrak{m}^{(0)}(u,x) \\[3pt]
        &= \big[\hat{\mathfrak{b}}(u,x) - \mathfrak{b}(u,x)\big] + \big[\hat{\mathfrak{m}}^{(1),\mathrm{TL}}(u,x) - \mathfrak{m}^{(1)}(u,x)\big],
    \end{aligned}
\end{equation}
where 
\[
    \hat{\mathfrak{m}}^{(1),\mathrm{TL}}(u,x)
    = \big( \hat{m}^{(1)}(u,x),\; h_{\mathrm{TL}}\hat{\nabla}_0 m^{(1)}(u,x),\; \dots,\; h_{\mathrm{TL}}\hat{\nabla}_d m^{(1)}(u,x) \big)^{\top}
\]

Now we put, 
\begin{align*}
    &Y_{t_0,T_0}^{(0)} = m^{(0)} \Big( \tTz, X_{t_0,T_0}^{(0)} \Big) + \varepsilon_{t_0, T_0}^{(0)} = m^{(1)} \left( \tTz,  X_{t_0,T_0}^{(0)}\right) + b\left(\tTz, X_{t_0,T_0}^{(0)}\right) + \varepsilon_{t, T_0}^{(0)}
\end{align*}
which means
\begin{align}
    \mbf{Y}_{T_0}^{(0)} = \mbf{m}^{(1)} + \mbf{b} + \tilde{\varepsilon}
\end{align}
where $\mbf{b} = \left(b\left(\frac{1}{T_0},\,X_{1, T_0}^{(0)}\right), \dots , b\left(\frac{T_0}{T_0},\,X_{T_0, T_0}^{(0)}\right)\right)^{\top}, \; \tilde{\varepsilon} = \left( \varepsilon_{1, T_0}^{(0)}, \dots , \,\varepsilon_{T_0, T_0}^{(0)}\right)^{\top}$. 

Now, we can write 
\begin{equation}\label{eqthm4.2hatb-b}
    \begin{aligned}
        \hat{\frakb}-\frakb
        & = (\mbf{D}_{\tl}^{\top} \mbf{W}_{\tl} \mbf{D}_{\tl})^{-1} \mbf{D}_{\tl}^{\top} \mbf{W}_{\tl} \left(\mbf{Y}_{T_0}^{(0)} - \hat{\mbf{m}}^{(1)} \right) -  \frakb \\
        & = (\mbf{D}_{\tl}^{\top} \mbf{W}_{\tl} \mbf{D}_{\tl})^{-1} \mbf{D}_{\tl}^{\top} \mbf{W}_{\tl} \left(\mbf{Y}_{T_0}^{(0)} - \hat{\mbf{m}}^{(1)} - \mbf{D}_{\tl} \frakb\right) \\
        & = (\mbf{D}_{\tl}^{\top} \mbf{W}_{\tl} \mbf{D}_{\tl})^{-1} \mbf{D}_{\tl}^{\top} \mbf{W}_{\tl} \left(\mbf{m}^{(1)} + \mbf{b} + \tilde{\varepsilon} - \hat{\mbf{m}}^{(1)} - \mbf{D}_{\tl} \frakb\right)
    \end{aligned}
\end{equation}

Also, by triangle inequality, 
\begin{equation}\label{ineqthm4.2hatb-b}
    \begin{aligned}
        & \left\|(\mbf{D}_{\tl}^{\top} \mbf{W}_{\tl} \mbf{D}_{\tl})^{-1} \mbf{D}_{\tl}^{\top} \mbf{W}_{\tl} \left(\mbf{m}^{(1)} + \mbf{b} + \tilde{\varepsilon} - \hat{\mbf{m}}^{(1)} - \mbf{D}_{\tl} \frakb \right)\right\|_{2} \\[0.2cm]
        & \quad \le \left\|(\frac{1}{T_0}\mbf{D}_{\tl}^{\top} \mbf{W}_{\tl} \mbf{D}_{\tl})^{-1} \frac{1}{T_0} \mbf{D}_{\tl}^{\top} \mbf{W}_{\tl} \left(\hat{\mbf{m}}^{(1)} -  \mbf{m}^{(1)} \right) \right\|_{2} \\[0.2cm]
        & \qquad + \left\|(\frac{1}{T_0}\mbf{D}_{\tl}^{\top} \mbf{W}_{\tl} \mbf{D}_{\tl})^{-1}\frac{1}{T_0} \mbf{D}_{\tl}^{\top} \mbf{W}_{\tl} \left(\mbf{b} -  \mbf{D}_{\tl} \frakb\right) \right\|_{2} \\[0.2cm]
        & \qquad + \left\|(\frac{1}{T_0}\mbf{D}_{\tl}^{\top} \mbf{W}_{\tl} \mbf{D}_{\tl})^{-1} \frac{1}{T_0}\mbf{D}_{\tl}^{\top} \mbf{W}_{\tl} \left(\tilde{\varepsilon} \right) \right\|_{2}
    \end{aligned}
\end{equation}

We first look at the property of $(\frac{1}{T_0}\mbf{D}_{\tl}^{\top} \mbf{W}_{\tl} \mbf{D}_{\tl})^{-1}$. 
Define $(d+2) \times (d+2)$ matrix $\mbf{M}_{\tl} (u,x)$ as
\begin{equation*}
    \begin{aligned}
        &[ \mbf{M}_{\tl} (u,x)]_{0,0} = \mu_{0,0}^{{\tl}}(u,x) = \frac{1}{T_0} \sum_{t_0=1}^{T_0} K_{h_{\tl}}\left(u-\frac{t_0}{T_0}\right) \int_{S}  \prod_{j=1}^{d} K_{h_{\tl}}(x^j -z^j ) dz \\
        &[\mbf{M}_{\tl} (u,x)]_{1,1} = \mu_{1,1}^{{\tl}}(u,x) = \frac{1}{T_0} \sum_{t_0=1}^{T_0}\left( \frac{\frac{t_0}{T_0} -u}{h_{\tl}} \right)^2 K_{h_{\tl}}\left(u-\frac{t_0}{T_0}\right) \int_{S} \prod_{j=1}^{d} K_{h_{\tl}} (x^j -z^j )  dz \\
        &[\mbf{M}_{\tl} (u,x)]_{k+1,k+1} = \mu_{k+1,k+1}^{{\tl}} (u,x) = \frac{1}{T_0} \sum_{t_0=1}^{T_0} K_{h_{\tl}}\left(u-\frac{t_0}{T_0}\right) \int_{S} \left( \frac{z^k-x^k}{h_{\tl}}\right)^2   \prod_{j=1}^{d} K_{h_{\tl}}(x^j -z^j)dz \\
        &[ \mbf{M}_{\tl}(u,x)]_{0,1} =[\mbf{M}_{\tl} (u,x)]_{1,0} = \mu_{0,1}^{{\tl}} (u,x) =\mu_{1,0}^{{\tl}} (u,x) \\& \quad= \frac{1}{T_0} \sum_{t_0=1}^{T_0} \left( \frac{\frac{t_0}{T_0}-u}{h_{\tl}} \right) K_{h_{\tl}}\left(u-\frac{t_0}{T_0}\right) \int_{S} \prod_{j = 1}^d K_{h_{\tl}}(x^j-z^j)dz \\
        &[ \mbf{M}_{\tl}(u,x)]_{0,k+1} =[\mbf{M}_{\tl} (u,x)]_{k+1,0} = \mu_{0,k+1}^{{\tl}} (u,x) =\mu_{k+1,0}^{{\tl}} (u,x) \\& \quad= \frac{1}{T_0} \sum_{t_0=1}^{T_0} K_{h_{\tl}}\left(u-\frac{t_0}{T_0}\right)\int_{S}   \left( \frac{z^k-x^k}{h_{\tl}}\right) \prod_{j = 1}^d K_{h_{\tl}}(x^j-z^j)  dz \\
        &[ \mbf{M}_{\tl} (u,x)]_{1,k+1} =[\mbf{M}_{\tl} (u,x)]_{k+1,1} = \mu_{1,k+1}^{\tl}(u,x) =\mu_{k+1,1}^{\tl}(u,x) \\& \quad= \frac{1}{T_0} \sum_{t_0=1}^{T_0} \left( \frac{\tTz -u}{h_{\tl}} \right) K_{h_{\tl}}\left(u-\tTz \right)\int_{S} \left( \frac{z^k-x^k}{h_{\tl}}\right)  \prod_{j = 1}^d K_{h_{\tl}}(x^j-z^j)  dz\\
        &[\mbf{M}_{\tl} (u,x)]_{k+1,k'+1} =[\mbf{M}_{\tl} (u,x)]_{k'+1,k+1} = \mu_{k+1,k'+1}^{{\tl}}(u,x) =\mu_{k'+1,k+1}^{{\tl}}(u,x) \\& \quad= \frac{1}{T_0} \sum_{t_0=1}^{T_0}  K_{h_{\tl}}\left(u-\frac{t_0}{T_0}\right) \int_{S} \left( \frac{z^k-x^k}{h_{\tl}}\right)\left( \frac{z^{k'}-x^{k'}}{h_{\tl}}\right)  \prod_{j = 1}^d K_{h_{\tl}}(x^j-z^j)  dz 
    \end{aligned}
\end{equation*}
where $z= (z_1, \dots, z_d)^{\top}$, $k, k'=1, \dots, d$, $k \ne k'$.
Then, by the similar arguments used in the proof of Theorem 3.2 (v), 
\begin{equation*}
    \inf_{u\in[0,1],\,x\in S} \lambda_{\min}(\mbf{M}_{\tl}(u,x)) >0
\end{equation*}
and following the similar step, we can derive 
\begin{equation*}
    \hat{\lambda}_{\inf} := \inf_{u \in [0,1], x \in S} \lambda_{\min} \left( \frac{1}{T_0} \mbf{D}_{\tl}^{\top} \mbf{W}_{\tl} \mbf{D}_{\tl}\right) > \varepsilon.
\end{equation*}

Uniformly over $(u,x)\in[0,1]\times S$, then, 
\begin{equation}\label{ineqthm4.2}
    \begin{aligned}
    & \left\|(\frac{1}{T_0}\mbf{D}_{\tl}^{\top} \mbf{W}_{\tl} \mbf{D}_{\tl})^{-1} \frac{1}{T_0} \mbf{D}_{\tl}^{\top} \mbf{W}_{\tl} \left(\hat{\mbf{m}}^{(1)} -  \mbf{m}^{(1)} \right) \right\|_{2} \le C \, \left\|\frac{1}{T_0} \mbf{D}_{\tl}^{\top} \mbf{W}_{\tl} \left(\hat{\mbf{m}}^{(1)} -  \mbf{m}^{(1)} \right) \right\|_{2} \\[0.2cm]
    & \left\|(\frac{1}{T_0}\mbf{D}_{\tl}^{\top} \mbf{W}_{\tl} \mbf{D}_{\tl})^{-1}\frac{1}{T_0} \mbf{D}_{\tl}^{\top} \mbf{W}_{\tl} \left(\mbf{b} -  \mbf{D}_{\tl} \frakb\right) \right\|_{2} \le C \, \left\|\frac{1}{T_0} \mbf{D}_{\tl}^{\top} \mbf{W}_{\tl} \left(\mbf{b} -  \mbf{D}_{\tl} \frakb\right) \right\|_{2}\\[0.2cm]
    & \left\|(\frac{1}{T_0}\mbf{D}_{\tl}^{\top} \mbf{W}_{\tl} \mbf{D}_{\tl})^{-1} \frac{1}{T_0}\mbf{D}_{\tl}^{\top} \mbf{W}_{\tl} \left(\tilde{\varepsilon} \right) \right\|_{2} \le C \, \left\| \frac{1}{T_0}\mbf{D}_{\tl}^{\top} \mbf{W}_{\tl} \left(\tilde{\varepsilon} \right) \right\|_{2}
\end{aligned}
\end{equation}

For the RHS of the second inequality, 
\begin{equation}\label{ineqthm4.2b-dtlb}
    \begin{aligned}
        \Bigg\| &\frac{1}{T_0}\mbf{D}_{\tl}^{\top} \mbf{W}_{\tl} \left(\mbf{b} -  \mbf{D}_{\tl} \frakb\right)\Bigg\|_{2} \\
        &\le \Bigg\|\frac{1}{T_0} \mbf{D}_{\tl}^{\top} \mbf{W}_{\tl} \left(\mbf{b} -  \mbf{D}_{\tl} \frakb\right)- \E \Big[ \frac{1}{T_0}\mbf{D}_{\tl}^{\top} \mbf{W}_{\tl} \left(\mbf{b} -  \mbf{D}_{\tl} \frakb\right) \Big]\Bigg\|_{2} +\Bigg\| \E \Big[\frac{1}{T_0}\mbf{D}_{\tl}^{\top} \mbf{W}_{\tl} \left(\mbf{b} -  \mbf{D}_{\tl} \frakb\right)  \Big] \Bigg\|_{2}
    \end{aligned}
\end{equation}

Since we assumed $m^{(0)}$ and $m^{(1)}$ three times continuously differentiable and the third-order partial derivatives is bounded, the bias function $b(u,x) = m^{(0)}(u,x) - m^{(1)}(u,x)$ is also three times continuously differentiable and the third-order partial derivatives is bounded.
Since, $b(u,x)$ satisfies the condition of theorem 3.1,  we can establish the uniform convergence rate for the first component of the bias estimation error. This is a critical step where our main assumption on the smoothness of the bias function $b$ comes into play. A direct application of uniform convergence results for kernel averages yields:
\begin{equation}\label{eqthm4.2.VAR}
    \begin{aligned}
        \sup_{u\in [0,1], x \in S}\Bigg\| \frac{1}{T_0}\mathbf{D}_{\tl}^{\top} \mathbf{W}_{\tl} \left(\mathbf{b} -  \mathbf{D}_{\tl} \mathfrak{b}\right)- \E \Big[ \frac{1}{T_0}\mathbf{D}_{\tl}^{\top} \mathbf{W}_{\tl} \left(\mathbf{b} -  \mathbf{D}_{\tl} \mathfrak{b}\right) \Big]\Bigg\|_{2} \\ = O_p \left( \eta_{2,b} h_{TL}^2 \sqrt{\frac{\log T_0 }{T_0 h_{\tl}^{d+1}}} \right).
    \end{aligned}
\end{equation}

The appearance of the $\eta_{2,b}$ term in the stochastic rate above is a direct consequence of our core assumption and constitutes a key part of our contribution. We briefly elaborate on how this factor arises.

In our case, the vector of random variables is $\mathbf{b} - \mathbf{D}_{\tl} \mathfrak{b}$. The $t_0$-th element of this vector is the local Taylor approximation error of the bias function $b$:
\[
R_{t_0, T_0} := b(Z_{t_0,T_0}^{(0)}) - b(z) - \nabla b(z)^{\top}(Z_{t_0,T_0}^{(0)}-z).
\]
For the locally linear estimator, the magnitude of this remainder term is determined by the second-order derivatives of $b$. Specifically, within the compact support of the kernel, where $\|Z_{t_0,T_0}^{(0)}-z\|_2 = O(h_{\tl})$, the remainder is bounded as:
\[
|R_{t_0, T_0}| \le C \cdot \sup_{u,x}\|\nabla^2 b(u,x)\|_F \cdot \|Z_{t_0,T_0}^{(0)}-z\|^2.
\]
By our assumption that $\sup_{u,x}\|\nabla^2 b(u,x)\|_F = \eta_{2,b}$, the magnitude of the random variable being averaged is bounded by $|R_{t_0}| = O(\eta_{2,b} h_{\tl}^2)$.

Also, applying similar arguments used in the proof of Theorem 3.2 ,
\begin{equation}\label{eqthm4.2.exp}
    \begin{aligned}
        \sup_{u\in [0,1], x \in S}\Bigg\| \E \Big[\frac{1}{T_0}\mbf{D}_{\tl}^{\top} \mbf{W}_{\tl} \left(\mbf{b} -  \mbf{D}_{\tl} \frakb\right) \Big] \Bigg\|_{2}  = O \left( \frac{ \eta_{2,b}}{T_0^r h_{\tl}^{d-1}}+h_{\tl}^2  \eta_{2,b} \right)
    \end{aligned}
\end{equation}

Define $Z_{t_0. T_0}^{(0)} = (\frac{t_0}{T_0}, \, X_{t_0, T_0}^{(0)})^{\top}, \; Z_{t_0}^{(0)} \left( \frac{t_0}{T_0} \right) = (\frac{t_0}{T_0} , \, X_{t_0}^{(0)} (\frac{t_0}{T_0}) )^{\top}, \; z = (u,x)^{\top}$
\begin{equation*}
    \begin{aligned}
        &\frac {1} {T_0} \mbf{D}_{\tl}^{\top} \mbf{W}_{\tl} (\mbf{b}_{\tl}\, - \, \mbf{D}_{\tl} \frakb_{\tl}(u,x)) = \\& 
        \begin{pmatrix}
            \tfrac{1}{T_0} \displaystyle\sum_{t_0=1}^T K_{h_{\tl}} \!\left(u - \tfrac{t_0}{T_0}\right) \prod_{j=1}^d K_{h_{\tl}} \!\left(x^j - X_{t_0,T_0}^{j,(0)}\right)\! \left( b(Z_{t_0. T_0}^{(0)}) - b(z) - \nabla b(z)^{\top}\! \left( Z_{t_0. T_0}^{(0)} - z\right)  \right) \\[0.5em]
            \tfrac{1}{T_0} \displaystyle\sum_{t_0=1}^{T_0}\! \left( \tfrac{\tfrac{t_0}{T_0}-u}{{h_{\tl}}} \right) \! K_{h_{\tl}} \!\left(u - \tfrac{t_0}{T_0}\right) \prod_{j=1}^d\! K_{h_{\tl}} \!\left(x^j - X_{t_0,T_0}^{j,(0)}\right)\! \left(  b(Z_{t_0. T_0}^{(0)}) - b(z) - \nabla b(z)^{\top}\! \left( Z_{t_0. T_0}^{(0)} - z\right)  \right)\\[0.5em] 
            \frac{1}{T_0} \displaystyle\sum_{t_0=1}^{T_0}\! \left( \tfrac{X_{t_0,T_0}^{1, (0)}-x^{1}}{{h_{\tl}}} \right)\! K_{h_{\tl}} \!\left(u - \tfrac{t_0}{T_0}\right) \prod_{j=1}^d K_{h_{\tl}} \!\left(x^j - X_{t
            _0,T_0}^{j, (0)}\right) \!\left(b(Z_{t_0. T_0}^{(0)}) - b(z) - \nabla b(z) ^{\top}\! \left( Z_{t_0. T_0}^{(0)} - z\right)  \right) \\[0.5em] \vdots \\[0.5em]
            \frac{1}{T_0} \displaystyle\sum_{t_0=1}^{T_0}\! \left( \tfrac{X_{t_0,T_0}^{d, (0)}-x^{d}}{{h_{\tl}}} \right)\! K_{h_{\tl}}\! \left(u - \tfrac{t_0}{T_0}\right) \prod_{j=1}^d K_{h_{\tl}} \!\left(x^j - X_{t
            _0,T_0}^{j,(0)}\right) \!\left(  b(Z_{t_0. T_0}^{(0)}) - b(z) - \nabla b(z) ^{\top}\! \left( Z_{t_0. T_0}^{(0)} - z\right)  \right)
        \end{pmatrix} 
    \end{aligned}
\end{equation*}

Define
\[
Z^{(0)}_{t_0,T_0}=\Big(\tfrac{t_0}{T_0},\,X^{(0)}_{t_0,T_0}\Big)^{\top}, 
\qquad 
Z^{(0)}_{t_0}\!\Big(\tfrac{t_0}{T_0}\Big)=\Big(\tfrac{t_0}{T_0},\,X^{(0)}_{t_0}(\tfrac{t_0}{T_0})\Big)^{\top},
\qquad 
z=(u,x)^{\top}.
\]
Let
\[
w_{t_0}^{(0)}(z;h_{\mathrm{TL}})
:=K_{h_{\mathrm{TL}}}\!\Big(u-\tfrac{t_0}{T_0}\Big)\,
\prod_{j=1}^{d}K_{h_{\mathrm{TL}}}\!\Big(x^{j}-X^{j,(0)}_{t_0,T_0}\Big),
\]
and define the $(d+2)$-vector of local regressors
\[
\psi\!\big(Z^{(0)}_{t_0,T_0};z,h_{\mathrm{TL}}\big)
:=\Big(1,\;\tfrac{\tfrac{t_0}{T_0}-u}{h_{\mathrm{TL}}},\;\tfrac{X^{1,(0)}_{t_0,T_0}-x^{1}}{h_{\mathrm{TL}}},\;\dots,\;\tfrac{X^{d,(0)}_{t_0,T_0}-x^{d}}{h_{\mathrm{TL}}}\Big)^{\top}.
\]
Then, with $b(z)=m^{(0)}(z)-m^{(1)}(z)$ and $\nabla b(z)$ the gradient in $(u,x)$,
\begin{align*}
    \frac{1} {T_0} \, \mathbf{D}_{\mathrm{TL}}^{\top} \mathbf{W}_{\mathrm{TL}} &\big(\mathbf{b}_{\mathrm{TL}} - \mathbf{D}_{\mathrm{TL}} \mathfrak{b}_{\mathrm{TL}}(u,x) \big) \\ 
    &=\frac{1}{T_0}\sum_{t_0=1}^{T_0} w_{t_0}^{(0)}(z;h_{\mathrm{TL}})\, \psi\!\big(Z^{(0)}_{t_0,T_0};z,h_{\mathrm{TL}}\big)\, \Big\{b\!\big(Z^{(0)}_{t_0,T_0}\big)-b(z)-\nabla b(z)^{\top}\!\big(Z^{(0)}_{t_0,T_0}-z\big)\Big\}.
\end{align*}

To find the uniform convergence rate of this vector, we must bound each of its elements. The first element is a standard kernel average. The other elements, which correspond to the slope coefficients, are slope-weighted kernel averages. A slope-weighted kernel term such as $(\frac{v}{h})K_h(v)$ can be rewritten as $\frac{1}{h}\tilde{K}(\frac{v}{h})$, where $\tilde{K}(v) = vK(v)$. Since the kernel $K$ is bounded, Lipschitz, and has compact support (Assumption C5), the transformed kernel $\tilde{K}$ inherits these same essential properties.

Consequently, all elements of the vector can be analyzed as kernel averages using the same uniform convergence machinery. They can therefore be shown to have the same order of magnitude. For clarity and brevity, it suffices to derive the convergence rate for the first (intercept) element, as the analysis for all other elements follows analogously.

Our goal is to derive the rate for the expectation of this first element.
\begin{equation*}
    \begin{aligned}
        \E \left[\frac{1}{T_0} \sum_{t_0=1}^T K_{h_{\tl}} \left(u - \frac{t_0}{T_0}\right) \prod_{j=1}^d K_{h_{\tl}} \left(x^j - X_{t_0,T_0}^{j,(0)}\right) \left( b(Z_{t_0. T_0}^{(0)}) - b(z) - \nabla b(z)^{\top} \left( b(Z_{t_0. T_0}^{(0)} - z\right)  \right) \right] \\ = O \left( \frac{ \eta_{2,b}}{T_0^r h_{\tl}^{d-1}}+h_{\tl}^2  \eta_{2,b} \right)
    \end{aligned}
\end{equation*}

Suppose $\bar{K} : \bbR \rightarrow \bbR$ be a Lipschitz continuous function that we defined in the proof of Theorem 4.1. Then, 
\begin{equation*}
    \begin{aligned}
        \E \left[\frac{1}{T_0} \sum_{t_0=1}^T K_{h_{\tl}} \left(u - \frac{t_0}{T_0}\right) \prod_{j=1}^d K_{h_{\tl}} \left(x^j - X_{t_0,T_0}^{j, (0)}\right) \left( b(Z_{t_0. T_0}^{(0)}) - b(z) - \nabla b(z)^{\top} \left( b(Z_{t_0. T_0}^{(0)} - z\right)  \right) \right] \\
        = Q_1(u,x) + Q_2(u,x) + Q_3(u,x) + Q_4(u,x)
    \end{aligned}
\end{equation*}
with 
\begin{align*}
    Q_i(u,x) = \frac{1}{T_0} \sum_{t_0=1}^{T_0} K_{h_{\tl}} \left( u - \frac{t_0}{T_0} \right) q_i(u,x)
\end{align*}
and
\begin{align*}
    q_{1}(z) \;=\; \E \Bigg[ \prod_{j=1}^{d}\bar{K}_{h_{\tl}} (x^j - X_{t_0, T_0}^{j, (0)})\, \Bigg\{ \prod_{j=1}^{d} K_{h_{\tl}}\Big(x^j - X_{t_0, T_0}^{j, (0)} \Big)\;-\;  \prod_{j=1}^{d} K_{h_{\tl}} \Big(x^j - X_{t_0}^{(0)}\Big(\frac{t_0}{T_0}\Big)\Big) \Bigg\} \\
    \;\times\; \Big\{ b(Z_{t_0. T_0}^{(0)}) - b(z) - \nabla b(z)^{\top} \left( Z_{t_0. T_0}^{(0)} - z\right) \Big\} \Bigg],
\end{align*}

\begin{align*}
    q_{2}(z) \;=\; \E &\Bigg[  \prod_{j=1}^{d} \bar{K}_{h_{\tl}} (x^j - X_{t_0, T_0}^{j, (0)})\,\prod_{j=1}^{d} K_{h_{\tl}} \Big(x^j - X_{t_0}^{j, (0)}\Big(\frac{t_0}{T_0}\Big)\Big) \\
    &\times \Bigg\{ b(Z_{t_0. T_0}^{(0)}) - b(Z_{t_0}^{(0)}(t_0 / T_0 )) - \nabla b(z)^{\top} \left( Z_{t_0. T_0}^{(0)} - Z_{t_0}^{(0)}(t_0/T_0)\right) \Bigg\} \Bigg],
\end{align*}

\begin{align*}
    q_{3}(z) &\;=\; \E \Bigg[ \Bigg\{  \prod_{j=1}^{d}\bar{K}_{h_{\tl}} (x^j - X_{t_0, T_0}^{j, (0)} ) - \prod_{j=1}^{d}\bar{K}_{h_{\tl}} \Big(x^j - X_{t_0}^{j, (0)}\Big(\frac{t_0}{T_0}\Big) \Big) \Bigg\} \\
    &\times \prod_{j=1}^{d} K_{h_{\tl}} \Big( x^j - X_{t_0}^{j,(0)}\Big(\frac{t_0}{T_0}\Big) \Big) \Bigg\{  b(Z_{t_0}^{(0)}(t_0/T_0)) - b(z) - \nabla b(z)^{\top} \left(Z_{t_0}^{(0)}(t_0/T_0) - z\right) \Bigg\} \Bigg],
\end{align*}

\begin{align*}
    q_{4}(z) \;=\; \E \Bigg[ \prod_{j=1}^{d}K_{h_{\tl}} \Big(x^j - X_{t_0}^{j, (0)}\Big(\frac{t_0}{T_0}\Big) \Big) \Bigg\{  b(Z_{t_0}^{(0)}(t_0/T_0)) - b(z) - \nabla b(z)^{\top} \left(Z_{t_0}^{(0)}(t_0/T_0) - z\right) \Bigg\} \Bigg].
\end{align*}

\textbf{Analysis of $Q_1(u,x)$ and $Q_3(u,x)$}
Using the assumption $\sup_{u,x} \|\nabla^2 b(u,x)\|_{F} = \eta_{2,b}$, 
\begin{equation*}
    \Big| b(Z_{t_0,T_0}^{(0)}) \;-\; b(x) - \nabla b(z)^{\top} (Z_{t_0,T_0}^{(0)} - z) \Big| \le \frac{1}{2} \eta_{2,b}h_{\tl}^2.
\end{equation*} 
Since $X_{t_0,T_0}^{(0)}$ is locally stationary, $K$ is Lipschitz ($\because$(C5)), we can infer that 
\begin{align*}
    \sup_{u\in [0,1],x\in S}|Q_1 (u,x) | &\le \frac{C\,\eta_{2,b}}{{T_0}^r h_{\tl}^{d+r-2}}
\end{align*}
uniformly in $u$ and $x$.
The analysis for $\sup_{u\in [0,1],x\in S}|Q_3(u,x)|\le \frac{C\,  \eta_{2,b}}{{T_0}^r h_{\tl}^{d+r-2}}$ is analogous. 

\textbf{Analysis of $Q_2(u,x)$}

Now, we consider 
\begin{align*}
    Q_2&(u,x) = \frac{1}{T_0} \sum_{t_0=1}^{T_0} K_{h_{\tl}} \Big( u - \tTz \Big) \\ &\qquad \times \E \Bigg[ \prod_{j=1}^d \bar{K}_{h_{\tl}} (x^j - X_{t_0,T_0}^{j,(0)})\, \prod_{j=1}^d K_{h_{\tl}}\Big(x^j - X_{t_0}^{j,(0)}\Big(\tTz \Big)\Big)\\ &\qquad  \times \Bigg\{ b(Z_{t_0,T_0}^{(0)}) - b\Big( Z_{t_0}^{(0)} (t_0/T_0)\Big) - \nabla b(z)^{\top} \left(Z_{t_0,T_0}^{(0)} -Z_{t_0}^{(0)} (t_0/T_0)\right)  \Bigg\} \Bigg].
\end{align*}

\begin{equation*}
    \begin{aligned}
        b(Z_{t_0,T_0}^{(0)}) &- b\Big( Z_{t_0}^{(0)}(t_0/T_0)\Big) - \nabla b(z)^{\top} \left(Z_{t_0,T_0}^{(0)} -Z_{t_0}^{(0)}(t_0/T_0)\right) \\
        &= b(Z_{t_0,T_0}^{(0)}) - b\Big( Z_{t_0}^{(0)}(t_0/T_0)\Big) - \nabla b\Big( Z_{t_0}^{(0)} (t_0/T_0)\Big)^{\top} \Big(Z_{t_0,T_0}^{(0)} -Z_{t_0}^{(0)}(t_0/T_0)\Big) \\
        &\qquad \qquad+  \nabla b\Big( Z_{t_0}^{(0)}(t_0/T_0)\Big)^{\top} \Big(Z_{t_0,T_0}^{(0)} -Z_{t_0}^{(0)}(t_0/T_0)\Big) - \nabla b(z)^{\top} \left(Z_{t_0,T_0}^{(0)} -Z_{t_0}^{(0)}(t_0/T_0)\right) 
    \end{aligned}
\end{equation*}
Notice that $Z_{t_0,T_0}^{(0)} -Z_{t_0}^{(0)}(t_0/T_0) = \Big(0, X_{t_0,T_0}^{(0)} - X_{t_0} ^{(0)}(t_0/T_0)\Big)^{\top}$. Since $ X_{t_0,T_0}^{(0)}$ is locally stationary, $\Big| X_{t_0,T_0}^{j,(0)} - X_{t_0}^{j,(0)}\Big(\tTz\Big) \Big| \le \frac{1}{T_0} U_{t_0,T_0}^{(0)}(\frac{t_0}{T_0})$. As, $b$ is twice continuously differentiable, for sufficiently large constant $C_{b,2}$, 
\begin{equation*}
    \begin{aligned}
        \Bigg|b(Z_{t_0,T_0}^{(0)}) - b \Big( Z_{t_0}^{(0)}(t_0/T_0)\Big) - \nabla b(z)^{\top} \left(Z_{t_0,T_0}^{(0)} -Z_{t_0}^{(0)}(t_0/T_0)\right)\Bigg| \\\le C_{b,2} \cdot  \eta_{2,b} \left( \frac 1 {T_0^r} \right)^2 + C_{b,2} \cdot  \eta_{2,b} h_{\tl} \cdot  \left( \frac 1 {{T_0}^r} \right) \\ 
        \le 2C_{b,2}\cdot \left( \frac { \eta_{2,b}h_{\tl}} {{T_0}^r}\right)
    \end{aligned}
\end{equation*}
The inequality holds for sufficiently large constant $C_{b,2}$ because $\frac{1}{T_0^r\,h_{\tl}^{d+r}} = o(1)$, $\sup_{u,x} \|\nabla^2 b(u,x)\|_{F} = \eta_{2,b}$. 

Similarly, for the sufficiently large constant $C$ which may vary line by line, using the boundedness of $K, \;\bar{K}_{h_{\tl}}$, 
\begin{align*}
    \sup_{u,x }\Big| Q_2(u,x) \Big| &\le \frac{C}{T_0} \sum_{t_0=1}^{T_0} K_{h_{\tl}} \Big( u - \tTz \Big) \cdot \left( \frac{1}{h_{\tl}^d}\right)\left( \frac { \eta_{2,b}h_{\tl}} {T_0^r}\right) \\
    &\le \frac{C\,  \eta_{2,b}} {T_0^r h_{\tl}^{d-1}}.
\end{align*}

\textbf{Analysis of $Q_4(u,x)$}

Finally, we consider
\begin{equation*}
    \begin{aligned}
        Q_4(u,x) = &\frac{1}{T_0} \sum_{t_0=1}^{T_0} K_{h_{\tl}} \Big( u - \tTz \Big) \\ 
        &\E \Bigg[ \prod_{j=1}^d K_{h_{\tl}} \Big(x^j - X_{t_0}^{j,(0)} \Big(\tTz\Big) \Big) \Bigg\{b\left(Z_{t_0}^{(0)} \left(\frac{t_0}{T_0}\right) \right) - b(z) - \nabla b(z)^{\top} \left(Z_{t_0}^{(0)} \left( \frac {t_0}{T_0} \right) - z\right)\Bigg\} \Bigg] \\
        & = \frac{1}{T_0} \sum_{t_0=1}^{T_0} K_{h_{\tl}} \Big( u - \tTz \Big) \int_S \prod_{j=1}^d K_h (x^j-y^j) \\& \qquad \qquad  \left\{b\left(\tTz, y\right) - b(u,x) - \nabla b(u,x)^{\top} \begin{pmatrix}
            \tTz - u \\[0.1em] y^1-x^1 \\ \vdots \\ y^d - x^d
        \end{pmatrix} \right\} f_{X_{t_0}^{(0)}(\tTz)}(y) dy
    \end{aligned}
\end{equation*}
for $y = (y_1, \dots , y_d)^{\top}$. We also define 
\begin{equation*}
    \begin{aligned}
        & G(\varphi,y) =  \left\{b\left(\varphi, y\right) - b(u,x) - \nabla b(u,x)^{\top} \begin{pmatrix}
            \varphi - u \\ y^1-x^1 \\ \vdots \\ y^d - x^d
        \end{pmatrix}  \right\} f_{X_{t_0}^{(0)}(\varphi)}(y) \;\;\text{ for } \varphi \in [0,1]\\
        &H\left(\varphi, x\right) = \int_S \prod_{j=1}^{d} K_{h_{\tl}} (x^j-y^j) G \left(\varphi, y\right) dy \\
        & g(\varphi) = K_{h_{\tl}} (u - \varphi) H \left(\varphi, x \right).
    \end{aligned}
\end{equation*}

Then, 
\begin{equation*}
    \begin{aligned}
        & \left| \frac{1}{T_0} \sum_{t_0=1}^{T_0} K_{h_{\tl}} \Big( u - \tTz \Big) H \left(\tTz, x \right) \right| \\
        & = \left| \frac{1}{T_0} \sum_{t=1}^{T} K_{h_{\tl}} \Big( u - \tTz \Big)  H \left(\tTz, x \right) - \int_{[0,1]} K_{h_{\tl}} (u-\varphi) H(\varphi, x) d\varphi \right| + \left|  \int_{[0,1]} K_{h_{\tl}} (u-\varphi) H(\varphi, x) d\varphi \right|
    \end{aligned}
\end{equation*}
We first bound the first absolute term. 
\begin{equation*}
    \begin{aligned}
        \left| \frac{1}{T_0} \sum_{t_0=1}^{T_0} K_{h_{\tl}} \Big( u - \tTz \Big)  H \left(\tTz, x \right) - \int_{[0,1]} K_{h_{\tl}} (u-\varphi) H(\varphi, x) d\varphi \right| &= \left|\frac{1}{T_0} \sum_{t_0=1}^{T_0} g \left( \tTz \right) - \int_{0}^1 g(\varphi) d \varphi \right| \\
        & \le \frac{1}{T_0} \sup_{\varphi \in [0,1]} |g'(\varphi)| = O\left( \frac{ \eta_{2,b}}{T_0} \right)
    \end{aligned}
\end{equation*}
The last equality holds because 
\begin{equation*}
    \begin{aligned}
        g'(\varphi) &= \frac{d} {d\varphi} g(\varphi) = \frac{d}{d\varphi} \left[ \frac1 {h_{\tl}} K(\frac{u-\varphi}{h_{\tl}})H(\varphi, x) \right]  \\
        &= - \frac{1}{{h_{\tl}}^2} K'\left( \frac{u-\varphi}{h_{\tl}}\right) H(\varphi,x) + \frac{1}{h_{\tl}} K \left( \frac{u - \varphi}{h_{\tl}}\right) \partial_{\varphi}H(\varphi, x)
    \end{aligned}
\end{equation*}
For the first term of $g'(\varphi)$, we can substitute the formula into
\begin{equation*}
    \begin{aligned}
        - \frac{1}{h_{\tl}^2} K'\left( \frac{u-\varphi}{h_{\tl}}\right) H(\varphi,x) = -\frac{1}{h_{\tl}^2} \left( K'(v)H(u-h_{\tl}v, x) \right). 
    \end{aligned}
\end{equation*}
Since $\sup_{u,x} \|\nabla^2 b(u,x)\|_{F} =  \eta_{2,b}$, 
\begin{equation*}
    \begin{aligned}
        &H(u-h_{\tl}v, x) \\&\qquad = \int_S \prod_{j=1}^d K_{h_{\tl}} (x^j-y^j) \left\{b\left(u-h_{\tl}v, y\right) - b(u,x) - \nabla b(u,x)^{\top} \begin{pmatrix}
            -h_{\tl}v \\ y^1-x^1 \\ \vdots \\ y^d - x^d
        \end{pmatrix} \right\} f_{X_{t_0}^{(0)}(u-h_{\tl}v)}(y) dy \\
        &\qquad = O( \eta_{2,b}h_{\tl}^2).
    \end{aligned}
\end{equation*}
So, 
\begin{equation*}
    - \frac{1}{h_{\tl}^2} K'\left( \frac{u-\varphi}{h_{\tl}}\right) H(\varphi,x) = O( \eta_{2,b}).     
\end{equation*}

For the second term $\frac{1}{h_{\tl}} K \left( \frac{u - \varphi}{h_{\tl}}\right) \partial_{\varphi}H(\varphi, x)$,
\begin{equation*}
    \begin{aligned}
        \partial_{\varphi}H(\varphi, x) = \int_S \prod_{j=1}^d K_{h_{\tl}}(x^j - y^j) \, \frac{\partial G}{\partial \varphi}(\varphi, y) \, dy.
    \end{aligned}
\end{equation*}
And since $\sup_{u,x} \|\nabla^2 b(u,x)\|_{F} = \eta_{2,b}$, $\partial_0 b(\varphi, y) - \partial_0 b(u, x) = O(\eta_{2,b} h_{\tl})$. This leads to 
\begin{equation*}
    \begin{aligned}
        \frac{\partial G}{\partial \varphi}(\varphi, y) = (\partial_0 b(\varphi, y) - \partial_0 b(u, x))\cdot f_{X_{t_0}^{(0)}(v)}(y) = O( \eta_{2,b}{h_{\tl}}).
    \end{aligned}
\end{equation*}
As a result, $|g'(\varphi) | \le C\cdot  \eta_{2,b} $. 

For the term $ \left|  \int_{[0,1]} K_{h_{\tl}} (u-\varphi) H(\varphi, x) d\varphi \right|$, using the assumption $\sup_{u,x} \|\nabla^2 b(u,x)\|_{F} = \eta_{2,b}$ in the last equality, 
\begin{equation*}
    \begin{aligned}
        \int_{[0,1]} & K_{h_{\tl}} (u-\varphi) H(\varphi, x) d\varphi    \\ 
        & = \int_{[0,1] } \int_{S } K_{h_{\tl}} (u-\varphi) \prod_{j=1}^d K_{h_{\tl}}(x^j - y^j) \\ & \qquad \qquad \qquad \qquad\left\{b\left(\varphi, y\right) - b(u,x) - \nabla b(u,x)^{\top} \begin{pmatrix}
            \varphi - u \\ y^1-x^1 \\ \vdots \\ y^d - x^d
        \end{pmatrix} \right\} f(\varphi, y) \; dy d \varphi \\
        & = \int_{[0,1] \times S} K(w_0) \prod_{j=1}^d K(w_j) \left\{ \frac{h_{\tl}^2}{2} (w^{\top} \nabla^2 b(z) w) +o(h_{\tl}^2) \right\} \\ & \qquad \qquad \qquad \qquad \left\{ f(u,x) + \nabla f(u,x)^{\top} \cdot h_{\tl}w+o(h_{\tl}) \right\} \; dw \\
        & \left(\because \text{integration by substitution}, w_0 = \frac{u - \varphi}{h_{\tl}}, \; w_j = \frac{x^j-y^j}{h_{\tl}}, \; w = (w_0, w_1, \dots, w_d)^{\top}\right)\\
        & = O( \eta_{2,b} h_{\tl}^2) 
    \end{aligned}
\end{equation*}
This completes the proof of 
\begin{align}
    \sup_{u \in [0,1], x \in S} \Abs{Q_4(u,x)}{} = O\left(\frac{ \eta_{2,b}}{T_0} +  \eta_{2,b} h_{\tl}^2\right)
\end{align}

Combining the results on $Q_1(u,x), \, Q_2(u,x), \, Q_3(u,x), \, Q_4(u,x)$ yields \eqref{eqthm4.2.exp}.

Now, by \eqref{ineqthm4.2b-dtlb}, \eqref{eqthm4.2.VAR}, \eqref{eqthm4.2.exp}, 
\begin{align}\label{eqthm4.2b-dtlb}
    \Big\| \Big(\mbf{D}_{\tl}^{\top} \mbf{W}_{\tl }\mbf{D}_{\tl}\Big)^{-1} \mbf{D}_{\tl}^{\top} \mbf{W}_{\tl}\, (\mbf{b}-\mbf{D}_{\tl}\frakb) \Big\|_{2} =O_p\!\Bigg( \eta_{2,b}h_{\tl}^2 \sqrt{\frac{\log T_0}{T_0 h_{\tl}^{\,d+1}}} + \frac{ \eta_{2,b}}{T_0^r h_{\tl}^{d-1}}+ \eta_{2,b} h_{\tl}^2 \Bigg).
\end{align}
This followed from a second-order Taylor expansion of $b$ with $\sup_{u.x} \|\nabla^2 b(u,x) \|_F = \eta_{2,b}$ mentioned in the statement. For the term driven by the source-fit error, note that the source vector $\hat{\frakm}^{(1)}$ uses $h_1$-scaling while $\frakb$ uses $h_{\tl}$-scaling; hence
\begin{equation*}
    \Big\| \Big(\mbf{D}_{\tl}^{\top}\mbf{W}_{\tl} \mbf{D}_{\tl}\Big)^{-1} \mbf{D}_{\tl}^{\top}\mbf{W}_{\tl} \, (\mbf{m}^{(1)}-\hat{\mbf{m}}^{(1)}) \Big\|_{2} \;\le\; \Big(\sqrt{\tfrac{\log T_0}{T_0 h_{\tl}^{\,d+1}}}\Big) \sup_{u,x}\Big\|\hat{\frakm}^{(1)}(u,x)-\frakm^{(1)}(u,x)\Big\|_{2}
\end{equation*}
by the uniform multivariate locally linear rate applied to the source sample \eqref{eqthm4.2normhatm1-m1}.

Lastly 
\begin{align*}
    \Big\|\Big(\mbf{D}_{\tl}^{\top} \mbf{W}_{\tl} \mbf{D}_{\tl} \Big)^{-1} \mbf{D}_{\tl}^{\top} \mbf{W}_{\tl}\,\tilde{\varepsilon}\Big\|_{2} = O_p\!\Bigg( \sqrt{\frac{\log T_0}{T_0 h_{\tl}^{\,d+1}}}\Bigg),
\end{align*}
holds. 

Combining the three pieces,
\begin{equation}\label{eqthm4.2normhatb-b}
    \sup_{u,x}\Big\|\hat{\frakb}(u,x)-\frakb(u,x)\Big\|_{2} =O_p\!\Big(\sqrt{\frac{\log T_0}{T_0 h_{\tl}^{\,d+1}}} +\frac{\eta_{2,b}}{T_0^{r}h_{\tl}^{\,d-1}} +\eta_{2,b}h_{\tl}^{2}\Big) +\Big(\sqrt{\frac{\log T_0}{T_0 h_{\tl}^{\,d+1}}}\Big)\Delta_1,
\end{equation}
where $\Delta_1:=O_p\!\Big(\sqrt{\frac{\log T_1}{T_1 h_1^{\,d+1}}}+\frac{1}{T_1^{r}h_1^{\,d-1}}+h_1^{2}\Big)$. 

Then by \eqref{eqthm4.2normhatm1-m1}, \eqref{eqthm4.2hatmtl-m0}, and \eqref{eqthm4.2normhatb-b}, 
\begin{equation}
    \begin{aligned} 
        &\sup_{u\in[0,1],\,x\in S}\Big|\hat{m}^{\tl}(u,x)-m^{(0)}(u,x)\Big| \\ &= O_p\!\left(\sqrt{\frac{\log T_0}{T_0 h_{\tl}^{\,d+1}}} +\frac{\eta_{2,b}}{T_0^{r}h_{\tl}^{\,d-1}} +\eta_{2,b}h_{\tl}^{2}+ \sqrt{\frac{\log T_1}{T_1 h_1^{\,d+1}}}+\frac{1}{T_1^{r}h_1^{\,d-1}}+h_1^{2} \right). 
    \end{aligned} 
\end{equation}
as claimed.
\end{proof}

\section{Appendix C}
\subsection{Additional Simulation Results}
\label{app:simulation_boxplots}

This appendix reports the complete boxplots corresponding to the simulation results in Section~5.2. 
Each panel displays the empirical distribution of the median grid errors across fifty replications for every value of $\gamma$ and for all four estimators. 
The boxplots further confirm the main conclusions of the analysis. 
For all bias families, the locally linear estimators exhibit smaller median errors and less variability than the Nadaraya–Watson estimators, demonstrating the improved stability of local linear smoothing under locally stationary dependence. 
Near $\gamma=0$, the transfer-learning estimators achieve the greatest accuracy and the smallest dispersion, indicating effective information transfer between the domains. 
As the absolute value of $\gamma$ increases, both bias and variability grow, and the transfer-learning methods gradually approach the performance of the target-only estimators, consistent with the theoretical phase transition described in Theorems~4.1 and~4.2.

\begin{figure}[H]
    \centering
    \includegraphics[width=0.85\textwidth]{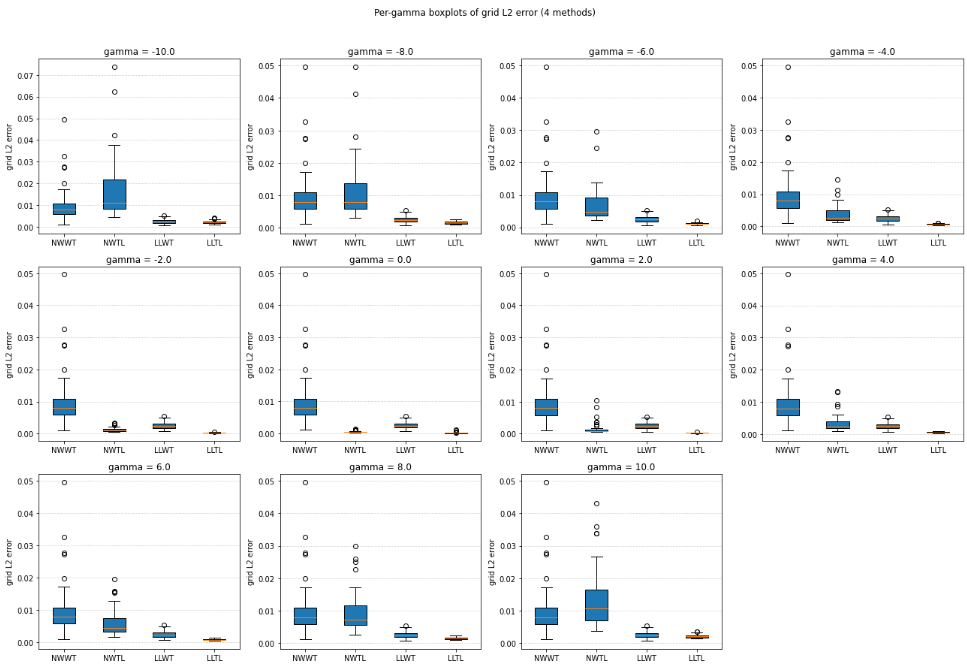}
    \caption{Boxplots of median grid errors for the quadratic bias family across fifty replications. 
    Each box corresponds to one value of $\gamma$ and compares the estimators fitted on the target sample with their transfer-learning counterparts.}
    \label{fig:quad_boxplot_appendix}
\end{figure}

\begin{figure}[H]
    \centering
    \includegraphics[width=0.85\textwidth]{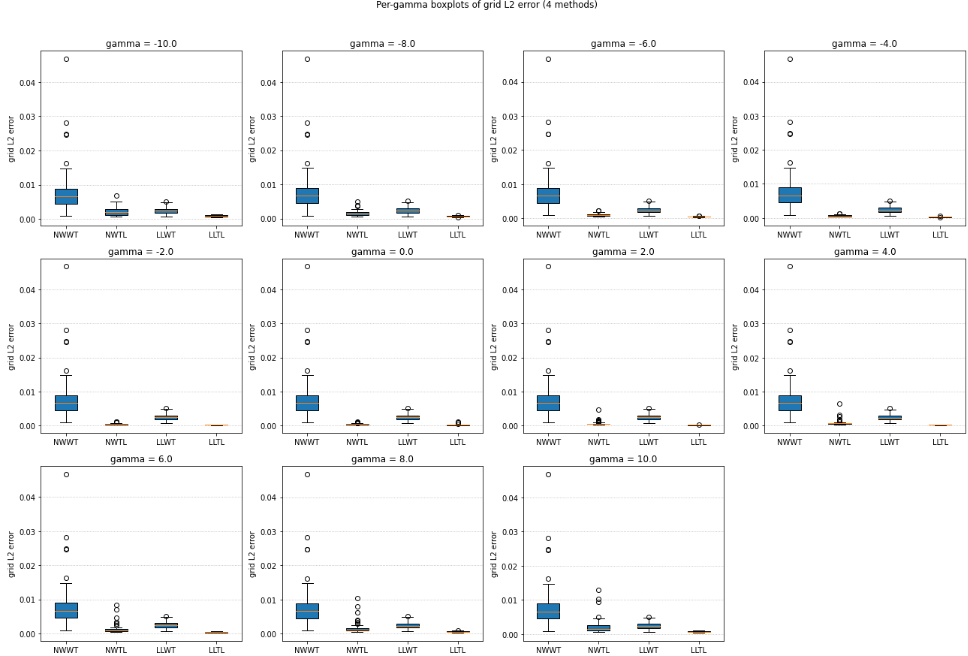}
    \caption{Boxplots of median grid errors for the cubic bias family across fifty replications. 
    The results confirm that the transfer-learning estimators provide substantial improvement near $\gamma=0$ and gradually converge to the target-only estimators as curvature increases.}
    \label{fig:cubic_boxplot_appendix}
\end{figure}

\begin{figure}[H]
    \centering
    \includegraphics[width=0.85\textwidth]{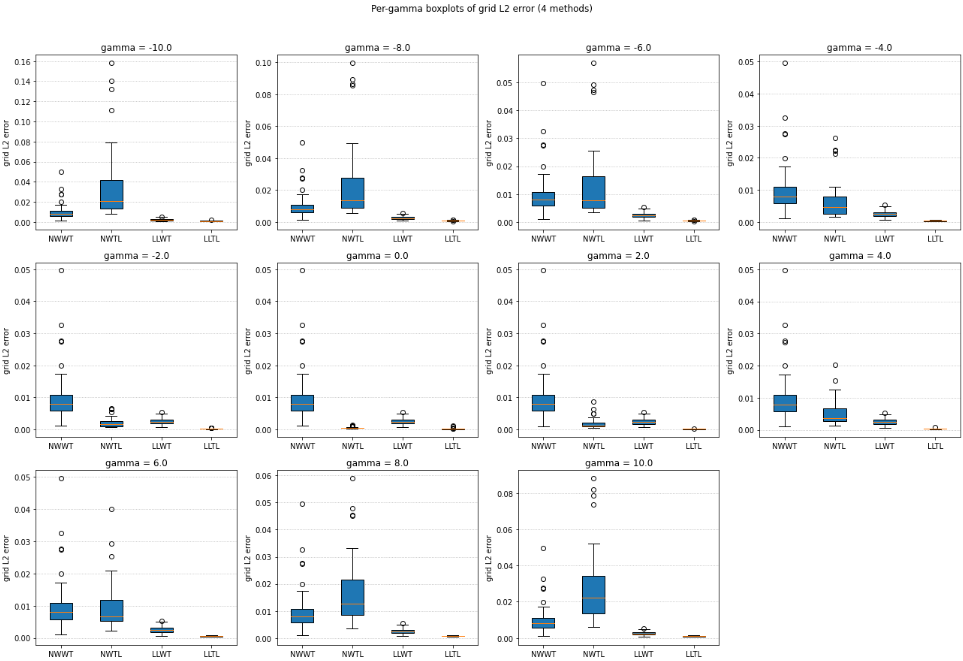}
    \caption{Boxplots of median grid errors for the exponential bias family across fifty replications. 
    The locally linear transfer-learning estimator shows the smallest error dispersion, reflecting its robustness to curvature misspecification.}
    \label{fig:exp_boxplot_appendix}
\end{figure}

\label{appendix:empirical_figures}

This appendix collects the graphical results from the empirical application described in Section 6. Each pair of panels corresponds to one experiment and estimator. The plots visualize target observations (dots), baseline fits (blue), and transfer-learning fits orange over the evaluation period.  They complement, but do not replace, the quantitative comparisons reported in the main text.
\subsection{Empirical analysis: Diesel and gasoline price relationships}

\begin{figure}[H]
\centering
\includegraphics[width=0.8\linewidth]{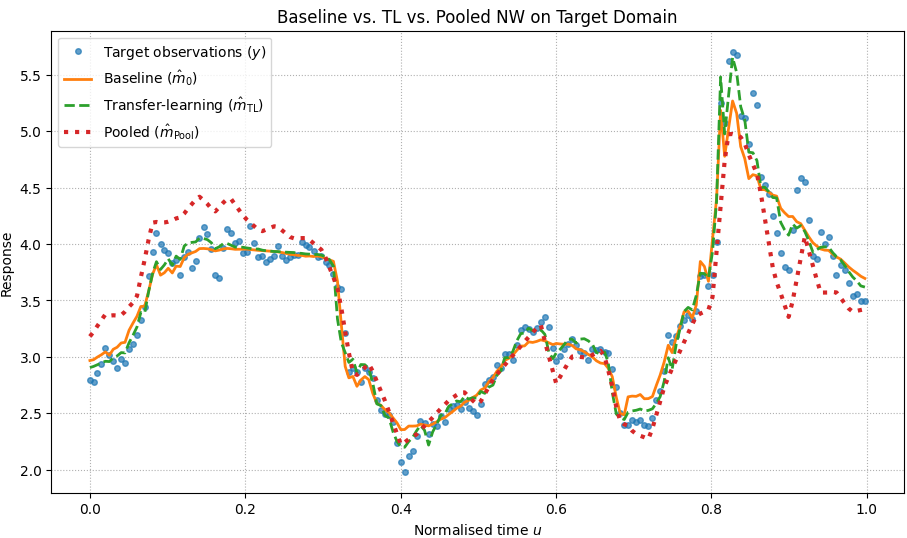}
\caption{Diesel price with WTI crude as covariate—Nadaraya–Watson estimator.}
\label{fig:diesel_WTI_nw}
\end{figure}

\begin{figure}[H]
\centering
\includegraphics[width=0.8\linewidth]{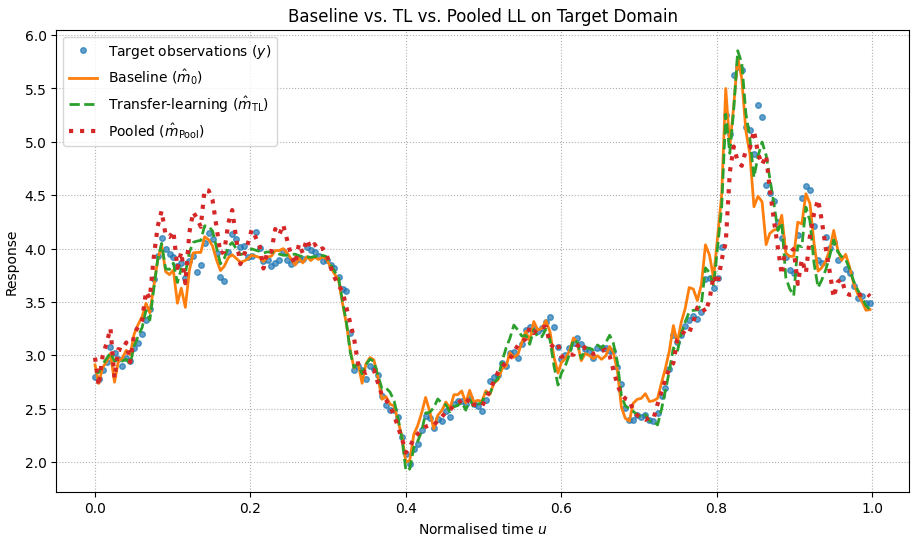}
\caption{Diesel price with WTI crude as covariate—Locally linear estimator.}
\label{fig:diesel_WTI_ll}
\end{figure}

\begin{figure}[H]
\centering
\includegraphics[width=0.8\linewidth]{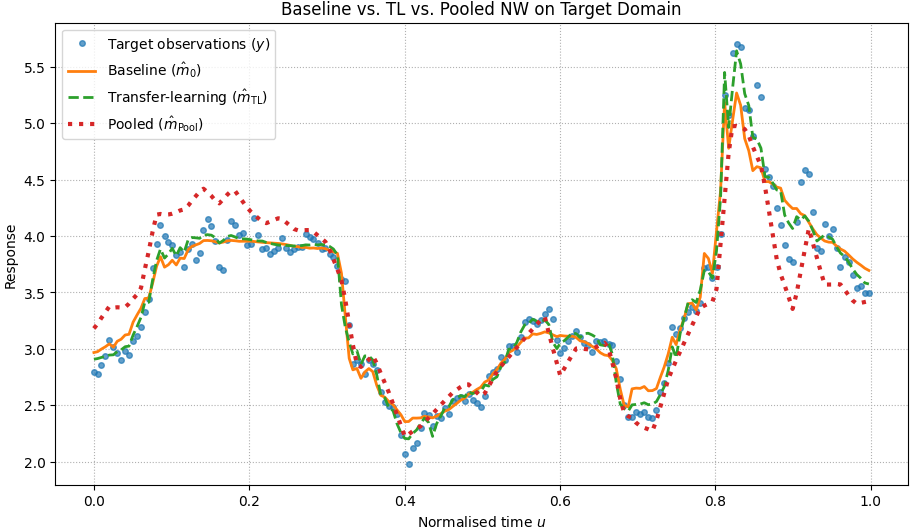}
\caption{Diesel price with Brent crude as covariate—Nadaraya–Watson estimator.}
\label{fig:diesel_brent_nw}
\end{figure}

\begin{figure}[H]
\centering
\includegraphics[width=0.8\linewidth]{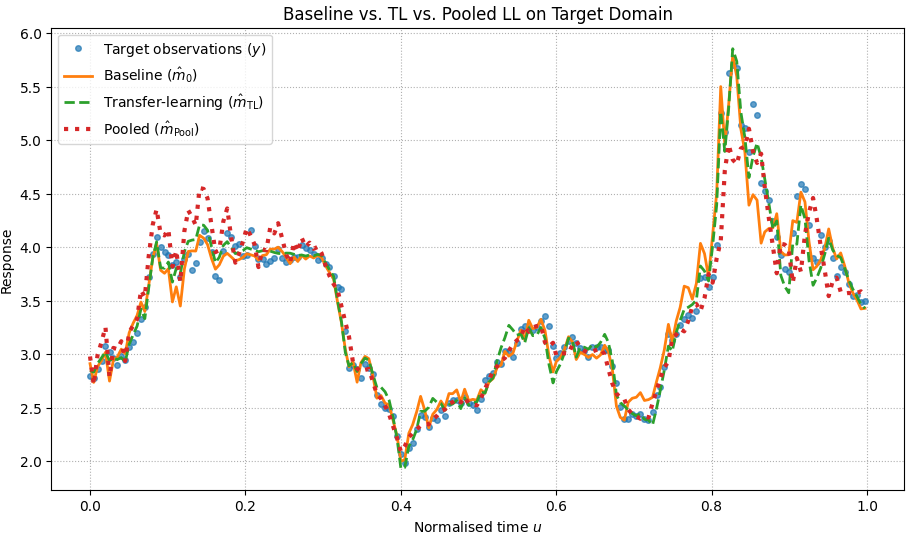}
\caption{Diesel price with Brent crude as covariate—Locally linear estimator.}
\label{fig:diesel_brent_ll}
\end{figure}

\begin{figure}[H]
\centering
\includegraphics[width=0.8\linewidth]{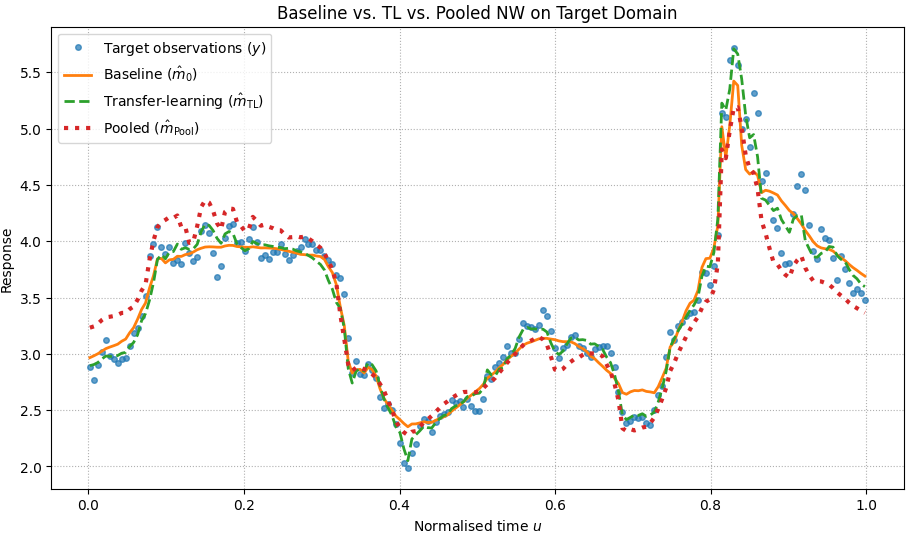}
\caption{Diesel price with gasoline as covariate—Nadaraya–Watson estimator.}
\label{fig:diesel_gasoline_nw}
\end{figure}

\begin{figure}[H]
\centering
\includegraphics[width=0.8\linewidth]{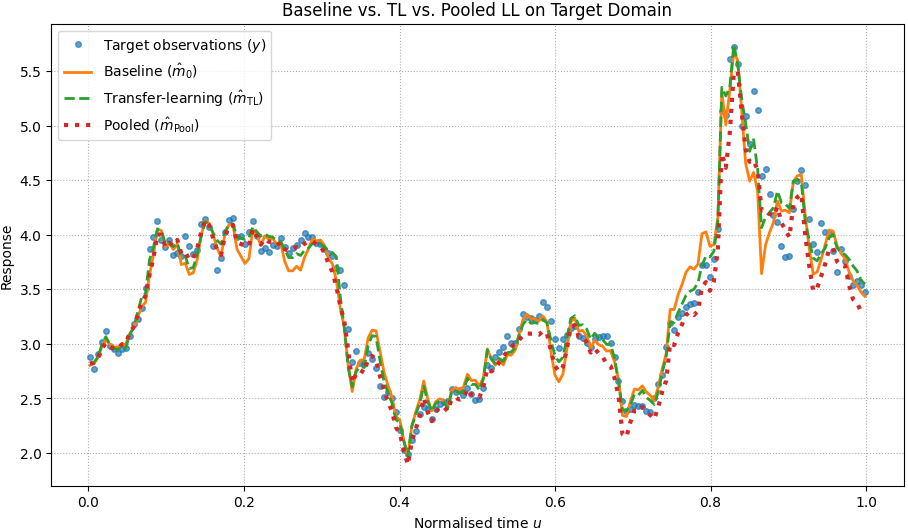}
\caption{Diesel price with gasoline as covariate—Locally linear estimator.}
\label{fig:diesel_gasoline_ll}
\end{figure}

\begin{figure}[H]
\centering
\includegraphics[width=0.8\linewidth]{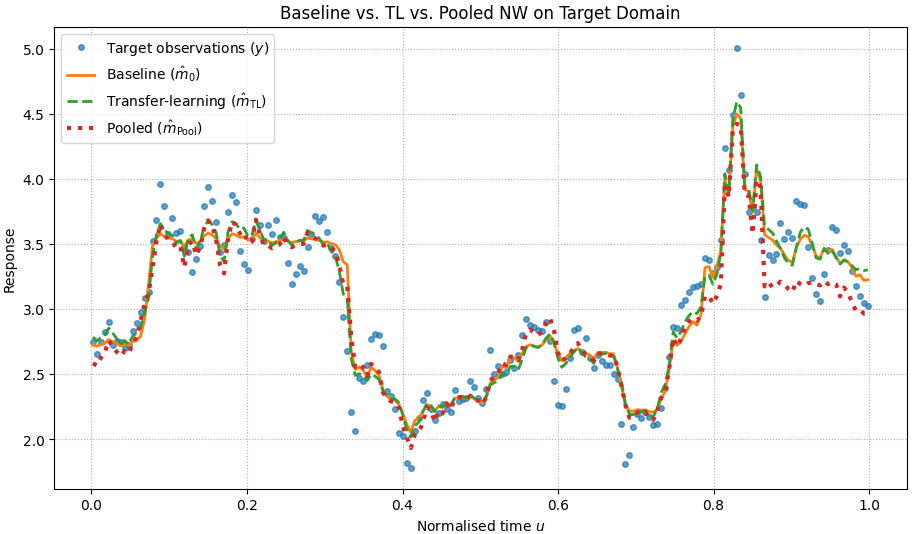}
\caption{Gasoline price with diesel as covariate—Nadaraya–Watson estimator.}
\label{fig:gasoline_diesel_nw}
\end{figure}

\begin{figure}[H]
\centering
\includegraphics[width=0.8\linewidth]{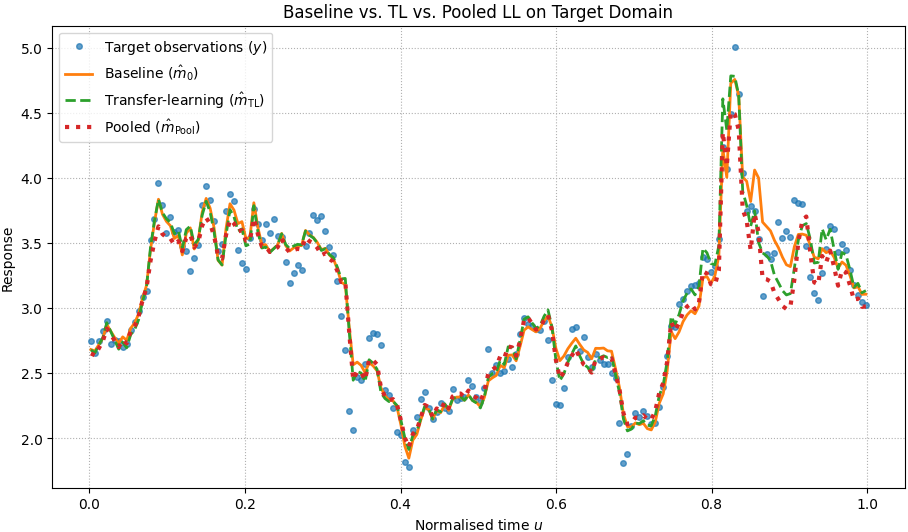}
\caption{Gasoline price with diesel as covariate—Locally linear estimator.}
\label{fig:gasoline_diesel_ll}
\end{figure}

\begin{figure}[H]
\centering
\includegraphics[width=0.8\linewidth]{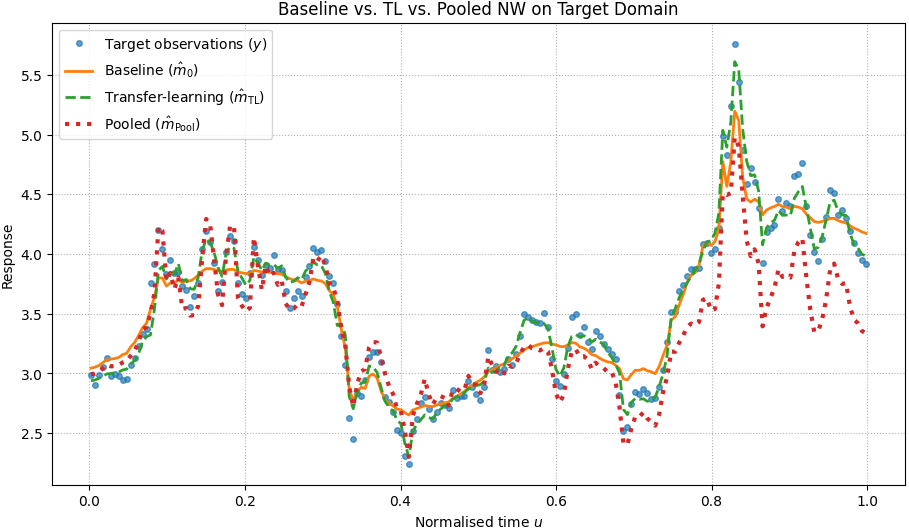}
\caption{Premium gasoline price with regular gasoline as covariate—Nadaraya–Watson estimator.}
\label{fig:premium_regular_nw}
\end{figure}

\begin{figure}[H]
\centering
\includegraphics[width=0.8\linewidth]{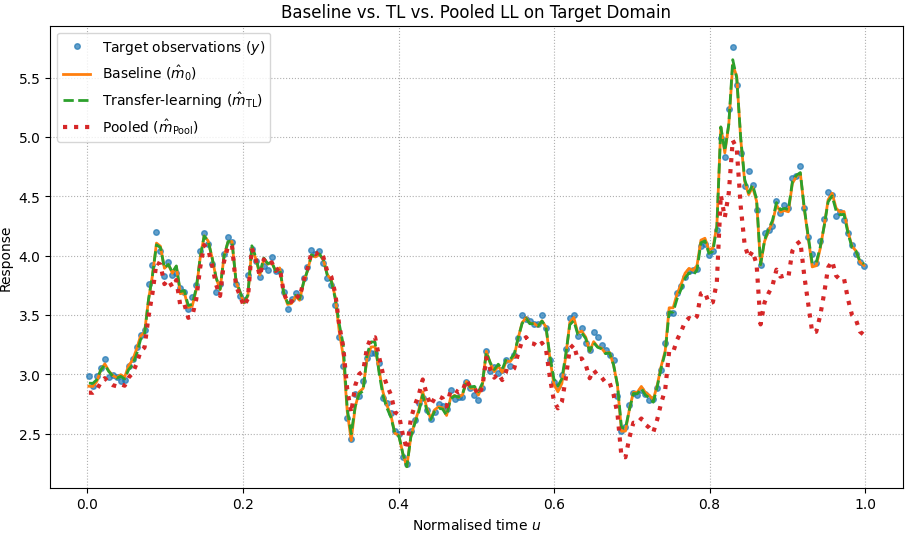}
\caption{Premium gasoline price with regular gasoline as covariate—Locally linear estimator.}
\label{fig:premium_regular_ll}
\end{figure}

To examine cross-market and cross-product dependencies in U.S. fuel prices, we estimate nonparametric regressions using both the Nadaraya–Watson and locally linear estimators under various covariate pairings. Figures~\ref{fig:diesel_WTI_nw}–\ref{fig:premium_regular_ll} illustrate how each model captures the dynamic relationship between diesel, gasoline, and crude oil prices.

\begin{figure}[H]
\centering
\includegraphics[width=0.5\linewidth]{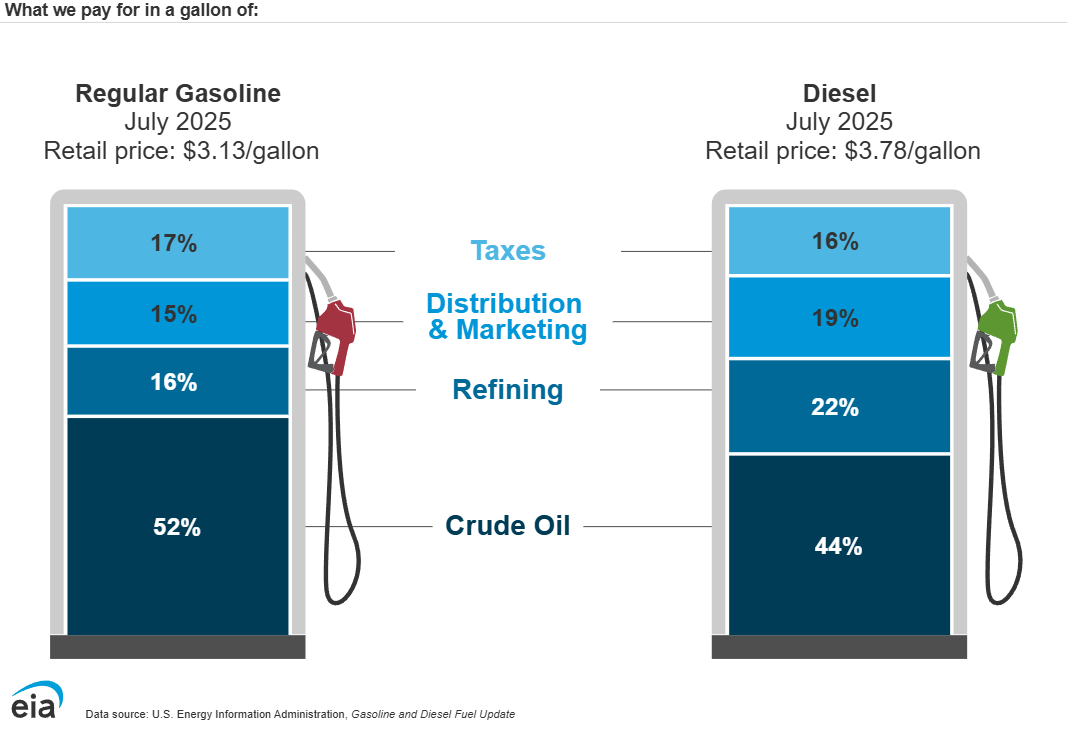}
\caption{Decomposition of U.S. retail gasoline and diesel prices (July 2025). Source: U.S. Energy Information Administration \cite{EIA_FuelBreakdown}.}
\label{fig:eia_breakdown}
\end{figure}

The nonparametric surfaces reveal pronounced nonlinear dependence between diesel and crude oil prices, particularly when Brent or WTI serve as covariates. The locally linear estimator consistently yields smoother fits and reduced boundary bias, aligning with the theoretical advantage discussed in Section~4. Similar patterns appear in the gasoline–diesel and premium–regular gasoline pairs, where the locally linear estimator provides a sharper local adjustment to changing price dynamics. Overall, the empirical results confirm that transfer learning and local linear smoothing enhance fit accuracy when cross-market linkages are strong, in agreement with the theoretical predictions.

\appendix

\end{document}